\documentclass[oneside,english]{amsart}
\usepackage[T1]{fontenc}
\usepackage[latin9]{inputenc}
\usepackage{amsthm,amssymb}
\usepackage[pdftex]{graphicx}
\makeatletter
\numberwithin{equation}{section}
\numberwithin{figure}{section}

\makeatother

\usepackage{babel}
\begin{document}

\title[\tiny{Boundary Value Problems for harmonic functions on domains in Sierpinski gaskets}]{Boundary Value Problems for harmonic functions on domains in Sierpinski gaskets}
%    author one information
\author{Shiping Cao}
\address{School of Physics, Nanjing University, Nanjing, 210093, P.R. China.}
\curraddr{} \email{shipingcao@hotmail.com}
\thanks{}

%    author two information
\author{Hua Qiu$^*$}
\address{Department of Mathematics, Nanjing University, Nanjing, 210093, P. R. China.}
\curraddr{} \email{huaqiu@nju.edu.cn}
\thanks{$^*$ Corresponding author. }
\thanks{The research of the second author was supported by the Nature Science Foundation of China, Grant 11471157.}

%    For articles to be published after 1 January 2010, you may use
%    the following version:
%\subjclass[2010]{Primary }


\date{}

\dedicatory{}
\begin{abstract}
We study boundary value problems for harmonic functions on certain domains in the level-$l$ Sierpinski gaskets $\mathcal{SG}_l$($l\geq 2$) whose boundaries are Cantor sets. We give explicit analogues of the Poisson integral formula to recover harmonic functions from their boundary values. Three types of domains, the left half domain of $\mathcal{SG}_l$ and the upper and lower domains generated by horizontal cuts of $\mathcal{SG}_l$ are considered at present.   We characterize harmonic functions of finite energy and obtain their energy estimates in terms of their boundary values. This paper settles several open problems raised in previous work. 
\end{abstract}
\maketitle

\section{Introduction }

A \textit{Dirichlet problem} is the problem of finding a function which is harmonic in the interior of a given domain that takes continuous prescribed values on the boundary of the domain. The solvability of this problem depends on the geometry of the boundary.  
For a bounded domain $D$ with sufficiently smooth boundary $\partial D$, the Dirichlet problem is always solvable, and the general solution is given by
$$u(x)=\int_{\partial D}f(s)\partial_nG(x,s)ds$$
where $G(x,y)$ is the \textit{Green's function} for $ D$, $\partial_nG(x,s)$ is the normal derivative of $G(x,y)$ along the boundary and the integration is performed on the boundary. The integral kernel $\partial_nG(x,s)$ is called the \textit{Poisson kernel} for $D$. 

With a well developed theory of Laplacians on \textit{post-critically finite (p.c.f.)} sets, originated by Kigami  {[}Ki1,Ki2{]}, it is natural to look for analogous results in the fractal context. Harmonic functions on p.c.f. self-similar sets are of finite dimension. Due to the self-similar construction of the fractal, the  Dirichlet problem on the entire fractal always reduces to solving systems of linear equations and multiplying matrices. However, for the boundary value problem on bounded subsets of fractals, the knowledge remains far from clear. 

\begin{figure}[h]
\begin{center}
\includegraphics[width=4.7cm]{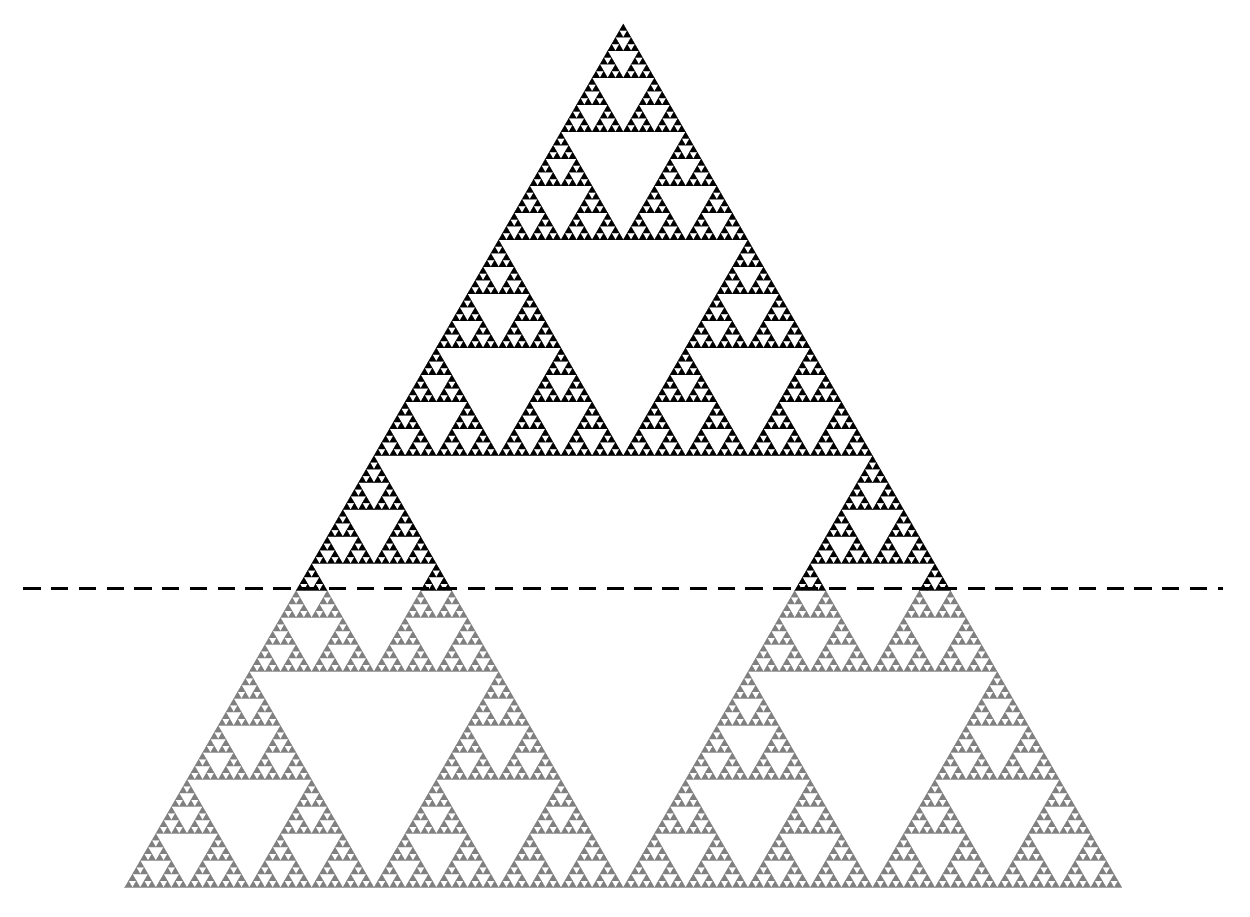}
\includegraphics[width=4.7cm]{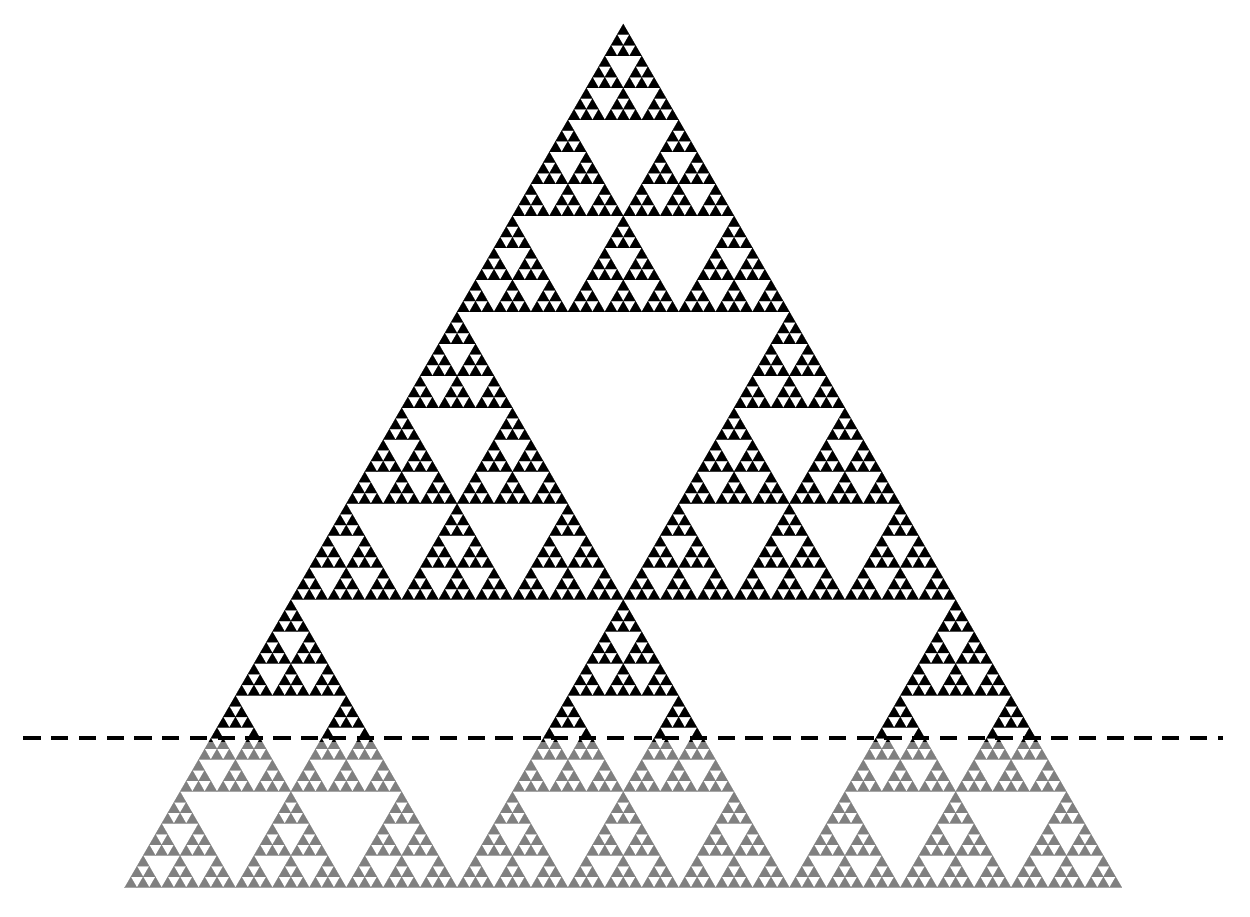}
\begin{center}
\begin{picture}(0,0) \thicklines
\put(-74.6,5.4){$(a)$}
\put(64.6,5.4){$(b)$}
\end{picture}

\textbf{Figure 1.1. Upper and lower domains in $\mathcal{SG}$ and $\mathcal{SG}_3$.}
\end{center}\end{center}
\end{figure}

The study of such problem was initiated in [S1] by Strichartz, where the upper domain generated by a horizontal cut of the Sierpinski gasket $\mathcal{SG}$ was considered. See Figure 1.1(a). Later it was continued in [OS] and [GKQS] to general case. In general, the boundary consists of a Cantor set together with the upper boundary  vertex of $\mathcal{SG}$. An explicit harmonic extension algorithm is given for solving the Dirichlet problem on such domains and the harmonic functions of finite energy are characterized in terms of their boundary values. The main tool is the \textit{Haar series expansion} of the boundary values on the Cantor set with respect to the normalized Hausdorff measure by symmetry consideration. Since the only \textit{generator} of the Haar basis is antisymmetric, one can localize the harmonic extension of this generator to any small scale along the boundary to get other basis harmonic functions. This observation plays a key role in their proof.  However, as pointed out in [GKQS], the results could not be extended to other fractals, even for the level-$3$ Sierpinski gasket $\mathcal{SG}_3$ on the base of their approach. See Figure 1.1(b). The reason is that for $\mathcal{SG}_3$ there exists a generator which is symmetric rather than antisymmetric whose harmonic extension could not be localized to small scales. On the other hand, the problem becomes much harder if we consider the domain lying below the horizontal cut instead. Except the very special case that the domains are made up of $2^m$ adjacent triangles of size $2^{-m}$ lying on the bottom line of $\mathcal{SG}$(in this case, the boundary is a finite set), we have little knowledge.

\begin{figure}[h]
\begin{center}
\includegraphics[width=2.0cm]{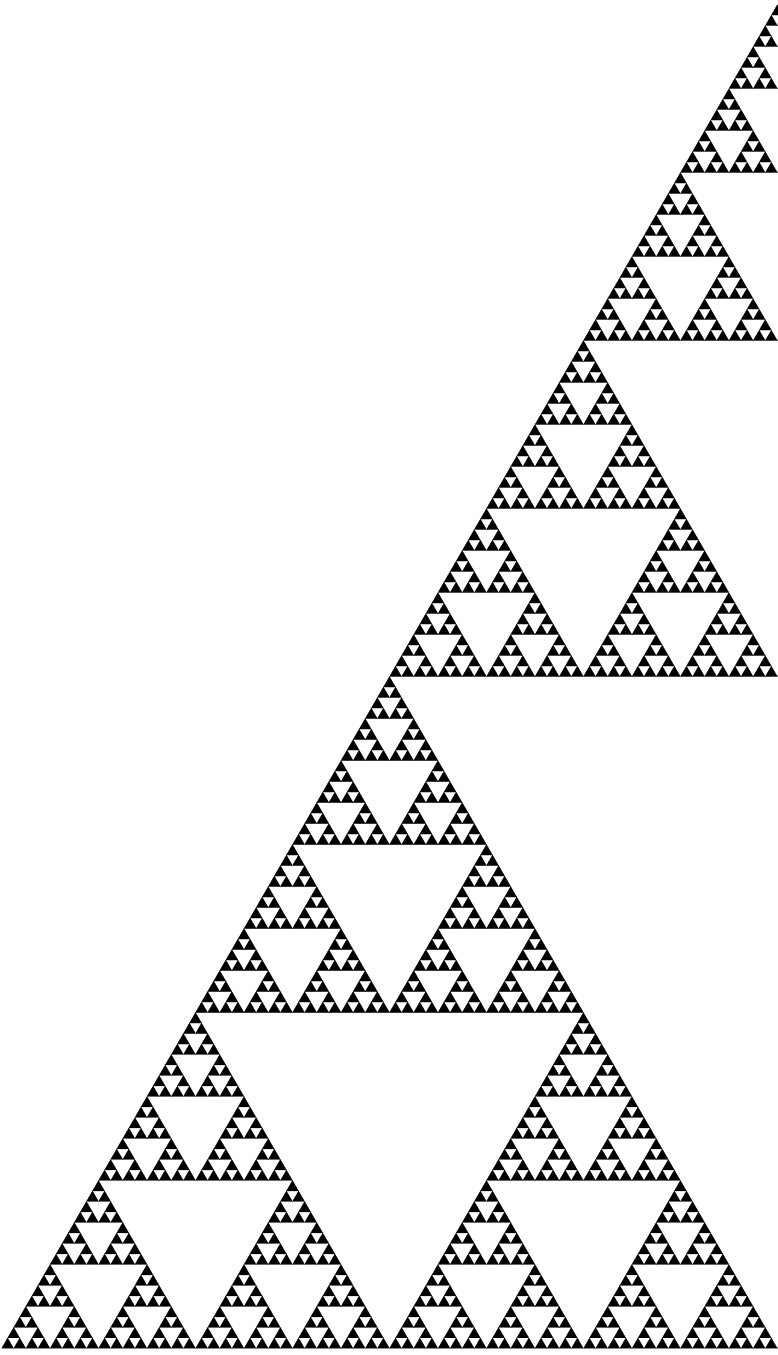}\hspace{2cm}
\includegraphics[width=2.0cm]{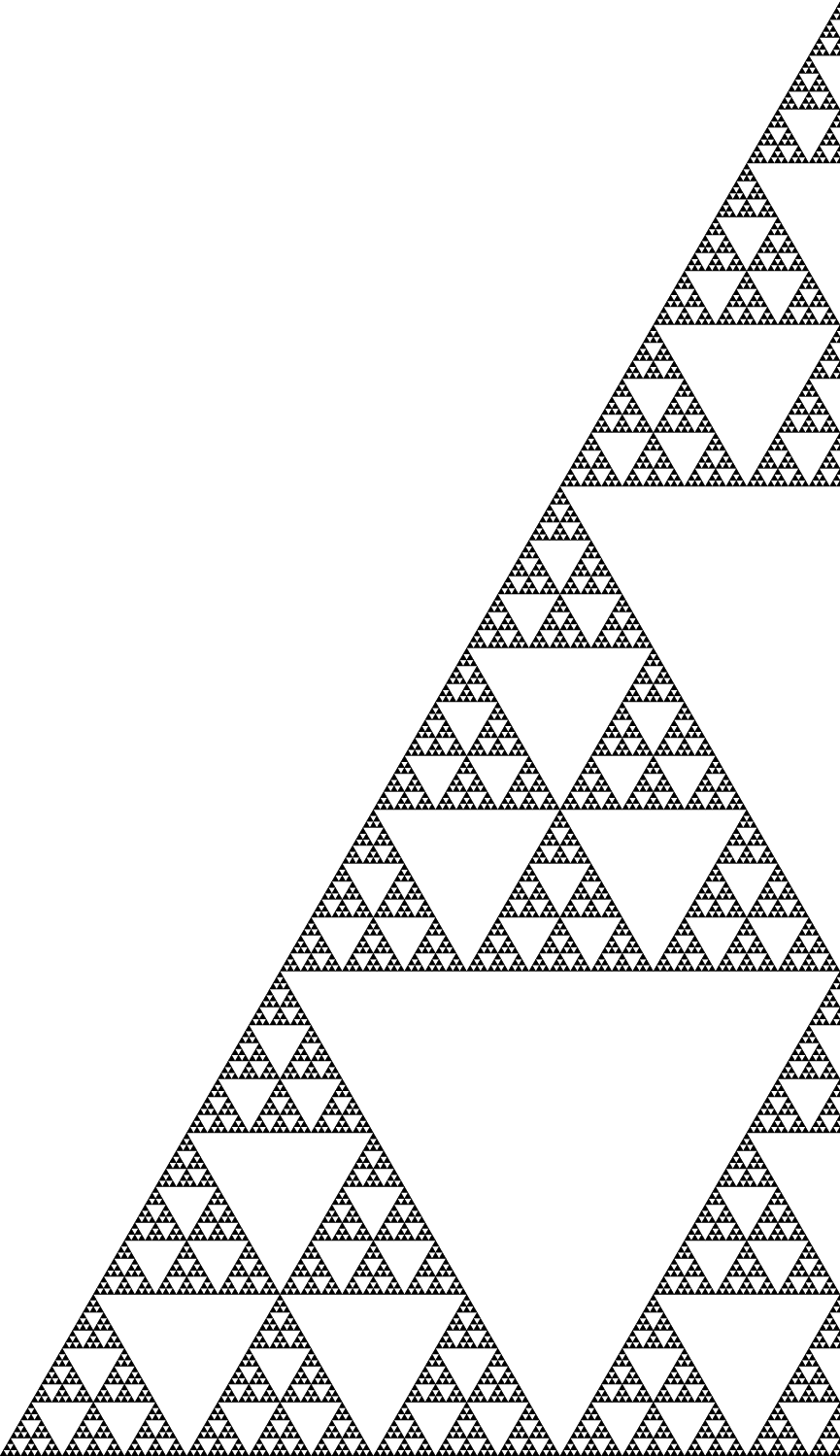}
\begin{center}
\begin{picture}(0,0) \thicklines
\put(-74.6,5.4){$(a)$}
\put(64.6,5.4){$(b)$}
\end{picture}

\textbf{Figure 1.2. Half domains in $\mathcal{SG}$ and $\mathcal{SG}_3$.}
\end{center}\end{center}
\end{figure}

 Recently, there is another natural choice of domain, namely the left half part of $\mathcal{SG}$ generated by a vertical cut along one of the symmetry lines of the gasket, becoming be interested. See Figure 1.2(a). It is the simplest example whose boundary is given as a level set of a harmonic function. In the $\mathcal{SG}$ setting, the boundary of the half domain consists of a countably infinite set of points, which makes it possible to study the Dirichlet problem by solving systems of countably infinite linear equations and multiplying infinite matrices. See [LS] for a satisfactory discussion on this domain, including an explicit harmonic extension algorithm, the characterization of harmonic functions of finite energy, and an explicit Dirichlet to Neumann map for harmonic functions. However, if consider the left half domain of level-$l$ Sierpinski gasket $\mathcal{SG}_l$ for $l\geq 3$ instead,  the approach in [LS] is not applicable. Comparing to the $\mathcal{SG}$ case, the essential difference  is that the boundary of the left half part of $\mathcal{SG}_l$ becomes a Cantor set together with the single left boundary vertex. See Figure 1.2(b).
 
 In the following, we will use \textit{upper domain}, \textit{lower domain} and \textit{half domain} to denote the above three types of domains for simplicity respectively. They are probably the simplest domains which should be handled in $\mathcal{SG}$. In this paper, we will consider the analogues of them in level-$l$ Sierpinski gasket. We will give explicit harmonic extension algorithms for all the three types of domains as well as the energy estimates for harmonic functions in terms of the boundary values (except the energy estimate for lower domains). This answers the questions raised in the above mentioned papers. As mentioned above, we need to introduce some new techniques to overcome the difficulties we met before. In fact, for each interior point $x$ in the domain $\Omega$ under consideration, we need to find the certain measure $\mu_x$ along the boundary $\partial\Omega$  analogous to $\partial_n G(x,s)ds$ in Euclidean case, so that
 $$u(x)=\int_{\partial{\Omega}}f(s)d\mu_x(s).$$ We observe that the measure $\mu_x$ is closely related to the normal derivatives of some special harmonic functions along the boundary of $\Omega$, which is crucial to our approach.   
 
 Nevertheless, these three types of domains are still the simplest domains in fractals with fractal boundary. We hope our results introduce different ideas and give insight into more general techniques for solving the Dirichlet problem and even the other boundary value problems on more general fractal domains.

\subsection{Preliminaries and the solvability of Dirichlet problems}

Let $l\geq 2$, recall that the \textit{level-$l$ Sierpinski gasket} $\mathcal{SG}_l$ is the unique nonempty compact subset of $\mathbb{R}^2$ satisfying $\mathcal{SG}_l=\bigcup_{i=0}^{\frac{l^2+l-2}{2}}F_i\mathcal{SG}_l$ with $F_i$'s being contraction mappings defined as $F_i(z)=l^{-1}z+d_{l,i}$ with suitable $d_{l,i}\in \mathbb{R}^2$. The set $V_0$ consisting of the three vertices $q_0, q_1, q_2$ of the smallest triangle containing $\mathcal{SG}_l$ is called the boundary. For convenience, we renumber $\{F_i\}_{i=0}^{\frac{l^2+l-2}{2}}$ so that $F_i(q_i)=q_i$ for $i=0,1,2$. $\mathcal{SG}_2$ is the standard \textit{Sierpinski gasket} (denoted by $\mathcal{SG}$ for simplicity). For $\mathcal{SG}_3$, in addition to $F_0, F_1, F_2$, we denote by $F_3(z)=\frac{1}{3}z+\frac{1}{3}(q_{1}+q_2)$, $F_4(z)=\frac{1}{3}z+\frac{1}{3}(q_{0}+q_2)$ and $F_5(z)=\frac{1}{3}z+\frac{1}{3}(q_{0}+q_1)$ the remaining three mappings, see Figure 1.3. These fractals have a well-developed theory of Laplacians, which allow us to perform analysis on them. In this paper, We will first describe the situation in more detail in the case of $\mathcal{SG}_3$ for half domains and upper domains, and $\mathcal{SG}$ for lower domains, then extend the considerations to general $\mathcal{SG}_l$ case. 

\begin{figure}[h]
\begin{center}
\includegraphics[width=3.7cm]{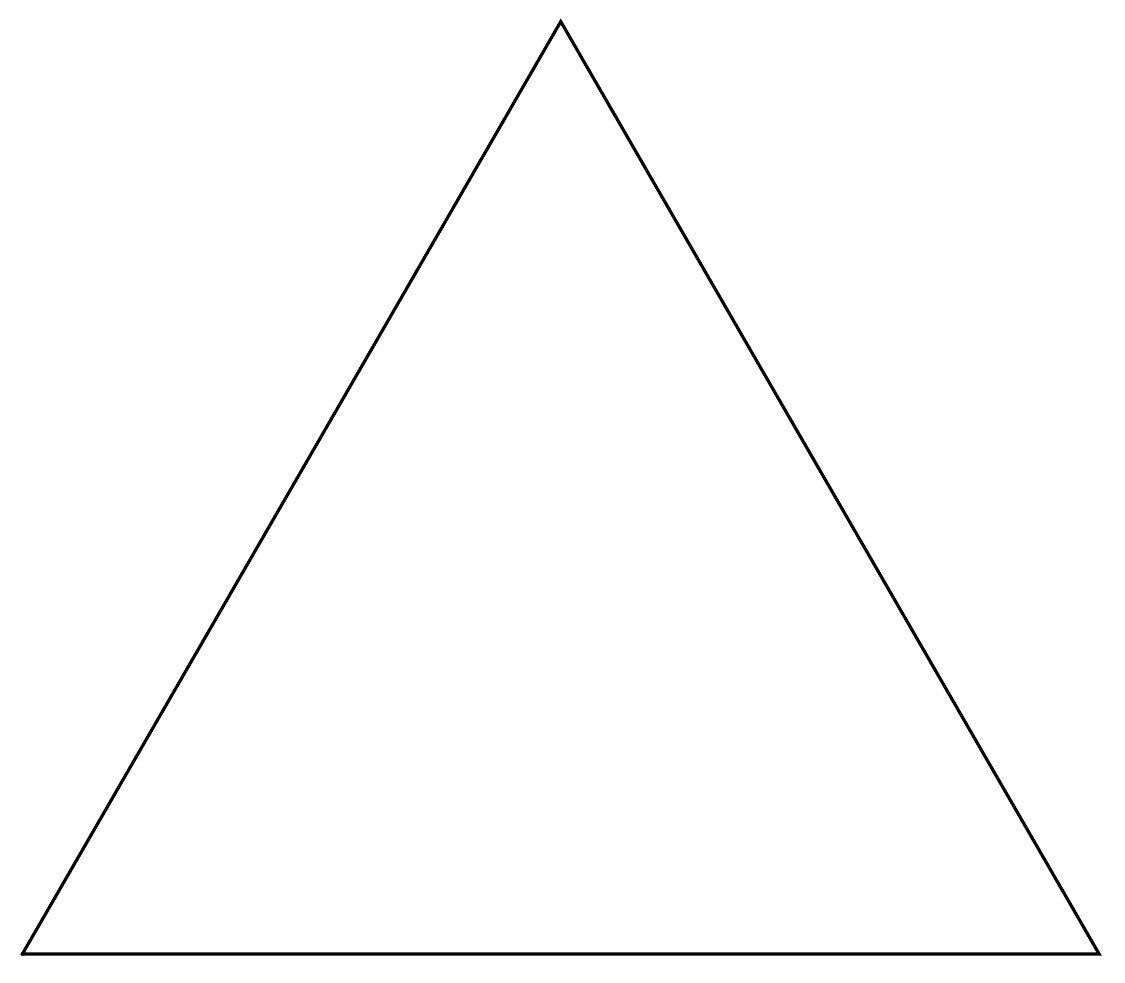}\hspace{0.3cm}
\includegraphics[width=3.7cm]{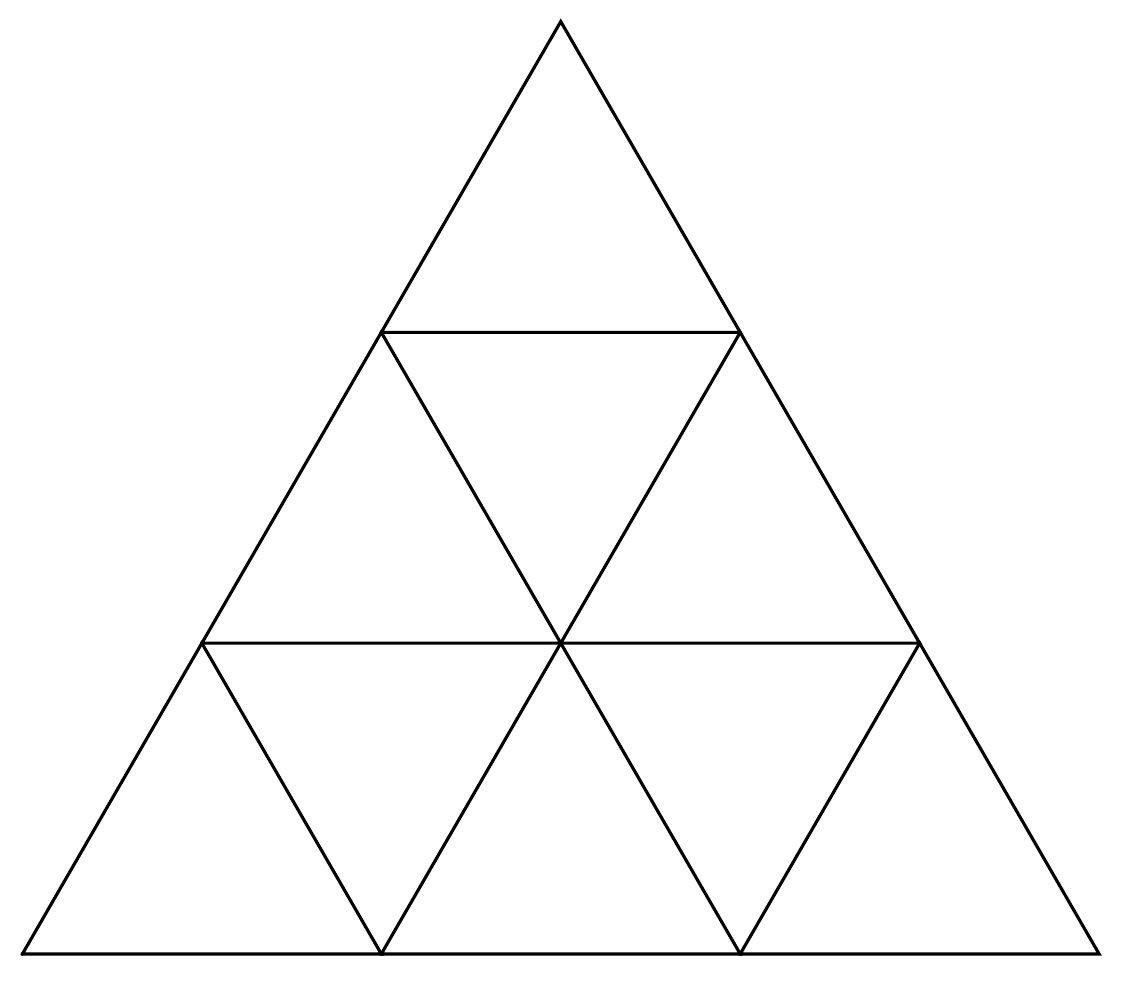}\hspace{0.3cm}
\includegraphics[width=3.7cm]{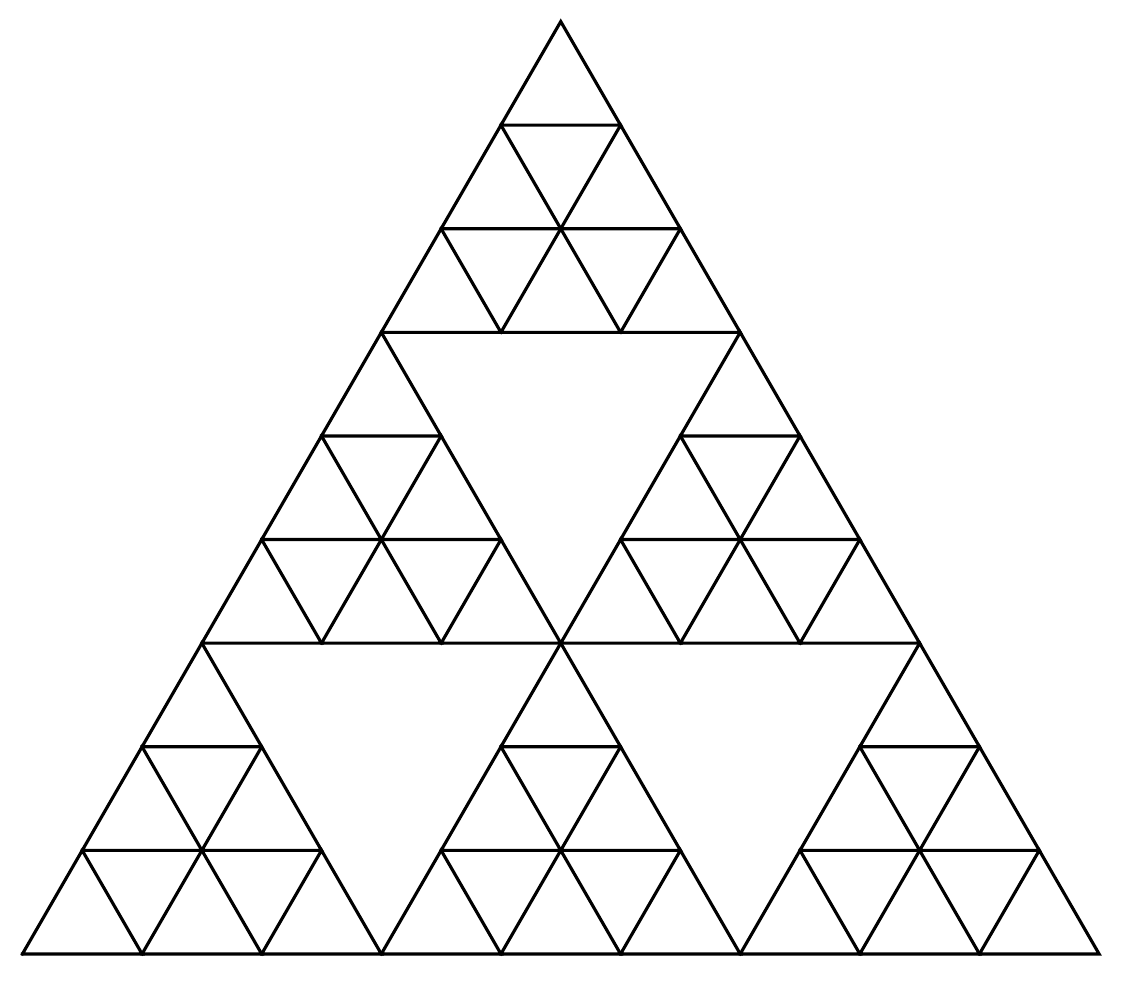}
\setlength{\unitlength}{1cm}
\begin{center}
\begin{picture}(0,0) \thicklines
\put(2.4,3.8){$q_0$}
\put(0.4,0.5){$q_1$}
\put(4.3,0.5){$q_2$}
\put(1.7,1.9){$F_5$}
\put(2.9,1.9){$F_4$}
\put(2.3,2.9){$F_0$}
\put(2.3,0.9){$F_3$}
\put(1.15,0.9){$F_1$}
\put(3.5,0.9){$F_2$}
\end{picture}
\textbf{Figure 1.3. $\Gamma_1,\Gamma_2,\Gamma_3$ of $\mathcal{SG}_3$.}
\end{center}
\end{center}
\end{figure}

We introduce some necessary notations. Readers can refer to textbooks [Ki3] and [S3] for precise definitions and known facts.  For $m\geq 1$, let $W_m=\{0,1,\cdots,\frac{l^2+l-2}{2}\}^m$ be the collection of words with length $m$ and $W_0=\{\emptyset\}$. Write $W_*=\bigcup_{m=0}^\infty W_m$, and denote the length of  $w\in W_*$ by $|w|$. For $w=w_1w_2\cdots w_m\in W_m$, we define $F_w=F_{w_1}\circ F_{w_2}\circ\cdots \circ F_{w_m}$ and call $F_w \mathcal{SG}_l$ a \textit{$m$-cell} of $\mathcal{SG}_l$.
Using the contraction mappings $F_i$, we  inductively define sets of vertices $V_m$ by $V_m=\bigcup_{i=0}^{\frac{l^2+l-2}{2}} F_iV_{m-1}$ and write $V_*=\bigcup_{m\geq 0}V_m$.  Denote $x\sim_m y$ if and only if $x\neq y$ and $x,y\in V_m$ belong to a same $m$-cell. The vertices $V_m$, together with the edge relation $\sim_m$, form a graph $\Gamma_m$ that approximates to $\mathcal{SG}_l$. See Figure 1.3 for an illustration for $\mathcal{SG}_3$.

For $m\geq 0$, the natural \textit{discrete resistance form} on $\Gamma_m$ is given by 
$$\mathcal{E}_m(u,v)=r^{-m}\sum_{x\sim_my}\big(u(x)-u(y)\big)\big(v(x)-v(y)\big)$$ for $u,v$ being functions defined on $V_m$, where $r=\frac{3}{5}$ for $\mathcal{SG}$ and $r=\frac{7}{15}$ for $\mathcal{SG}_3$.  For a real-valued function $u$ defined on ${V_*}$, it is easy to check that the \textit{graph energies} $\mathcal{E}_m(u):=\mathcal{E}_m(u,u)$ is an increasing sequence so that $\lim_{m\rightarrow\infty}\mathcal{E}_m(u)$ exists if we allow the value $+\infty$. Define 
$$\mathcal{E}(u)=\lim_{m\rightarrow\infty}\mathcal{E}_m(u)$$ to be the \textit{energy} of the function $u$  and say that $u\in dom\mathcal{E}$ if and only if $\mathcal{E}(u)<\infty$. We can regard $dom\mathcal{E}\subset C(\mathcal{SG}_l)$ since  each function of finite energy admits a unique continuous extension to $\mathcal{SG}_l$. Moreover, $dom\mathcal{E}$ is dense in $C(\mathcal{SG}_l)$. There is a natural \textit{resistance form} on $\mathcal{SG}_l$ defined as 
$$\mathcal{E}(u,v)=\lim_{m\rightarrow\infty}\mathcal{E}_m(u,v)$$ for $u,v\in dom\mathcal{E}$. 

Let $\nu$ be the standard(with equal weights) \textit{self-similar probability measure} on $\mathcal{SG}_l$. The \textit{standard Laplacian} $\Delta$ could be defined using the weak formulation. Suppose $u\in dom\mathcal{E}$ and $f$ is continuous, say $u\in dom\Delta$ with $\Delta u=f$ if
$$\mathcal{E}(u,v)=-\int_{\mathcal{SG}_l}fvd\nu$$ holds for any $v\in dom_0\mathcal{E}$ where $dom\mathcal{E}_0=\{v: v\in dom\mathcal{E}, v|_{V_0}=0\}$. 

A function $h$ is \textit{harmonic} if it minimizes the energy from each level to its next level. All the harmonic functions form a $3$-dimensional space, and hence any given values on $V_0$ can uniquely determine a harmonic function on $\mathcal{SG}_l$. They are just the solutions of the equation $\Delta h=0$. In particular, there is  an explicit extension algorithm, which determines $h|_{V_1}$ in terms of $h|_{V_0}$ and inductively $h\circ F_w|_{V_1}$ in terms of $h\circ F_w|_{V_0}$ for any $w\in W_*$ in a same manner. See Figure 1.4 for the exact formula for $\mathcal{SG}_3$.  A harmonic function $h$ satisfies the \textit{mean value property}, that is, for each $m\geq 1$, 
\[\sum_{y\sim_m x}\big(h(x)-h(y)\big)=0,\forall x\in V_m\setminus V_0.\]

\begin{figure}[h]
\begin{center}
\includegraphics[width=6cm]{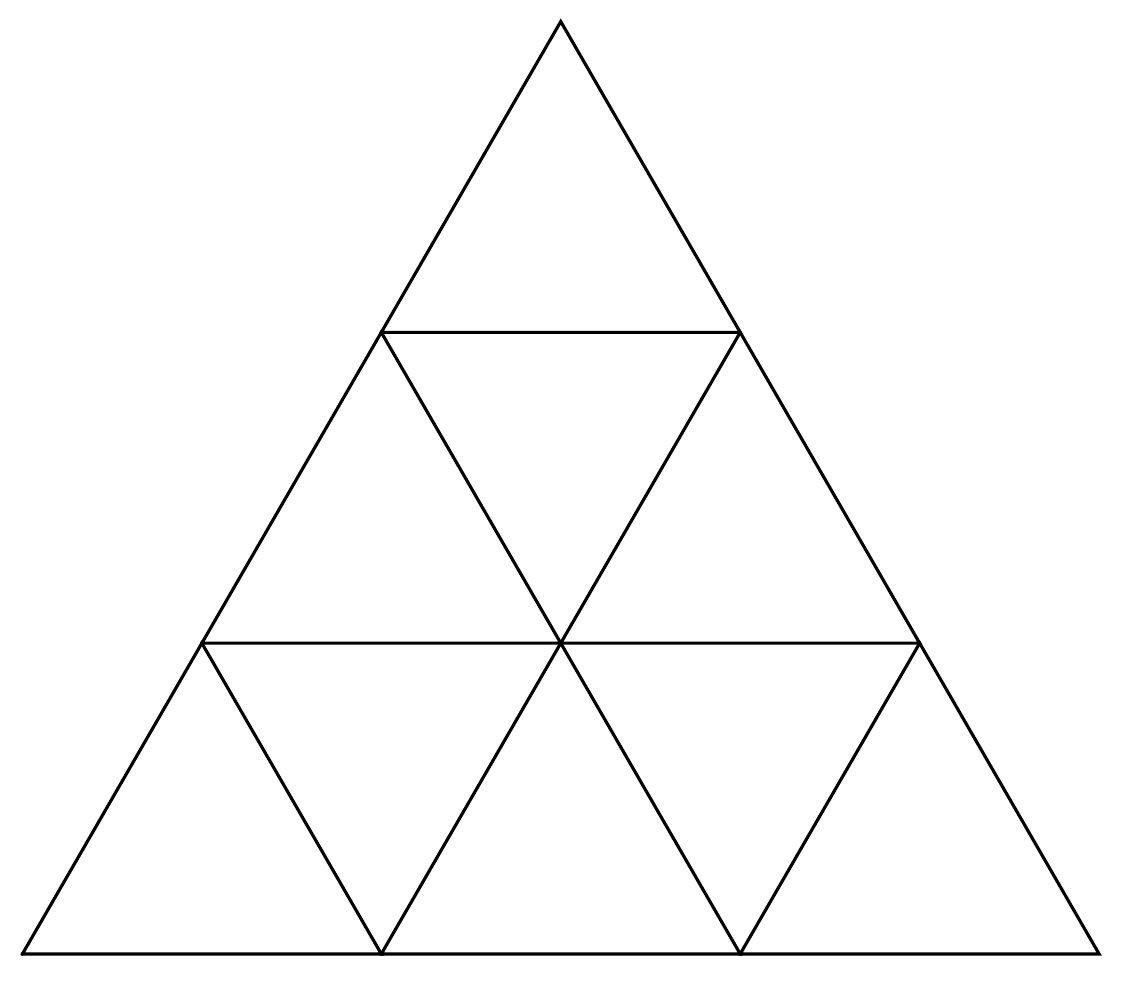}
\setlength{\unitlength}{1cm}
\begin{picture}(0,0) \thicklines
\put(-3.6,5.2){$h(q_0)$}
\put(-3.6,5.2){$h(q_0)$}
\put(-7.0,-0.15){$h(q_1)$}
\put(-7.0,-0.15){$h(q_1)$}
\put(-0.2,-0.15){$h(q_2)$}
\put(-0.2,-0.15){$h(q_2)$}
\put(-4.5,2.1){$\frac{h(q_0)+h(q_1)+h(q_2)}{3}$}
\put(-4.5,2.1){$\frac{h(q_0)+h(q_1)+h(q_2)}{3}$}
\put(-6.5,3.7){$\frac{8h(q_0)+4h(q_1)+3h(q_2)}{15}$}
\put(-6.5,3.7){$\frac{8h(q_0)+4h(q_1)+3h(q_2)}{15}$}
\put(-3.0,3.7){$\frac{8h(q_0)+3h(q_1)+4h(q_2)}{15}$}
\put(-3.0,3.7){$\frac{8h(q_0)+3h(q_1)+4h(q_2)}{15}$}
\put(-7.5,1.8){$\frac{4h(q_0)+8h(q_1)+3h(q_2)}{15}$}
\put(-7.5,1.8){$\frac{4h(q_0)+8h(q_1)+3h(q_2)}{15}$}
\put(-1.9,1.8){$\frac{4h(q_0)+3h(q_1)+8h(q_2)}{15}$}
\put(-1.9,1.8){$\frac{4h(q_0)+3h(q_1)+8h(q_2)}{15}$}
\put(-6.2,0.3){$\frac{3h(q_0)+8h(q_1)+4h(q_2)}{15}$}
\put(-6.2,0.3){$\frac{3h(q_0)+8h(q_1)+4h(q_2)}{15}$}
\put(-3.2,0.3){$\frac{3h(q_0)+4h(q_1)+8h(q_2)}{15}$}
\put(-3.2,0.3){$\frac{3h(q_0)+4h(q_1)+8h(q_2)}{15}$}
\end{picture}
\vspace{0.2cm}
\begin{center}
\textbf{Figure 1.4. Harmonic extension algorithm of $\mathcal{SG}_3$.}
\end{center}
\end{center}
\end{figure}

The \textit{normal derivative} of a function $u$ at a boundary point $q_i\in V_0$ is defined by
$$\partial_n(q_i)=\lim_{m\rightarrow\infty}r^{-m}\big(2u(q_i)-u(F_i^m q_{i+1})-u(F_i^m q_{i-1})\big)$$(cyclic notation $q_3=q_0$) providing the limit exists. For harmonic functions, these derivatives can be evaluated without taking limit.  We could localize the definition of normal derivative to any vertex in $V_*$. Let $x=F_wq_i$ be a boundary point of a $m$-cell $F_w\mathcal{SG}_l$. Define $\partial_n u(x)$ with respect to the cell $F_w\mathcal{SG}_l$ to be $r^{-m}\partial_n(u\circ F_w)(q_i)$. In this paper, we use the notations 
$\partial_n^\leftarrow,\partial_n^\rightarrow,\partial_n^\uparrow$ to represent the normal derivatives of different directions. In particular, $\partial_n^\uparrow u(q_0)=\partial_n u(q_0),\partial_n^\leftarrow u(q_1)=\partial_n u(q_1),\partial_n^\rightarrow u(q_2)=\partial_n u(q_2)$.
 For $u\in dom\Delta$, the sum of all normal derivatives of $u$ at each $x\in V_*\setminus V_0$ must vanish. This is called the \textit{matching condition}.

We have an analogue of the \textit{Gauss-Green's formula} in the fractal setting. Suppose $u\in dom\Delta$, then $\partial_n u(q_i)$ exists for all $q_i\in V_0$ and 
$$\mathcal{E}(u,v)=-\int_{\mathcal{SG}_l}(\Delta u)vd\nu+\sum_{q_i\in V_0}v(q_i)\partial_n u(q_i).$$
We also have a localized version of this formula,
$$\mathcal{E}_{A}(u,v)=-\int_{A}(\Delta u)vd\nu+\sum_{\partial A}v(x)\partial_n u(x)$$ for any \textit{simple set} $A$, which is defined as a finite union of cells.

Let $\Omega$ be a half, upper or lower domain in $\mathcal{SG}_l$. Consider the Dirichlet problem
\begin{equation}
\begin{cases}
\Delta u=0\text{ in } \Omega,\\
u|_{\partial \Omega}=f, f\in C(\partial \Omega).
\end{cases}
\end{equation}

\textbf{Proposition 1.1.} \textit{The Dirichlet problem (1.1) has a unique solution.}

\textit{Proof.} First, by Lemma 8.2 of [Ki5], if there exists a function $v\in dom\mathcal{E}$ such that $v|_{\partial \Omega}=f$, then a solution of (1.1) exists, which minimizes the energy on $\Omega$. For general case, notice that the set $dom\mathcal{E}|_{\partial \Omega}:=\{f\in C(\partial\Omega)|\exists v\in dom\mathcal{E}, v|_{\partial\Omega}=f \}$ is dense in $C(\partial\Omega)$, since $dom\mathcal{E}$ is dense in $C(\mathcal{SG}_l)$ and $\partial\Omega$ is a closed subset of $\mathcal{SG}_l$. Let $\{f_n\}$ be  a sequence of functions in  $ dom\mathcal{E}|_{\partial \Omega}$  converging uniformly to $f$, and $u_n$ be their corresponding solutions of (1.1). Then $\{u_n\}$ also uniformly converge to a function $u$ with $u|_{\partial \Omega}=f$ by the maximum principle for harmonic functions. It is easy to get that $u$ is harmonic in $\Omega$.

The uniqueness of the solution is an immediate consequence of the maximum principle.\hfill$\square$\\

\subsection{The organization of the paper.}

Throughout this paper, although in different situations, we always use the same symbol $\Omega$ to denote the domain and $X$ to denote the Cantor set boundary without causing any confusion. 

In Section 2, we solve the Dirichlet problem for the half domain in $\mathcal{SG}_3$. An explicit harmonic extension algorithm is provided. Let $f$ be the prescribed value on $\partial\Omega$. We only need to find the explicit formula for the values of the extension harmonic function $u$ on $V_1\cap\Omega$, since if we do so, then the value of $u$ in the $1$-cells contained in $\bar{\Omega}$ is determined by the harmonic extension algorithm, and then the problem of finding values of $u$ in the remaining region is essentially the same by dilation. An interesting phenomenon is that the solution could be expressed explicitly in terms of only a countably infinite set of points  which is dense in $\partial\Omega$.   We also characterize the energy estimate of solutions of finite energy in terms of their boundary values. 

We consider the analogous problem in the case of upper domain in $\mathcal{SG}_3$ in Section 3. Basing on the same reason, we need to find the explicit formula for finite number of crucial points then use dilation to continue. For the energy estimate, we also use the technique of Haar series expansion. But now we expand the boundary values with respect to a more natural probability measure rather than the normalized Hausdorff measure. 

In Section 4, we deal with the lower domain in $\mathcal{SG}$. Essentially the method is the same as before, but the situation is more complicated. We still obtain the explicit harmonic extension algorithm. However, it is unclear how to work out the energy estimate in term of the boundary values. 

Finally, we show our methods on the above three types of domains are still valid in the case of general $\mathcal{SG}_l$ and briefly state the outcomes. We present an intriguing correspondence between the normal derivatives and the boundary values of harmonic functions on the half domain of $\mathcal{SG}$, although we have no idea on how to extend it to general cases.

At the end of this section, we list some previous work on related topics. See [GKQS], [HKu], [J], [Ki4], [LRSU], [LS],  [OS] and the references therein. In particular, in [LRSU], some extension problems  on $\mathcal{SG}$ are studied which is to find a function with certain prescribed data such as values and derivatives at a finite set, that minimizes prescribed Sobolev types of norms. It is interesting to consider analogous problems on domains in this paper. We leave these as open problems for future research. The energy estimates considered in this paper characterize the restriction to the Cantor set boundary $X$ of functions of finite energy on $\Omega$. It is also interesting to characterize the traces on $X$ of functions in some other Sobolev spaces, such as $dom_{L^2}(\Delta^k)$ defined as $\{u\in L^2(\mathcal{SG}_l): \Delta^j u\in L^2(\mathcal{SG}_l), \forall j\leq k\}$. One could also consider how to extend a function of finite energy defined on $\Omega$ to a function of finite energy on the whole $\mathcal{SG}_l$, and analogous problems for other Sobolev spaces. Related problems are discussed in [LS], [GKQS]. The above mentioned Sobolev spaces are easily characterized in terms of expansions in eigenfunctions of the Laplacian, see [S2]. For half domains, as pointed out in [LS] in case of $\mathcal{SG}$ setting, essentially, there are no new eigenfunctions. A complete theory of the eigenspaces of the Laplacian on the upper domain in $\mathcal{SG}$ with $X$ equal to the bottom line segment is given in [Q].

\section{Dirichlet problem on the half domain of $\mathcal{SG}_3$}
In this section, we focus on solving the Dirichlet problem on the half domain of $\mathcal{SG}_3$. We will first give an extension algorithm for harmonic functions with continuous prescribed boundary values, then estimate the energies  of them in terms of  their boundary values. 

\subsection{Extension Algorithm} The domain $\Omega$ can be defined by a level set of an antisymmetric harmonic function, denoted by $h_a$, with boundary values
$h_a|_{\{q_0,q_1,q_2\}}=(0,1,-1)$, so that
$$
\Omega=\{x\in \mathcal{SG}_3\setminus V_0:h_a(x)>0\},$$
 and the boundary
$$\partial\Omega=\{q_1\}\cup X,\text{ with }X=\{x\in \mathcal{SG}_3:h_a(x)=0\}.$$
Let $\bar{\Omega}$ denote the closure of $\Omega$. It is easy to check that \[
\bar{\Omega}=F_1\mathcal{SG}_3\cup F_5\mathcal{SG}_3\cup F_0\bar{\Omega}\cup F_3\bar{\Omega}.
\]

As shown in Section 1, to solve the Dirichlet problem (1.1), we only need to find the explicit algorithm for the values of the harmonic function $u$ on $V_1\cap\Omega$.
For convenience, we use $x_\emptyset, y_\emptyset, z_\emptyset$ to represent the three ``crucial'' vertices in $V_1\cap\Omega$ with $$x_\emptyset=F_1q_2, y_\emptyset=F_1q_0, z_\emptyset=F_0q_1.$$ We also denote $p_\emptyset=F_3q_0$.

For $m\geq 0$, write \[\tilde{W}_m=\{0,3\}^m\text{ and }\tilde{W}_*=\bigcup\limits_{m=0}^{\infty}\tilde{W}_m.\] Obviously, $\tilde{W}_m\subset W_m$ and $\tilde{W}_*\subset W_*$. Denote \[x_w=F_wx_\emptyset, \text{ }y_w=F_wy_\emptyset,   \text{ }z_w=F_wz_\emptyset,\text{ }  p_w=F_wp_\emptyset\text{ for } w\in\tilde{W}_*.\] Obviously, $\{x_w, y_w, z_w\}_{w\in \tilde{W}_*}\subset V_*\cap \Omega$ and  $\{p_w\}_{w\in \tilde{W}_*}=V_*\cap X\setminus\{q_0\}$.

Now, we proceed to show how to determine the values of the harmonic function $u$ on $V_1\bigcap \Omega$ in terms of the boundary function $f$. From the matching condition at each vertex in $V_1\bigcap \Omega$, we have
\begin{equation}
\begin{cases}
\frac{7}{15}\partial_n^\leftarrow u(x_\emptyset)+2u(x_\emptyset)-u(y_\emptyset)-f(q_1)=0,\\
4u(y_\emptyset)-u(x_\emptyset)-u(z_\emptyset)-f(q_1)-f(p_\emptyset)=0,\\
\frac{7}{15}\partial_n^\leftarrow u(z_\emptyset)+2u(z_\emptyset)-u(y_\emptyset)-f(p_\emptyset)=0.
\end{cases}
\end{equation}
To make the equations (2.1) enough to determine the unknown, we need to represent  the normal derivatives at $x_\emptyset$ and $z_{\emptyset}$ in terms of $\{u(x_{\emptyset}),u(y_{\emptyset}),u(z_{\emptyset})\}$ and the boundary data $f$. 

We will prove that there exists a \textit{signed measure} on the boundary $\partial\Omega$ such that the normal derivative of $u$ at $q_1$ could be evaluated as the integral of $f$ with respect to this measure. Moreover, this signed measure is determined by the normal derivative of the antisymmetric harmonic function $h_a$ along the boundary $\partial\Omega$. See Figure 2.1 for the values of $h_a$ on $V_1\cap \bar{\Omega}$.\\

\begin{figure}[h]
\begin{center}
\includegraphics[width=2.5cm]{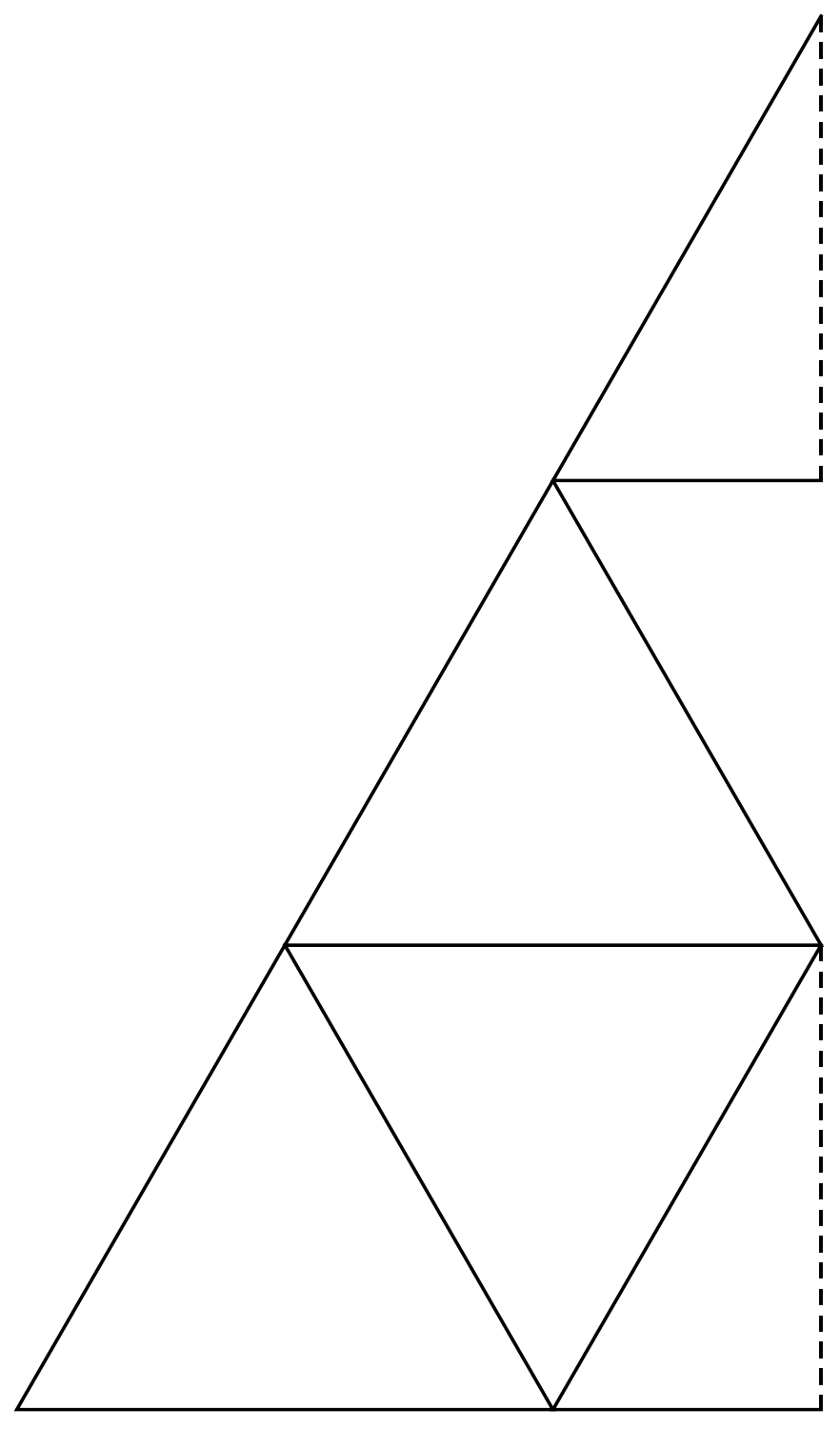}
\setlength{\unitlength}{1cm}
\begin{picture}(0,0) \thicklines
\put(-0.1,3.4){$0$}
\put(-0.1,0.5){$0$}
\put(-3,-0.1){$1$}
\put(-1.25,-0.3){$\frac{4}{15}$}
\put(-2.1,1.5){$\frac{1}{3}$}
\put(-1.4,3){$\frac{1}{15}$}
\end{picture}
\begin{center}
\vspace{0.4cm}
\textbf{Figure 2.1. The values of $h_a$ on $V_1\cap\bar{\Omega}$.}
\end{center}
\end{center}
\end{figure}

\textbf{Theorem 2.1.} \textit{Let $u$ be a solution of the Dirichlet problem (1.1). Then }
\begin{equation}
\partial^\leftarrow_n u(q_1)=3f(q_1)-\sum_{w\in \tilde{W}_*}\frac{6}{7}\mu_wf(p_w),
\end{equation}
\textit{where $\mu_w=\mu_{w_1}\mu_{w_2}\cdots \mu_{w_{|w|}}$ with $\mu_0=\frac{1}{7}$ and $\mu_3=\frac{4}{7}$. In addition, if $u\in dom\mathcal{E}_{\Omega}$, we have}
\[\mathcal{E}_{\Omega}(h_a,u)=\partial^\leftarrow_n u(q_1).\]

\textit{Proof.} Set $O_1=F_1{\mathcal{SG}_3}\bigcup F_5{\mathcal{SG}_3}$, and $O_m=\bigcup_{w\in \tilde{W}_*,|w|\leq m-1} F_wO_1$ for $m\geq 2$. Obviously, $\bar{\Omega}$ equals the closure of $\bigcup_{m\geq 1} O_m$. See Figure 2.2 for $O_1$ and $O_2$.

\begin{figure}[h]
\begin{center}
\includegraphics[width=3cm]{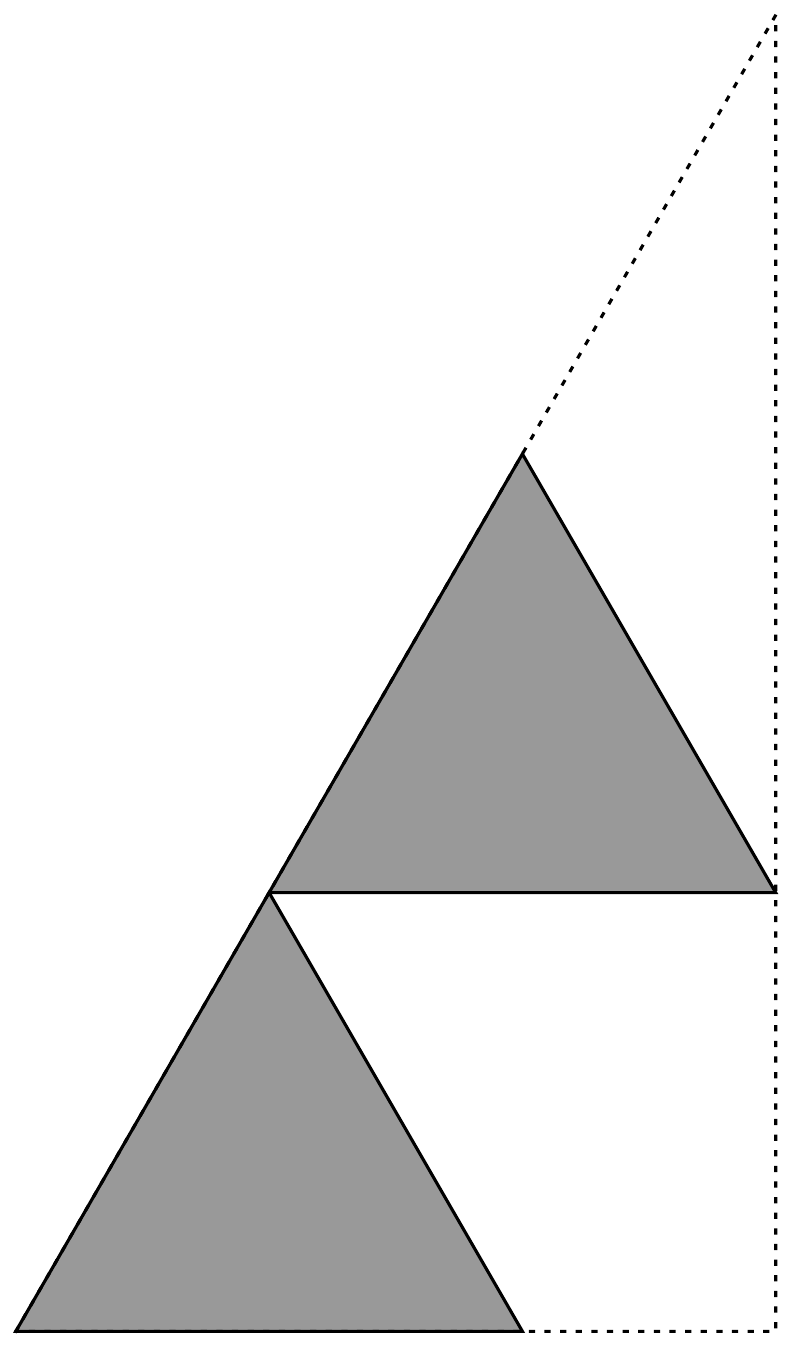}\hspace{1cm}
\includegraphics[width=3cm]{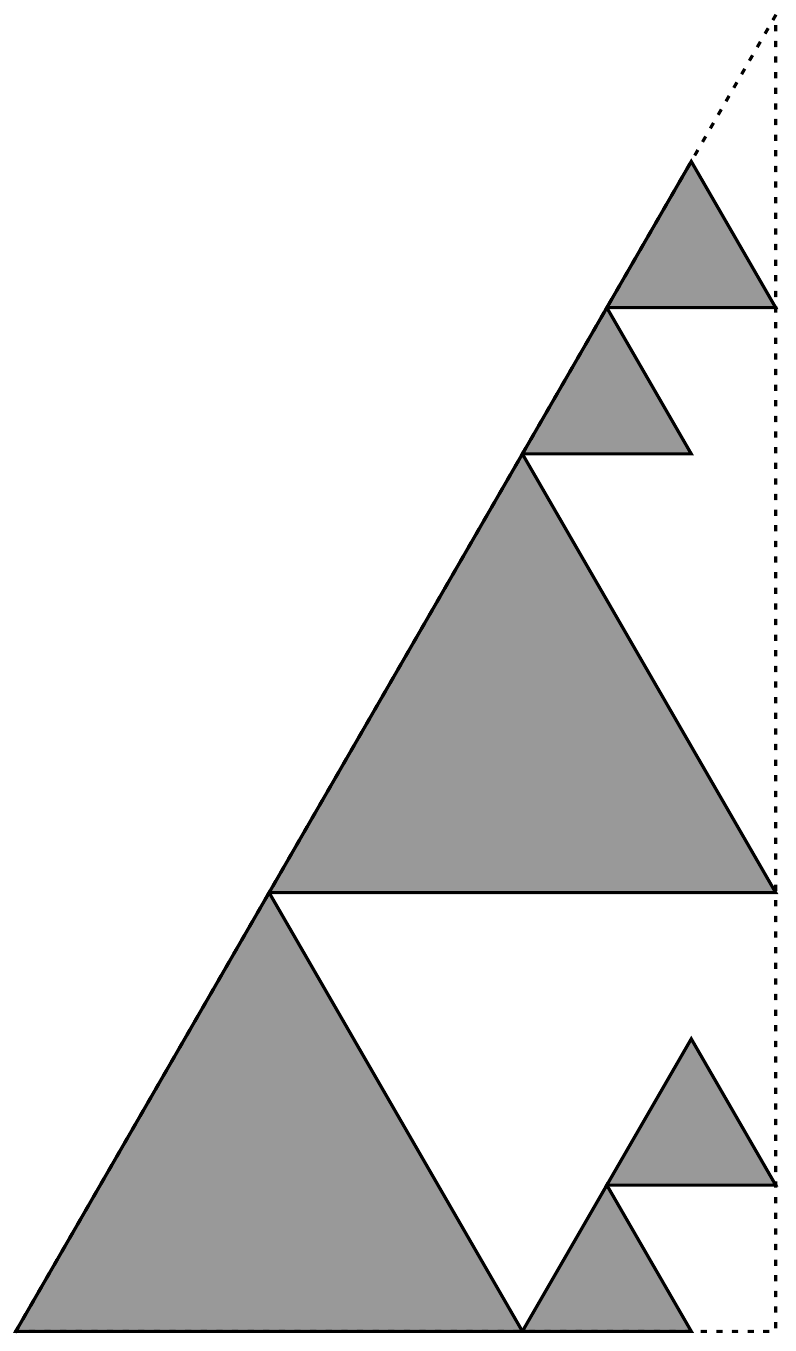}
\setlength{\unitlength}{1cm}
\begin{picture}(0,0) \thicklines
\put(-7.6,0){$q_1$}
\put(-5.4,-0.1){$x_\emptyset$}
\put(-6.65,1.75){$y_\emptyset$}
\put(-4.35,1.7){$p_\emptyset$}
\put(-5.65,3.45){$z_\emptyset$}

\put(-3.5,0){$q_1$}
\put(-1.4,-0.1){$x_\emptyset$}
\put(-2.55,1.75){$y_\emptyset$}
\put(-0.25,1.7){$p_\emptyset$}
\put(-1.58,3.4){$z_\emptyset$}

\put(-0.7,-0.1){$x_3$}
\put(-1.25,0.65){$y_3$}
\put(-0.25,0.6){$p_3$}
\put(-0.7,1.2){$z_3$}

\put(-0.65,3.3){$x_0$}
\put(-1.25,3.95){$y_0$}
\put(-0.25,3.9){$p_0$}
\put(-0.95,4.5){$z_0$}
\end{picture}
\begin{center}
\vspace{0.4cm}
\textbf{Figure 2.2. Simple sets $O_1,O_2$.}
\end{center}
\end{center}
\end{figure}

Applying the local Gauss-Green's formula on $O_m$, we get
\[\begin{aligned}
\mathcal{E}_{O_m}(h_a,u)
=&\partial_n^\leftarrow h_a(q_1)f(q_1)+\sum_{w\in \tilde{W}_*, |w|\leq m-1} \partial_n^\rightarrow h_a(p_w)f(p_w)\\& +\sum_{w\in \tilde{W}_{m-1}}\left(\partial_n^\rightarrow h_a(x_w)u(x_w)+\partial_n^\uparrow h_a(z_w)u(z_w)\right).
\end{aligned}\]
 It is easy to calculate the normal derivatives of $h_a$ at $p_w, x_w, z_w$, 
\[\partial_n^\rightarrow h_a(p_w)=-\frac{6}{7}\mu_w,\quad \partial_n^\rightarrow h_a(x_w)=-\frac{12}{7}\mu_w,\quad \partial_n^\uparrow h_a(z_w)=-\frac{3}{7}\mu_w.\]
So we have the estimate that
\[\begin{aligned}&\Big|\mathcal{E}_{O_m}(h_a,u)-\big(3f(q_1)-\sum_{w\in \tilde{W}_*}\frac{6}{7}\mu_wf(p_w)\big)\Big|
\\\leq &\Big|\sum_{w\in \tilde{W}_{m-1}} \big(\partial_n^\rightarrow h_a(x_w)+\partial_n^\uparrow h_a(z_w)\big)+\sum_{w\in \tilde{W}_*,|w|\geq m} \partial_n^\rightarrow h_a(p_w)\Big|\cdot\|f\|_\infty
\\=&\frac{30}{7}(\frac{5}{7})^{m-1}\|f\|_{\infty}.\end{aligned}\]
Thus if $u\in dom \mathcal{E}_\Omega$, by taking the limit we have $\mathcal{E}_\Omega (h_a,u)=3f(q_1)-\sum_{w\in \tilde{W}_*}\frac{6}{7}\mu_wf(p_w)$. 

For the rest part of the theorem, we introduce a sequence of harmonic functions $\{u_n\}_{n\geq 0}$  which are piecewise constant on $X$, defined as $u_n|_{F_{\tau}X}=f(p_{\tau}),\forall \tau\in \tilde{W}_n$. The existence of such functions is ensured by Proposition 1.1. Moreover, it is easy to check that $u_n$ uniformly converges to $u$ by the maximum principle for harmonic functions. Applying the Gauss-Green's formula, we have
\[\begin{aligned}\mathcal{E}_{O_m}(h_a,u_n)=&\partial_n^\leftarrow u_n(q_1)h_a(q_1)+\sum_{w\in \tilde{W}_*,|w|\leq m-1} \partial_n^\rightarrow u_n(p_w)h_a(p_w)\\& +\sum_{w\in \tilde{W}_{m-1}}\left(\partial_n^\rightarrow u_n(x_w)h_a(x_w)+\partial_n^\uparrow u_n(z_w)h_a(z_w)\right)\\=&\partial_n^\leftarrow u_n(q_1) +\sum_{w\in \tilde{W}_{m-1}}\left(\partial_n^\rightarrow u_n(x_w)h_a(x_w)+\partial_n^\uparrow u_n(z_w)h_a(z_w)\right).\end{aligned}\]

Noticing that for fixed $n$ and $\tau\in \tilde{W}_n$, we have $\partial_n^\rightarrow u_n(x_{\tau w})$, $\partial_n^\uparrow u_n(z_{\tau w})$ taking the same sign with $\sum_{w\in \tilde{W}_*,|w|=m} \big(\partial_n^\rightarrow u(x_{\tau w})+\partial_n^\uparrow u(z_{\tau w})\big)$ uniformly bounded as $u_n\circ F_\tau=c_1+c_2h_a$ for some constants $c_1,c_2$. In addition, $h_a(x_w)$ and $h_a(z_w)$ converge uniformly to $0$ as $|w|\to\infty$. Thus, letting $m\to\infty$, we get
\[\begin{aligned}\mathcal{E}_{\Omega}(h_a,u_n)&=\lim_{m\to\infty}\Big(\partial_n^\leftarrow u_n(q_1)+\sum_{\tau\in \tilde{W}_n}\sum_{w\in \tilde{W}_m} \big(\partial_n^\rightarrow u(x_{\tau w})h_a(x_{\tau w})+\partial_n^\uparrow u(z_{\tau w})h_a(z_{\tau w})\big)\Big)\\&=\partial_n^\leftarrow u_n(q_1).\end{aligned}\] 
Combining this equality with the first part of the proof, we then have
\[\partial_n^\leftarrow u_n(q_1)=3u_n(q_1)-\sum_{w\in \tilde{W}_*}\frac{6}{7}\mu_w u_n(p_w),\]
Taking $n\to\infty$, we  get $(2.2)$.
$\hfill\square$

\textbf{Remark.} One can regard the signed measure $3\delta_{q_1}-\sum_{w\in \tilde{W}_*}\frac{6}{7}\mu_w\delta_{p_w}$ as the normal derivative of $h_a$ on $\partial \Omega$. In this opinion, Theorem 2.1 is just a result of  the extended ``Guass-Green's formula'' acting on $h_a$ and $u$. 

In the following, we denote $\mu$ the probability measure $\sum_{w\in \tilde{W}_*}\frac{2}{7}\mu_w\delta_{p_w}$ on $X$. Thus we could write
\begin{equation}
\partial_n^\leftarrow u(q_1)=3f(q_1)-3\int_{X} fd\mu.
\end{equation}

Now, we have enough information to calculate the values $u(x_{\emptyset}),u(y_{\emptyset}),u(z_{\emptyset})$.

\textbf{Theorem 2.2.(Extension Algorithm)} \textit{ There exists a unique solution of the Dirichlet problem (1.1). In addition, we have }
\begin{eqnarray}
& u(x_\emptyset)=\frac{4}{15} f(q_1)+\frac{1}{15} f(p_\emptyset)+\frac{1}{30} \int_{X} f\circ F_0d\mu+\frac{19}{30}\int_{X} f\circ F_3d\mu,\\
& u(y_\emptyset)=\frac{1}{3} f(q_1)+\frac{1}{3} f(p_\emptyset)+\frac{1}{6} \int_{X} f\circ F_0d\mu+\frac{1}{6}\int_{X} f\circ F_3d\mu,\\
& u(z_\emptyset)=\frac{1}{15} f(q_1)+\frac{4}{15} f(p_\emptyset)+\frac{19}{30} \int_{X} f\circ F_0d\mu+\frac{1}{30}\int_{X} f\circ F_3d\mu.
\end{eqnarray}
\textit{Proof.} The existence and uniqueness of a solution of (1.1) has been shown in Proposition 1.1. Taking $\partial_n^{\leftarrow} u(x_\emptyset)=\frac{15}{7}\partial_n^{\leftarrow} (u\circ F_3)(q_1)=\frac{15}{7}(3u(x_\emptyset)-3\int_X f\circ F_3d\mu)$ and $\partial_n^{\leftarrow} u(z_\emptyset)=\frac{15}{7}\partial_n^{\leftarrow} (u\circ F_0)(q_1)=\frac{15}{7}(3u(z_\emptyset)-3\int_X f\circ F_0d\mu)$ into (2.1), and solving the system of linear equations, we get (2.4), (2.5) and (2.6). $\hfill\square$

\subsection{Energy Estimate}
In Theorem 2.2, we have shown that the harmonic function $u$ could be explicitly determined by its values at only countably infinite vertices $\{q_1\}\bigcup\{p_w\}_{w\in \tilde{W}_*}$. It is natural to hope that the energy estimate of $u$ also depends on the same values as well.

\textbf{Theorem 2.3.} \textit{Let $f\in C(\partial\Omega)$, write }
\[
Q(f)=\big(f(q_1)-f(p_\emptyset)\big)^2+\sum_{w\in \tilde{W}_*}(\frac{15}{7})^{|w|}\Big(\big(f(p_w)-f(p_{w0})\big)^2+\big(f(p_w)-f(p_{w3})\big)^2\Big).
\]
\textit{Then we have the energy estimate that}
\[C_1 Q(f) \leq\mathcal{E}(u)\leq C_2 Q(f),\] 
\textit{where $C_1,C_2$ are two positive constants independent of $f$.}

\textit{Proof.} Notice that 
\[\begin{aligned}
\mathcal{E}_{O_1} (u)=&\frac{15}{7}\Big(\big(f(q_1)-u(x_\emptyset)\big)^2+\big(f(q_1)-u(y_\emptyset)\big)^2+\big(u(x_\emptyset)-u(y_\emptyset)\big)^2\\
&+\big(f(p_\emptyset)-u(y_\emptyset)\big)^2+\big(f(p_\emptyset)-u(z_\emptyset)\big)^2+\big(u(y_\emptyset)-u(z_\emptyset)\big)^2\Big)\\
\geq&\frac{45}{28}\big(f(q_1)-f(p_{\emptyset})\big)^2,
\end{aligned}\]
where the equality holds when $u(x_\emptyset)=\frac{3}{4}f(q_1)+\frac{1}{4}f(p_\emptyset),u(y_\emptyset)=\frac{1}{2}f(q_1)+\frac{1}{2}f(p_\emptyset)$ and $u(z_\emptyset)=\frac{1}{4}f(q_1)+\frac{3}{4}f(p_\emptyset).$
 Similarly, for any $w\in\tilde{W}_*$, we also have
\[\mathcal{E}_{F_wO_1}(u)+\mathcal{E}_{F_{w0}O_1}(u)\geq c_1(\frac{15}{7})^{|w|}\big(f(p_w)-f(p_{w0})\big)^2,\]
and
\[\mathcal{E}_{F_wO_1}(u)+\mathcal{E}_{F_{w3}O_1}(u)\geq c_2(\frac{15}{7})^{|w|}\big(f(p_w)-f(p_{w3})\big)^2,\]
where $c_1,c_2$ are suitable positive constants.  Thus, we have
\[\begin{aligned}
\mathcal{E}_{\Omega} (u)&=\sum_{w\in \tilde{W}_*}\mathcal{E}_{F_wO_1}(u)=\mathcal{E}_{O_1}(u)+\sum_{w\in \tilde{W}_*}\big(\mathcal{E}_{F_{w0}O_1}(u)+\mathcal{E}_{F_{w3}O_1}(u)\big)\\
&=\frac{1}{3}\mathcal{E}_{O_1}(u)+\frac{1}{3}\sum_{w\in \tilde{W}_*}\big(2\mathcal{E}_{F_wO_1}(u)+\mathcal{E}_{F_{w0}O_1}(u)+\mathcal{E}_{F_{w3}O_1}(u)\big)\\
&\geq \frac{1}{3} \min\{c_1,c_2,\frac{45}{28}\}\cdot Q(f) .
\end{aligned}\]

Conversely, we assume without loss of generality that  $\mathcal{E}_\Omega(u)<\infty$, otherwise there is nothing  to prove. Consider  a piecewise harmonic function $v$ defined on $\bar{\Omega}$ assuming values $v|_{\partial \Omega}=f$ and $ v(x_w)=v(z_w)=f(p_w),\forall w\in\tilde{W_*}$, which is harmonic in $F_wO_1,\forall w\in\tilde{W_*}$. See Figure 2.3 for the value of this function.

\begin{figure}[h]
\begin{center}
\includegraphics[width=3.5cm]{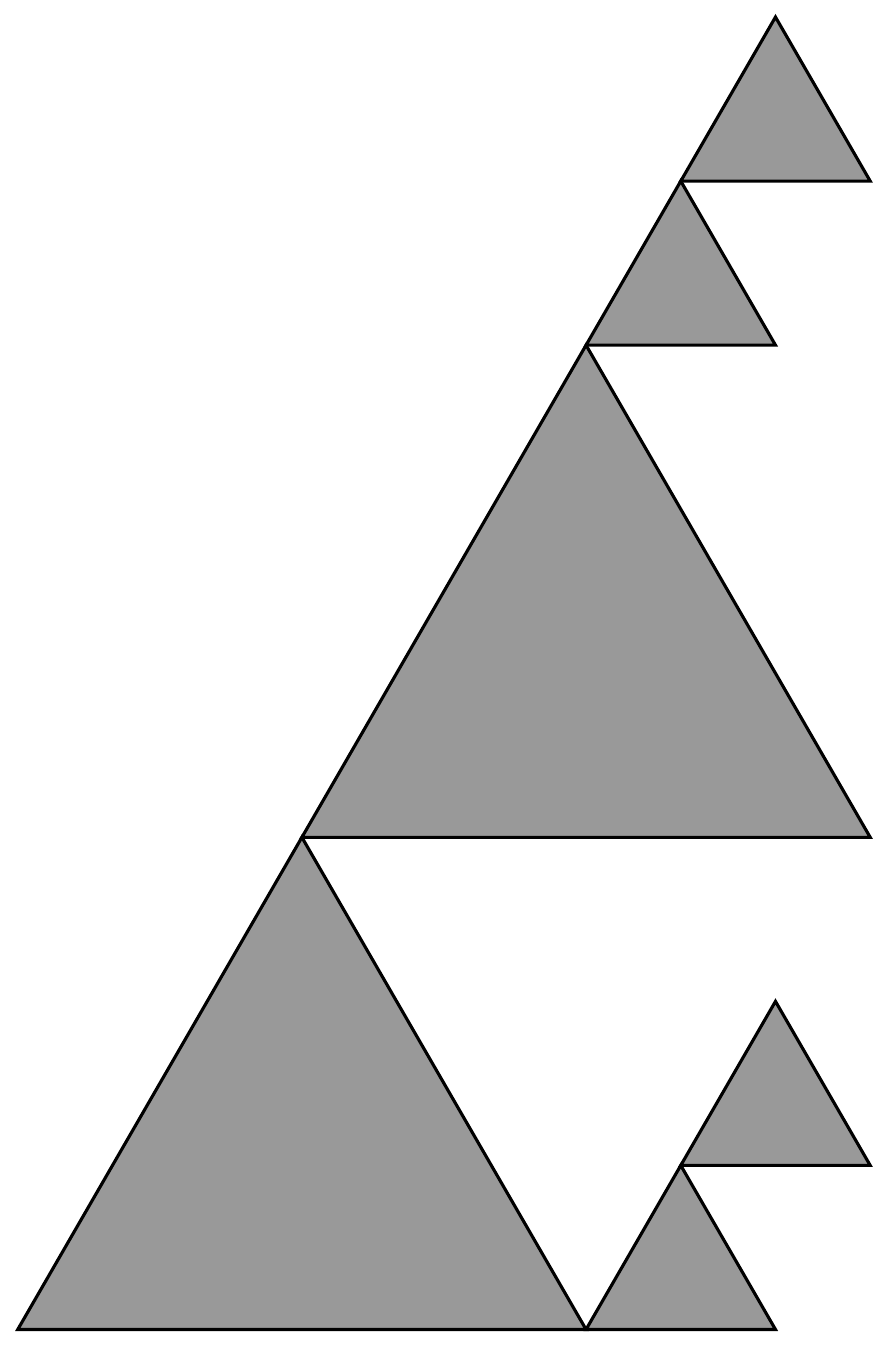}
\setlength{\unitlength}{1cm}
\begin{picture}(0,0) \thicklines
\put(-4.45,0){$f(q_1)$}
\put(-1.9,-0.25){$f(p_\emptyset)$}
\put(-2.3,3.9){$f(p_\emptyset)$}
\put(-0.3,1.95){$f(p_\emptyset)$}

\put(-0.95,5.25){$f(p_0)$}
\put(-0.7,3.9){$f(p_0)$}
\put(-0.35,4.5){$f(p_0)$}

\put(-0.7,0){$f(p_3)$}
\put(-1.0,1.43){$f(p_3)$}
\put(-0.3,0.65){$f(p_3)$}
\end{picture}
\begin{center}
\vspace{0.2cm}
\textbf{Figure 2.3. The values of $v$ on $O_2$.}
\end{center}
\end{center}
\end{figure}

It is easy to calculate the energy of $v$,
\[
\begin{aligned}
\mathcal{E}_\Omega(v)=&\sum_{w\in \tilde{W}_*}\mathcal{E}_{F_wO_1}(v)
\\=&\frac{15}{4}\big(f(q_1)-f(p_\emptyset)\big)^2\\&+\sum_{w\in \tilde{W}_*} \frac{7}{4}(\frac{15}{7})^{|w|+2}\Big(\big(f(p_w)-f(p_{w0})\big)^2+\big(f(p_w)-f(p_{w3})\big)^2\Big)
\\\leq&\frac{225}{28}Q(f).
\end{aligned}\]
On the other hand, $\mathcal{E}_{\Omega}(v)\geq\mathcal{E}_{\Omega}(u)$, as harmonic functions minimize the energy. 
\hfill$\square$

\section{Dirichlet problem on Upper Domains of $\mathcal{SG}_3$}
In this section, we deal with the Dirichlet problem on upper domains of $\mathcal{SG}_3$.
Prescribe that the boundary vertices $q_0,q_1,q_2\in \mathbb{R}^2$ take the following coordinates,
\[q_0=(\frac{1}{\sqrt{3}},1),\text{ }q_1=(0,0),\text{ }q_2=(\frac{2}{\sqrt{3}},0).\]
Then for each $0<\lambda\leq 1$, define the upper domain 
\[\Omega_\lambda=\{(x,y)\in \mathcal{SG}_3\setminus V_0|y>1-\lambda\},\]
together with the boundary
\[\partial\Omega_\lambda=\{q_0\}\cup X_\lambda,\text{ with }X_\lambda=\{(x,y)\in \mathcal{SG}_3|y=1-\lambda\}.\]
See Figure 3.1 for an illustration.
Denote $\bar{\Omega}_\lambda=\Omega_\lambda\bigcup \partial\Omega_\lambda$ the closure of $\Omega_\lambda$. In the following context, we write $X$ instead of $X_\lambda$ when there is no confusion. 

\begin{figure}[h]
\begin{center}
\includegraphics[width=4.5cm]{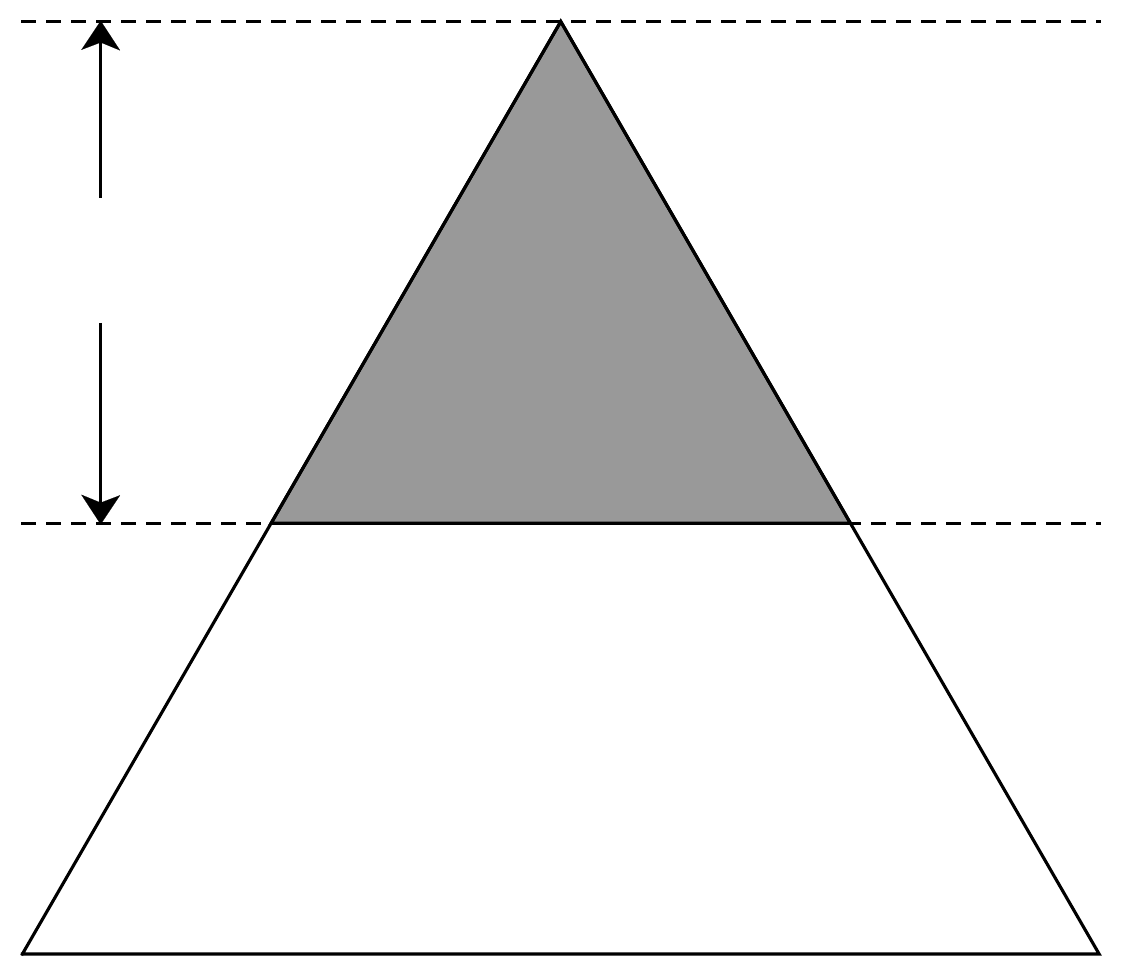}
\setlength{\unitlength}{1cm}
\begin{picture}(0,0) \thicklines
\put(-4.45,2.75){$\lambda$}
\put(-2.65,3.95){$q_0$}
\put(-4.95,0){$q_1$}
\put(-0.3,0){$q_2$}
\end{picture}
\vspace{0.2cm}
\begin{center}
\textbf{Figure 3.1. The upper domain. }
\end{center}
\end{center}
\end{figure}

For $0<\lambda \leq 1$, there is a unique representation
\begin{equation}
\lambda=\sum_{k=1}^{\infty} \iota_k3^{-m_k},
\end{equation}
with an integer sequence $0<m_1<m_2<\cdots$, and $\iota_k=1$ or $2$.
Denote
$$R\lambda=\sum_{k=2}^{\infty} \iota_k 3^{-(m_k-m_1)}.$$
Inductively, write $$\lambda_n=R^n\lambda=\sum_{k=n+1}^{\infty} \iota_k 3^{-(m_k-m_n)}.$$
Set $\lambda_0=\lambda$ and $m_0=0$ for the sake of formality.

It is easy to check the following  relationship between $\Omega_{\lambda_n}$ and $\Omega_{\lambda_{n+1}}$,
\begin{equation}
\bar{\Omega}_{\lambda_n}=
\begin{cases}
F_0^{m_{n+1}-m_n-1}(F_0\mathcal{SG}_3\cup F_4 \bar{\Omega}_{\lambda_{n+1}} \cup F_5 \bar{\Omega}_{\lambda_{n+1}})\\
\qquad\qquad\qquad\qquad\qquad\qquad\qquad\qquad \text{if }\iota_{n+1}=1,\\

F_0^{m_{n+1}-m_n-1}(F_0\mathcal{SG}_3\cup F_4\mathcal{SG}_3 \cup F_5\mathcal{SG}_3\\
\qquad\qquad\qquad\cup F_1 \bar{\Omega}_{\lambda_{n+1}} \cup F_2 \bar{\Omega}_{\lambda_{n+1}} \cup F_3 \bar{\Omega}_{\lambda_{n+1}})\\
\qquad\qquad\qquad\qquad\qquad\qquad\qquad\qquad \text{if }\iota_{n+1}=2.
\end{cases}
\end{equation}
For $1\leq i\leq 5$, write $p^\lambda_i=F_0^{m_1-1}F_iq_0$.  We omit the superscript $\lambda$ when there is no confusion caused. See Figure 3.2 for an illustration.

\begin{figure}[h]
\begin{center}
\includegraphics[width=5cm]{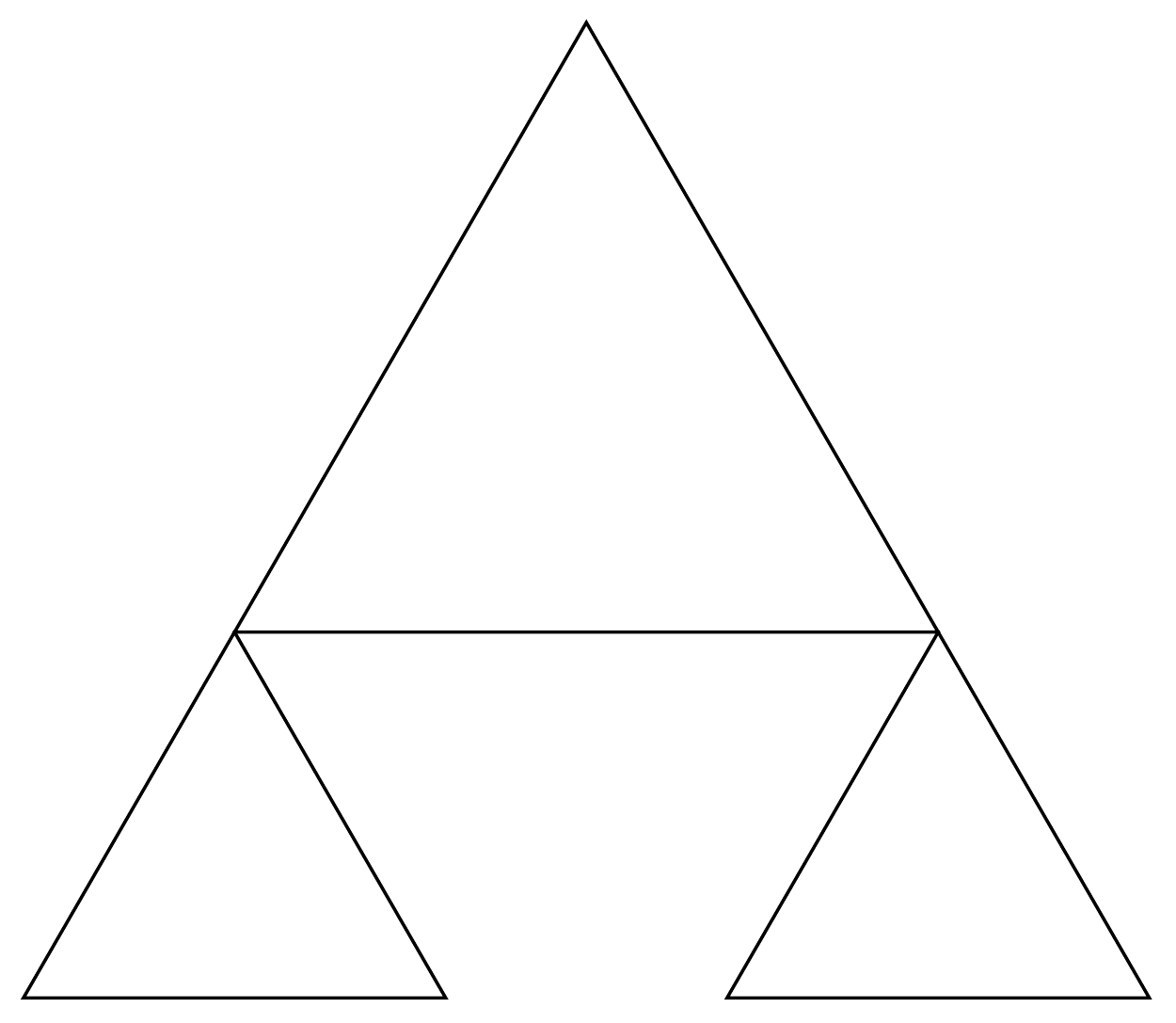}\hspace{1cm}
\includegraphics[width=5cm]{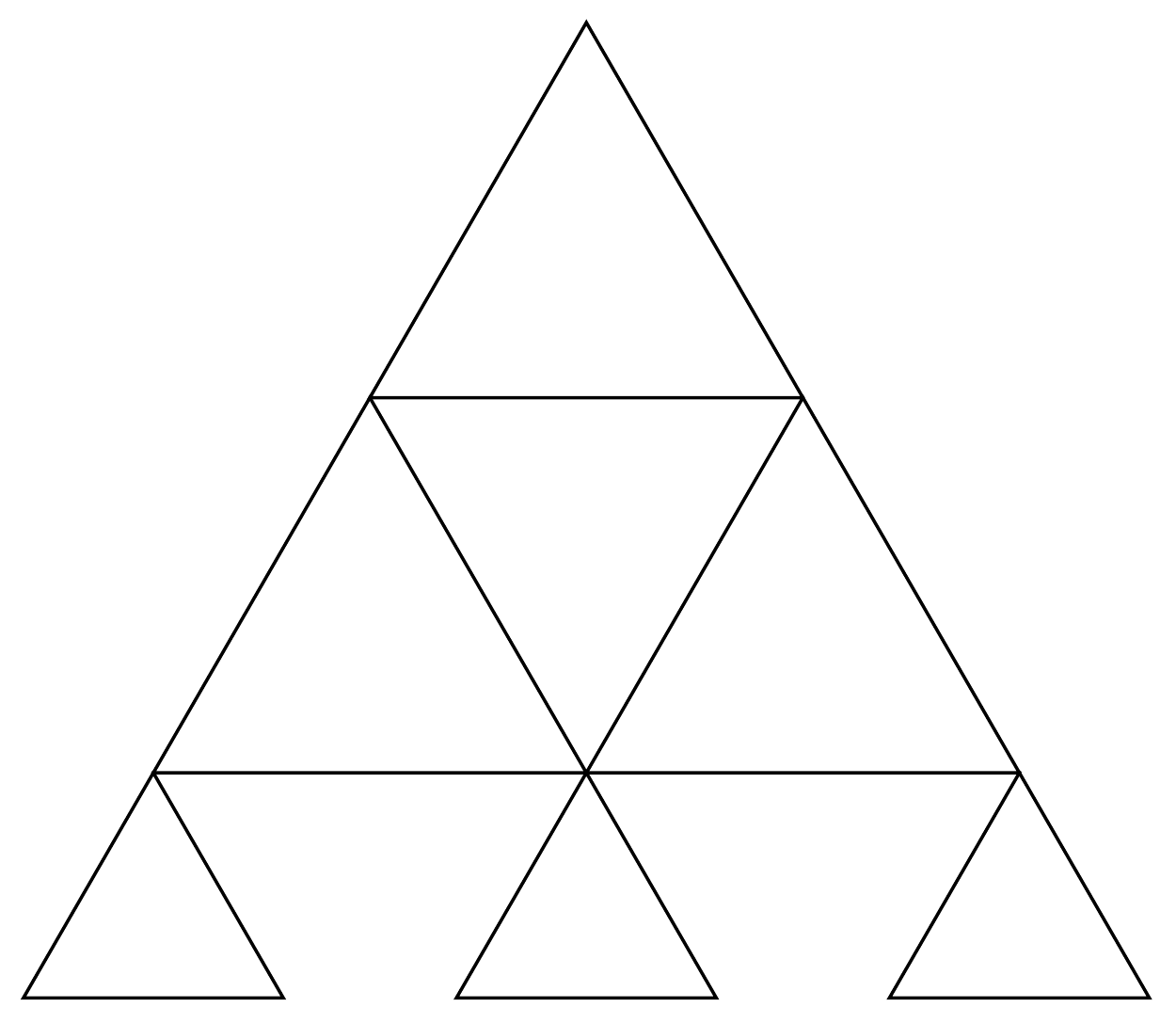}
\setlength{\unitlength}{1cm}
\begin{picture}(0,0) \thicklines
\put(-9.0,4.35){$q_0$}
\put(-10.7,1.7){$p_5$}
\put(-7.3,1.7){$p_4$}
\put(-11.3,0.3){$F_0^{m_1-1}F_5\Omega_{\lambda_1}$}
\put(-8.3,0.3){$F_0^{m_1-1}F_4\Omega_{\lambda_1}$}

\put(-2.9,4.35){$q_0$}
\put(-1.75,2.65){$p_4$}
\put(-4.05,2.65){$p_5$}
\put(-4.9,1.1){$p_1$}
\put(-0.85,1.1){$p_2$}
\put(-2.8,1.2){$p_3$}
\put(-5.6,0.2){$F_0^{m_1-1}F_1\Omega_{\lambda_1}$}
\put(-1.75,0.2){$F_0^{m_1-1}F_2\Omega_{\lambda_1}$}
\put(-3.7,0.4){$F_0^{m_1-1}F_3\Omega_{\lambda_1}$}

\put(-9.2,-0.4){\textbf{$Case 1:\iota_1=1$}}
\put(-3.6,-0.4){\textbf{$Case 2:\iota_1=2$}}
\end{picture}
\vspace{0.5cm}
\begin{center}
\textbf{Figure 3.2. The relationship between $\Omega_{\lambda_0}$ and $\Omega_{\lambda_1}$. }
\end{center}
\end{center}
\end{figure}

For convenience, let \[W^\lambda_n=\prod_{k=1}^{n} S_{\iota_k}, \text{ and } W^\lambda_*=\bigcup\limits_{n\geq 0} W^\lambda_n\text{, with $S_1=\{4,5\}$ and $S_2=\{1,2,3\}$},\] and write 
\[F^\lambda_w=F_0^{m_1-1}F_{w_1}\cdots F_0^{m_{n}-m_{n-1}-1}F_{w_n},\text{ for }w\in W^\lambda_n.\]
Then $\forall n\geq 0$, $$X=\bigcup_{w\in W_n^{\lambda}} X_w$$
where $X_w=F_w^\lambda\mathcal{SG}_3\cap X$.
It is easy to see that for $\lambda$ not a triadic rational, $X$ is homeomorphic to the space
$\Sigma^\lambda=\prod_{k=1}^{\infty} S_{\iota_k}$ equipped with the product topology. Otherwise, $X$ is a union of finite segments. 

Here we give an example to help readers to get familiar with the notations. 

\textbf{Example.} We plot the area $\Omega_\lambda$ for $\lambda=0.39=3^{-1}+3^{-2}+2\cdot 3^{-3}+\cdots$. See Figure 3.3.
In this case, $\iota_1=1,\iota_2=1,\iota_3=2$ and 
$$W^{\lambda}_1=\{4,5\}, W^{\lambda}_2=\{4,5\}^2=\{44, 45, 54, 55\}, W^{\lambda}_3=\{4,5\}^2\times \{1,2,3\}.$$

\begin{figure}[htbp]
\begin{center}
\includegraphics[width=4cm]{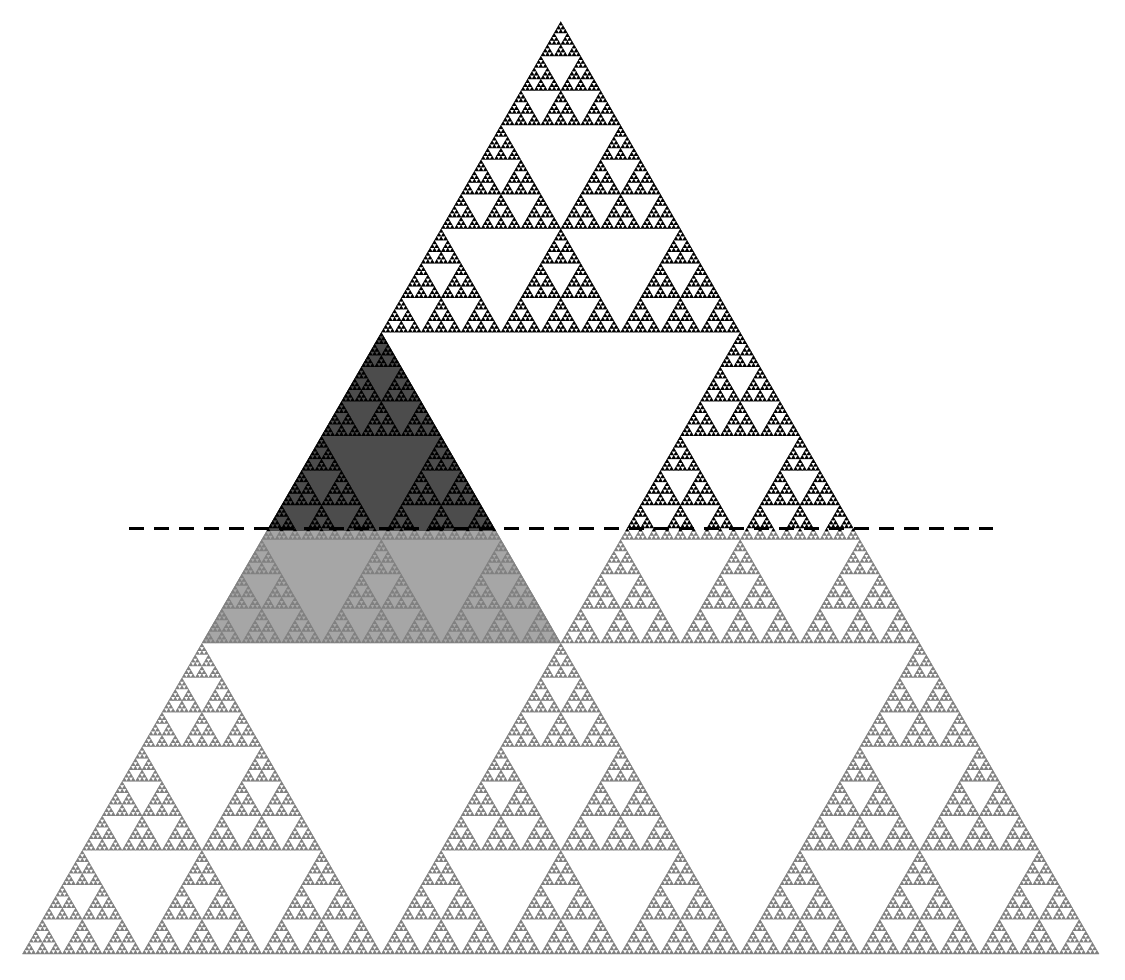}
\includegraphics[width=4cm]{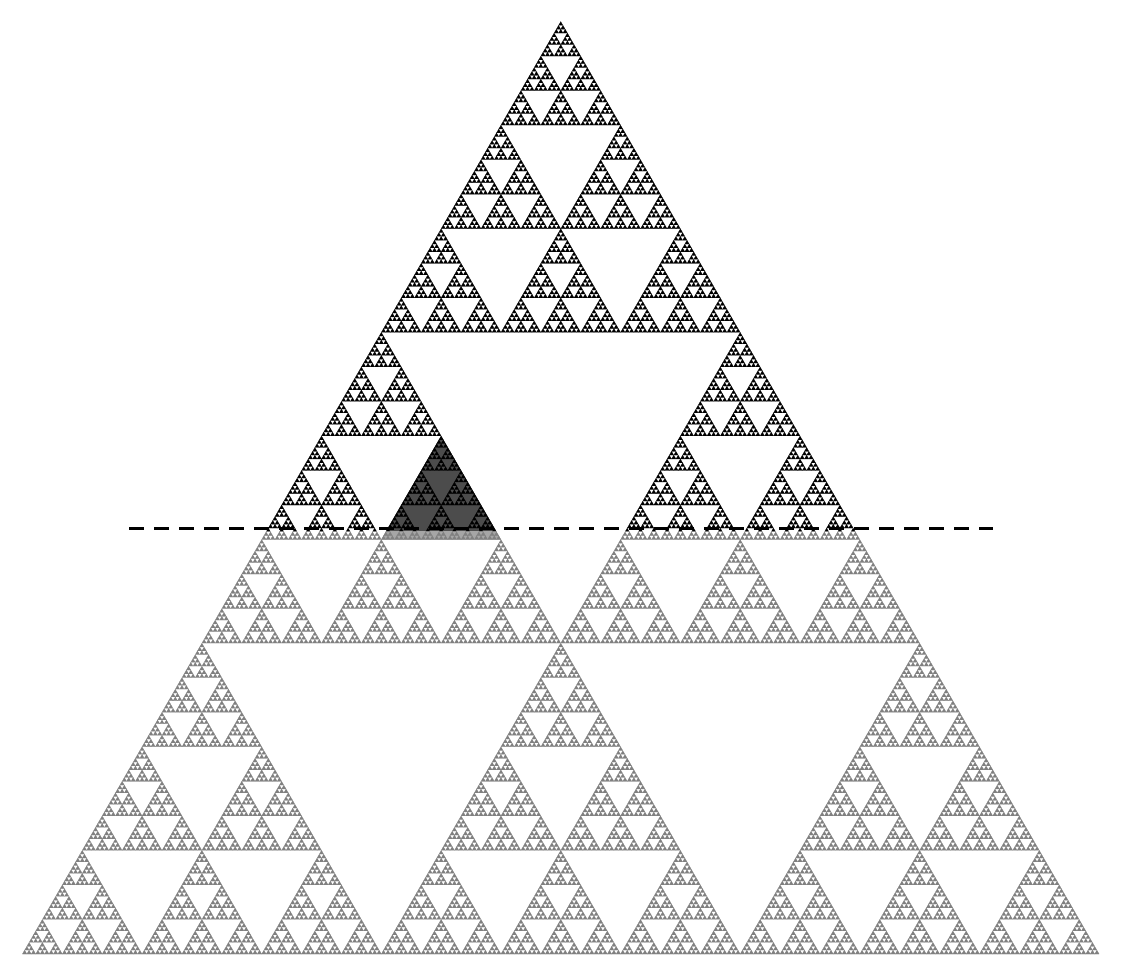}
\includegraphics[width=4cm]{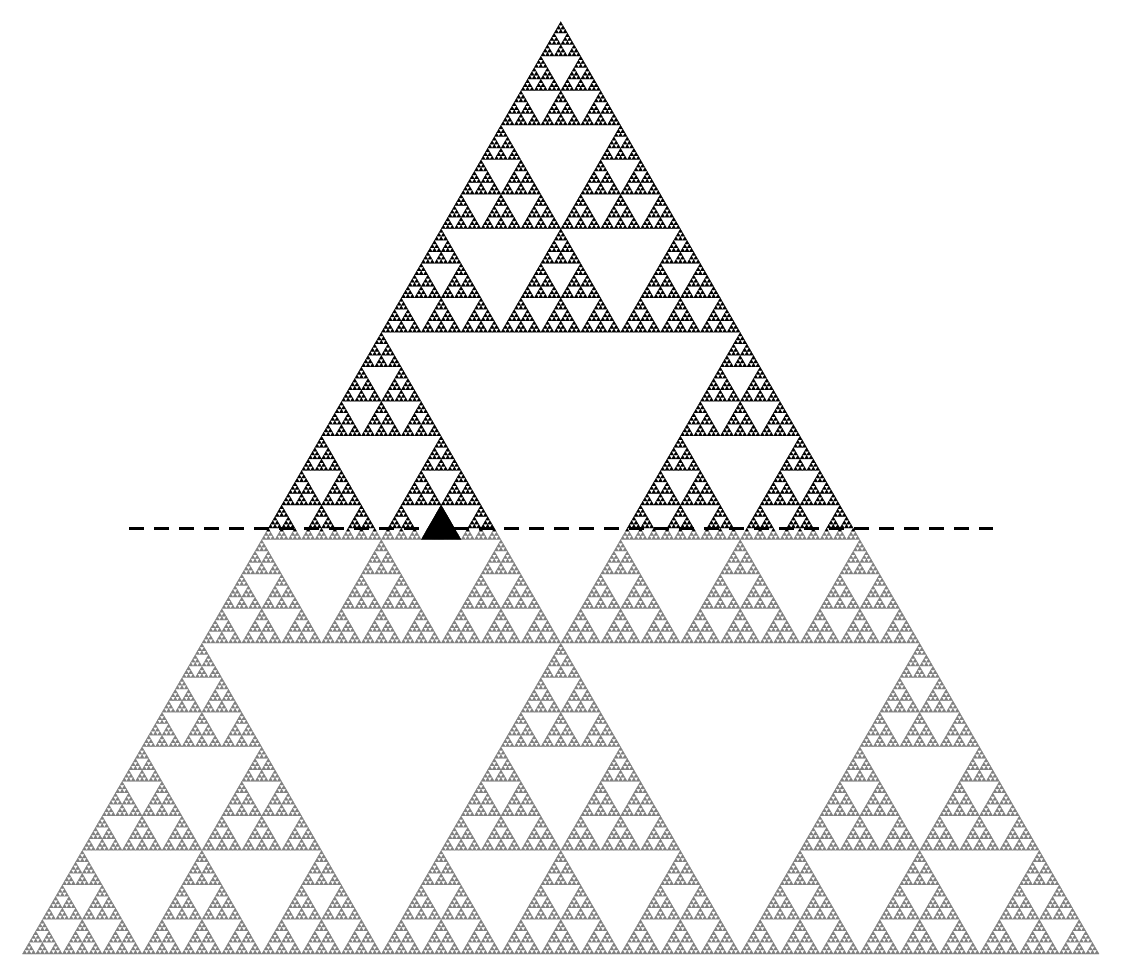}
\end{center}
\begin{center}
\textbf{Figure 3.3. The upper domain $\Omega_{0.39}$. The shaded regions are $F^\lambda_5\Omega_{\lambda_1}$, $F^\lambda_{54}\Omega_{\lambda_2}$, $F^\lambda_{543}\Omega_{\lambda_3}$. }
\end{center}
\end{figure}

\subsection{Extension Algorithm}

We still use $f$ to denote the boundary data on $\partial\Omega_\lambda$ and $u$ the harmonic solution of the Dirichlet problem (1.1). We only need to find an explicit algorithm for the values of  $u$ at $p_i$'s since if we do so, the problem of finding values of $u$ in the remaining region is essentially the same after dilation.  

From the  matching condition at each vertex $p_i$, we have the following system of equations.

Case 1($\iota_1=1$):
\begin{equation}
\begin{cases}
\partial_n^\uparrow u(p_5)+(\frac{15}{7})^{m_1}\big(2u(p_5)-u(p_4)-f(q_0)\big)=0,\\
\partial_n^\uparrow u(p_4)+(\frac{15}{7})^{m_1}\big(2u(p_4)-u(p_5)-f(q_0)\big)=0.
\end{cases}
\end{equation}

Case 2($\iota_1=2$):
\begin{equation}
\begin{cases}
4u(p_5)-u(p_1)-u(p_3)-u(p_4)-f(q_0)=0,\\
4u(p_4)-u(p_2)-u(p_3)-u(p_5)-f(q_0)=0,\\
\partial_n^\uparrow u(p_1)+(\frac{15}{7})^{m_1}\big(2u(p_1)-u(p_3)-u(p_5)\big)=0,\\
\partial_n^\uparrow u(p_2)+(\frac{15}{7})^{m_1}\big(2u(p_2)-u(p_3)-u(p_4)\big)=0,\\
\partial_n^\uparrow u(p_3)+(\frac{15}{7})^{m_1}\big(4u(p_3)-u(p_1)-u(p_2)-u(p_4)-u(p_5)\big)=0.
\end{cases}
\end{equation}

Due to the same consideration in Section 2, we need to express the normal derivatives $\partial_n^\uparrow u(p_i)$'s in the above equations in terms of $u(p_i)$'s and the boundary data $f$. Thus we turn to find the explicit representation of  $\partial^\uparrow_n u(q_0)$ in terms of the boundary values. 

For this purpose, we need to look at the normal derivative along the boundary of a special harmonic function on $\Omega_\lambda$, denoted by $h^\lambda_0$,   assuming value $1$ at $q_0$ and $0$ along $X$. We will write $h_0=h_0^\lambda$ for simplicity if there is no confusion.

\textbf{Lemma 3.1.} \textit{Denote $\alpha(\lambda)=h_{0}(p_4)$($=h_{0}(p_5)$ by symmetry) and $\eta(\lambda)=\partial^\uparrow_n h_{0}(q_0)=2(\frac{15}{7})^{m_1}\big(1-\alpha(\lambda)\big)$. Then}
\begin{equation}
\alpha(\lambda)=
\begin{cases}
\frac{1}{1+\eta(R\lambda)}\qquad\qquad\qquad \text{ if } \iota_1=1,\\
\frac{6+6\eta(R\lambda)+\eta(R\lambda)^2}{6+15\eta(R\lambda)+3\eta(R\lambda)^2}\qquad \text{if } \iota_1=2.
\end{cases}
\end{equation}
\textit{Denote} 
\[T^\lambda (x)=
\begin{cases}
\frac{1}{1+2(\frac{15}{7})^{m_2-m_1}(1-x)}\qquad\qquad\qquad\qquad\qquad \text{ if } \iota_1=1,\\
\frac{3+6(\frac{15}{7})^{m_2-m_1}(1-x)+2\big((\frac{15}{7})^{m_2-m_1}(1-x)\big)^2}{3+15(\frac{15}{7})^{m_2-m_1}(1-x)+6\big((\frac{15}{7})^{m_2-m_1}(1-x)\big)^2}\qquad \text{if } \iota_1=2, \end{cases}\]
\textit{then $\alpha(\lambda)=T^\lambda\big(\alpha(R\lambda)\big)$. In addition, }
\[\alpha(\lambda)=\lim_{n\to\infty} T^{\lambda}\circ T^{\lambda_1}\circ\cdots\circ T^{\lambda_n}(0).\]

\textit{Proof.} By (3.2) and the definition of $h_0$, we have
$h_0\circ F^\lambda_i=h_0(p_i)h^{R\lambda}_0$, $\forall i\in W^\lambda_1.$ Taking  $u=h_0$ and $\partial_n^\uparrow h_0(p_i)=(\frac{15}{7})^{m_1}h_0(p_i)\eta(R\lambda)$ into (3.3) or (3.4), noticing that $h_0$ is symmetric in $\Omega_\lambda$, we have,
\[\eta(R\lambda)\alpha(\lambda)+\big(2\alpha(\lambda)-\alpha(\lambda)-1\big)=0,\text{ if }\iota_1=1,\]
or
\[
\begin{cases}
4\alpha(\lambda)-\alpha(\lambda)-1-h_0(p_1)-h_0(p_3)=0,\\
\eta(R\lambda)h_0(p_1)+\big(2h_0(p_1)-\alpha(\lambda)-h_0(p_3)\big)=0,\\
\eta(R\lambda)h_0(p_3)+\big(4u(p_3)-2\alpha(\lambda)-2h_0(p_1)\big)=0,
\end{cases}\text{ if }\iota_1=2.
\]

Solving the above equations, we get (3.5) and  
\begin{equation}
\begin{cases}
h_0(p_1)=h_0(p_2)=\frac{6+\eta(R\lambda)}{6+15\eta(R\lambda)+3\eta(R\lambda)^2},\\
h_0(p_3)=\frac{6+2\eta(R\lambda)}{6+15\eta(R\lambda)+3\eta(R\lambda)^2},
\end{cases}
\text{ in case of } \iota_1=2.
\end{equation}
Moreover, substituting $\eta(R\lambda)=2(\frac{15}{7})^{m_2-m_1}\big(1-\alpha(R\lambda)\big)$ into (3.5), we have $\alpha(\lambda)=T^{\lambda}\big(\alpha(R\lambda)\big)$. Inductively,
\[\alpha(\lambda)=T^{\lambda}\circ T^{\lambda_1}\circ\cdots\circ T^{\lambda_{n-1}}\big(\alpha(\lambda_n)\big).\]

For the rest part of the theorem, we introduce a sequence of functions $u^\lambda_n$ which assume values 
$u^\lambda_n(F^\lambda_w p^{\lambda_n}_4)=u^\lambda_n(F^\lambda_wp^{\lambda_n}_5)=0,\forall w\in W^\lambda_n$ and $u^\lambda_n|_{\partial\Omega_\lambda}=h_0|_{\partial\Omega_\lambda}$, and take harmonic extension in the remaining region. A similar discussion yields that
\[
u^\lambda_n(p_4)=T^{\lambda}\circ T^{\lambda_1}\circ\cdots\circ T^{\lambda_{n-1}}(0).
\]
In fact,  $\forall i\in W^\lambda_1$, we have $u^\lambda_n\circ F^\lambda_i=u^\lambda_n(p_i)u^{R\lambda}_{n-1}$. Taking $u=u_n^\lambda$ and $\partial_n^\uparrow u^\lambda_n(p_i)=(\frac{15}{7})^{m_1}u^\lambda_n(p_i)\partial_n^\uparrow u^{R\lambda}_{n-1}(q_0)$ into (3.3) or (3.4), after a similar calculation we have $u^\lambda_n(p_4)=T^{\lambda}\big(u_{n-1}^{\lambda_1}(p_4^{\lambda_1})\big)$.

Looking at the region bounded by $\{q_0\}\cup\{F^{\lambda}_{w}p_4^{\lambda_n}, F^{\lambda}_{w}p_5^{\lambda_n}\}_{w\in W^\lambda_n}$, applying the maximum principle for harmonic functions, we get
\[
0\leq(h_0-u^\lambda_n)(p_4)\leq \max\limits_{w\in W^\lambda_n,i=4,5}(h_0-u^\lambda_n)(F^{\lambda}_{w}p_i^{\lambda_n}).
\]
Since for $i=4,5$, $u^\lambda_n(F^{\lambda}_{w}p_i^{\lambda_n})=0, \forall w\in W_n^\lambda$, and $h_0(F^{\lambda}_{w}p^{\lambda_n}_i)$ converges uniformly to $0$ as $|w|\to\infty$, we then have $u^\lambda_n(p_4)\to h_0(p_4)$ as $n\to\infty$.\hfill$\square$

\textbf{Remark.} In fact, we have 
\[\alpha(\lambda)=\lim_{n\to\infty} T^\lambda\circ T^{\lambda_1}\circ\cdots\circ T^{\lambda_n}(c),\]
for any fixed constant $-\infty<c<1$. We only need small changes in the proof, and readers may refer to Theorem 4.7 for a similar discussion.

$\alpha$ is an increasing function of $\lambda$ on $(\frac{1}{3},1]$, because $h_0^{\lambda_b}|_{\bar{\Omega}_{\lambda_a}}\geq h_0^{\lambda_a}$ if $\lambda_b\geq \lambda_a$. For $0<\lambda\leq\frac{1}{3}$, we have $\alpha(\lambda)=\alpha(3^{m_1-1}\lambda)$ by dilation. Thus
\begin{equation} 
0<\alpha(\lambda)\leq \alpha(1)=\frac{75 - \sqrt{2353}}{60}\approx 0.441538,
\end{equation}
where $\alpha(1)$ is the root of $x=\frac{3+6\cdot\frac{15}{7}(1-x)+2\cdot(\frac{15}{7})^2(1-x)^2}{3+15\cdot\frac{15}{7}(1-x)+6\cdot(\frac{15}{7})^2(1-x)^2}$.\\

\textbf{Definition 3.2.}\textit{ Let $0<\lambda\leq 1$, denote}
\[
\mu^\lambda_i=\frac{1}{2} \textit{\text{ for }} i=4,5,\textit{\text{ if }}\iota_1=1,
\]
\textit{or}
\[
\mu^\lambda_i=
\begin{cases}
\frac{6+\eta(R\lambda)}{18+4\eta(R\lambda)}\textit{\text{ for }} i=1,2,\\
\frac{3+\eta(R\lambda)}{9+2\eta(R\lambda)}\textit{\text { for }} i=3,
\end{cases}\textit{\text{ if }}  \iota_1=2.
\]
\textit{Define a probability measure $\mu^\lambda$ on $X$ by}
\begin{equation}
\mu^\lambda(X_w)=\prod_{k=1}^{|w|}\mu^{\lambda_{k-1}}_{w_k},\forall w\in W^\lambda_*.
\end{equation}
We can easily verify that 
\begin{equation}
\mu^\lambda\circ F^\lambda_w=(\prod_{k=1}^{n}\mu^{\lambda_{k-1}}_{w_k})\mu^{\lambda_n}, \forall w\in W^\lambda_n,
\end{equation}
and
\[\frac{1}{4}<\mu^\lambda_1=\mu^\lambda_2<\frac{1}{3},\text{ }  \frac{1}{3}<\mu^\lambda_3<\frac{1}{2},\text{ if }\iota_1(\lambda)=2,\]
since $\eta(R\lambda)>0$. 

\textbf{Lemma 3.3.}  $\forall w\in W_*^\lambda$, $\partial_n^\uparrow h_0(F_w^\lambda q_0)=\mu^\lambda(X_w)\eta(\lambda)$.

\textit{Proof.} By the local Gauss-Green's formula, we have
$\partial^\uparrow_n h_0(q_0)=\sum_{_{w\in W^\lambda_1}}\partial^\uparrow_nh_0(F^\lambda_wq_0).$

If $\iota_1=1$, then $\partial^\uparrow_nh_0(p_i)=\frac{1}{2}\eta(\lambda)=\mu^\lambda_i\eta(\lambda)$ by symmetry consideration. If $\iota_1=2$, then $\partial^{\uparrow}_n h_0(p_1):\partial^{\uparrow}_n h_0(p_2):\partial^{\uparrow}_n h_0(p_3)=h_0(p_1):h_0(p_2):h_0(p_3)$, which still leads to $\partial_n^\uparrow h_0(p_i)=\mu^\lambda_i\eta(\lambda)$ by (3.6). Applying the above discussion iteratively, we get 
\[
\begin{aligned}
\partial_n^\uparrow h_0(F^\lambda_wq_0)&=(\mu^{\lambda_{|w|-1}}_{w_{|w|}})\partial_n^\uparrow h_0(F^\lambda_{w_1w_2\cdots w_{|w|-1}}q_0)=\cdots\\
&=\eta(\lambda)(\prod_{k=1}^{|w|} \mu^{\lambda_{k-1}}_{w_k})=\mu^\lambda(X_w)\eta(\lambda), \forall w\in W^\lambda_*. \qquad \hfill\square
\end{aligned}
\]

\textbf{Theorem 3.4.}\textit{ Let $u$ be a solution of the Dirichlet problem (1.1). Then }
\begin{equation}
\partial^\uparrow_n u(q_0)=\eta(\lambda)\big(f(q_0)-\int_{X} fd\mu^{\lambda}\big).
\end{equation}
\textit{In addition, if $u\in dom\mathcal{E}_{\Omega_{\lambda}}$, we have}
\[\mathcal{E}_{\Omega_\lambda}(h_0,u)=\partial^\uparrow_n u(q_0)=\eta(\lambda)\big(f(q_0)-\int_{X} fd\mu^{\lambda}\big).\]

\textit{Proof.} Consider the simple set $O_{\lambda,n}$ whose boundary vertices are $\{q_0\}\bigcup\{F^\lambda_wq_0\}_{w\in W^\lambda_n}$. By the local Gauss-Green's formula, we have
\begin{equation}
\mathcal{E}_{O_{\lambda,n}}(h_0,u)=\partial^\uparrow_n h_0(q_0)u(q_0)+\sum_{_{w\in W^\lambda_n}}(-\partial^\uparrow_nh_0)(F^\lambda_wq_0)u(F^\lambda_wq_0).
\end{equation}

Applying Lemma 3.3 and comparing to the right side of (3.10), we get
\[
\begin{aligned}
\mathcal{E}_{O_{\lambda,n}}(h_0,u)-\eta(\lambda)\big(f(q_0)-\int_{X} fd\mu^{\lambda}\big)&=\eta(\lambda)\big(\int_{X} fd\mu^{\lambda}-\sum_{w\in W^\lambda_n}u(F^\lambda_wq_0)\mu^\lambda(X_w)\big)\\
&=\eta(\lambda)\int_{X} \big(f-\sum_{w\in W^\lambda_n}u(F^\lambda_wq_0)\chi_{_{X_w}}\big)d\mu^{\lambda}\\
&\to 0, \text{ as } n\to\infty,
\end{aligned}
\]
where $\chi_{_{X_w}}$ is the characteristic function of $X_w$. So for any $u\in dom\mathcal{E}_{\Omega_\lambda},$ we have
\[\mathcal{E}_{\Omega_\lambda}(h_0,u)=\eta(\lambda)\big(f(q_0)-\int_X fd\mu^\lambda \big).\]

To show the rest part of the theorem, define a sequence of harmonic functions $u_k$, which are piecewise constant on $X$, taking boundary values $u_k(q_0)=f(q_0)$ and $u_k|_{X_\tau}=u(F^\lambda_\tau q_0),\forall \tau\in W^\lambda_k.$ Applying the Gauss-Green's formula, we get that
\[
\begin{aligned}
\mathcal{E}_{O_{\lambda,n}}(h_0,u_k)&=\partial^\uparrow_n u_k(q_0)h_0(q_0)+\sum_{_{w\in W^\lambda_n}}(-\partial^\uparrow_nu_k)(F^\lambda_wq_0)h_0(F^\lambda_wq_0)\\
&=\partial^\uparrow_n u_k(q_0)+\sum_{\tau\in W^\lambda_k}\sum_{w\in W^{\lambda_k}_{n-k}}(-\partial^\uparrow_nu_k)(F^\lambda_{\tau w}q_0)h_0(F^\lambda_{\tau w}q_0)\\
&\to \partial^\uparrow_n u_k(q_0)\text{ as }n\to\infty,
\end{aligned}\]
where we use the fact that for each fixed $\tau\in W_k^\lambda$, $\sum_{w\in W^{\lambda_k}_{n-k}}|\partial^\uparrow_nu_k(F^\lambda_{\tau w}q_0)|$ is uniformly bounded and all $h_0(F^\lambda_{\tau w}q_0)$'s converge uniformly to $0$ as $n$ goes to infinity.
Noticing that $\lim_{n\to\infty} \mathcal{E}_{O_{\lambda,n}}(h_0,u_k)=\eta(\lambda)\big(u_k(q_0)-\int_{X} u_kd\mu^{\lambda}\big)$, we get
\[\partial_n^\uparrow u_k(q_0)=\eta(\lambda)\big(f(q_0)-\int_X u_k d\mu^\lambda\big).\]

Finally, letting $k\to\infty$, noticing that $u_k$ converges to $u$ uniformly in $\Omega_\lambda$, we obtain (3.10).\hfill$\square$
 
\textbf{Corollary 3.5.} 
$\mathcal{E}_{\Omega_\lambda} (h_0)=\eta(\lambda).$

\textit{Proof.} In the proof of Theorem 3.4, we have seen that $\lim_{n\to\infty}\mathcal{E}_{O_{\lambda,n}}(h_0,h_0)=\partial_n^\uparrow h_0(q_0)h_0(q_0)=\eta(\lambda)$.\hfill$\square$

By (3.7) and the fact that $\eta(\lambda)=2(\frac{15}{7})^{m_1}\big(1-\alpha(\lambda)\big)$, the energy of $h_0$ is estimated by
\[2\big(1-\alpha(1)\big)(\frac{15}{7})^{m_1}\approx 1.116924(\frac{15}{7})^{m_1}\leq\eta(\lambda)=\mathcal{E}_{\Omega_\lambda} (h_0)<2(\frac{15}{7})^{m_1}.\]
This result will be helpful in the energy estimate for general functions.\\

Combining Theorem 3.4  with equations (3.3) and (3.4), by solving linear equations, we could calculate the values of the solution $u$ at the ``crucial'' points $p_i$'s for $i\in W_1^{\lambda}$, which are sufficient to recover $u$ by iteration. 

\textbf{Theorem 3.6.(Extension Algorithm)}\textit{ There exists a unique solution of the Dirichlet problem (1.1). In addition, we have the following formulas for $u(p_i),i\in W^\lambda_1$. }

\textit{Case 1 ($\iota_1=1$):}
\[
\begin{aligned}
u(p_4)=&\frac{1}{1+\eta(R\lambda)}f(q_0)+\frac{2\eta(R\lambda)+\eta(R\lambda)^2}{3+4\eta(R\lambda)+\eta(R\lambda)^2}\int_{X_{R\lambda}} f\circ F^\lambda_4d\mu^{R\lambda}\\
&+\frac{\eta(R\lambda)}{3+4\eta(R\lambda)+\eta(R\lambda)^2}\int_{X_{R\lambda}} f\circ F^\lambda_5d\mu^{R\lambda}.
\end{aligned}
\]

\textit{Case 2 ($\iota_1=2$):}

\[
\begin{aligned}
u(p_1)=&\frac{1}{54+165\eta(R\lambda)+102\eta(R\lambda)^2+15\eta(R\lambda)^3}\Big(\big(54+39\eta(R\lambda)+5\eta(R\lambda)^2\big)f(q_0)\\
&+\big(60\eta(R\lambda)+76\eta(R\lambda)^2+15\eta(R\lambda)^3\big)\int_{X_{R\lambda}} f\circ F^\lambda_1d\mu^{R\lambda}+\\
&(30\eta(R\lambda)+\eta(R\lambda)^2)\int_{X_{R\lambda}} f\circ F^\lambda_2d\mu^{R\lambda}+\big(36\eta(R\lambda)+20\eta(R\lambda)^2\big)\int_{X_{R\lambda}} f\circ F^\lambda_3d\mu^{R\lambda}\Big),
\end{aligned}\]

\[\begin{aligned}
u(p_3)=&\frac{6+2\eta(R\lambda)}{6+15\eta(R\lambda)+3\eta(R\lambda)^2}f(q_0)+\frac{5\eta(R\lambda)+3\eta(R\lambda)^2}{6+15\eta(R\lambda)+3\eta(R\lambda)^2}\int_{X_{R\lambda}} f\circ F^\lambda_3d\mu^{R\lambda}
\\&+\frac{4\eta(R\lambda)}{6+15\eta(R\lambda)+3\eta(R\lambda)^2}\big(\int_{X_{R\lambda}} f\circ F^\lambda_1d\mu^{R\lambda}+\int_{X_{R\lambda}} f\circ F^\lambda_2d\mu^{R\lambda}\big),
\end{aligned}\]

\[\begin{aligned}
u(p_4)=&\frac{1}{54+165\eta(R\lambda)+102\eta(R\lambda)^2+15\eta(R\lambda)^3}\Big(\big(54+84\eta(R\lambda)+39\eta(R\lambda)^2+5\eta(R\lambda)^3\big)f(q_0)\\+&\big(24\eta(R\lambda)+12\eta(R\lambda)^2+\eta(R\lambda)^3\big)\int_{X_{R\lambda}} f\circ F^\lambda_1d\mu^{R\lambda}
\\+&\big(30\eta(R\lambda)+27\eta(R\lambda)^2+4\eta(R\lambda)^3\big)\int_{X_{R\lambda}} f\circ F^\lambda_2d\mu^{R\lambda}
\\+&\big(27\eta(R\lambda)+24\eta(R\lambda)^2+5\eta(R\lambda)^3\big)\int_{X_{R\lambda}} f\circ F^\lambda_3d\mu^{R\lambda}
\Big).
\end{aligned}\]
\textit{The formulas for $u(p_2),u(p_5)$ can be obtained symmetrically.}

\textit{Proof. } See Proposition 1.1 for the existence and uniqueness of the solution. Substituting $\partial_n^\uparrow u(p_i)=(\frac{15}{7})^{m_1} \eta(R\lambda) \big(u(p_i)-\int_{X_{R\lambda}} f\circ F^\lambda_i d\mu^{R\lambda}\big)\text{ for }i\in W^\lambda_1$ into (3.3) or (3.4), after solving linear equations, we get the result.\hfill$\square$

\subsection{Haar series expansion and energy estimate}

Now we consider the energy estimate for the harmonic solutions in terms of their boundary values. 
For a harmonic function $u$ with boundary value $u|_X=f$ in $L^2(X,\mu^\lambda)$, we will give an estimation of $\mathcal{E}_{\Omega_\lambda}(u)$ in terms of the Fourier coefficients of $f$ with respect to a Haar basis. 

\textbf{Definition 3.7.} (1) \textit{Assume $0<\lambda\leq 1$.
Define $\psi^{(1),\lambda}$ and $\psi^{(2),\lambda}$ to be piecewise constant functions on $X_\lambda$ such that},
\begin{equation}
\psi^{(1),\lambda}|_{X_5}=1, \psi^{(1),\lambda}|_{X_4}=-1, \textit{\text{ if }}  \iota_1=1,
\end{equation}

\begin{equation}
\begin{cases}
 \psi^{(1),\lambda}|_{X_1}=1, \psi^{(1),\lambda}|_{X_2}=-1,\psi^{(1), \lambda}|_{X_3}=0,\\ \psi^{(2),\lambda}|_{X_1}=\psi^{(2), \lambda}|_{X_2}=\mu^\lambda_3, \psi^{(2),\lambda}|_{X_3}=-2\mu^\lambda_1,
\end{cases}\textit{\text{ if }}  \iota_1=2,
\end{equation}
\textit{ and there is no $\psi^{(2),\lambda}$ in case of $\iota_1=1$.
In addition, for $w\in W^\lambda_n$, define $\psi_w^{(j),\lambda}$ supported in $X_w$ by}
\[
\psi^{(1),\lambda}_w=\psi^{(1),\lambda_n}\circ (F^\lambda_w)^{-1},
\]
\textit{and similarly}
\[\psi^{(2),\lambda}_w=\psi^{(2),\lambda_n}\circ (F^\lambda_w)^{-1}, \textit{\text{ if }} \iota_{n+1}=2.\]

(2) \textit{For $w\in W^\lambda_*$ and $j\leq \iota_{|w|+1}$, define a series of harmonic functions on $\Omega_\lambda$ by
\[h^{(j),\lambda}_w(q_0)=0,\textit{ and } h^{(j),\lambda}_w|_X=\psi^{(j),\lambda}_{w}.\]
}

We will write $\psi^{(j)}_w$ and $h^{(j)}_w$ instead of $\psi^{(j),\lambda}_w$ and $h^{(j),\lambda}_w$ when it causes no confusion. Set $\psi_\emptyset ^{(j),\lambda}=\psi ^{(j),\lambda}$, and  $h^{(j),\lambda}_\emptyset=h^{(j),\lambda}$ for the sake of formality.

It is easy to check that for $w\in W^\lambda_*$ and $j\leq \iota_{|w|+1}$, $h^{(j)}_w$ is supported in $F_w^\lambda\bar{\Omega}_{\lambda_n}$ with $h^{(j)}_w=h^{(j),\lambda_n}_\emptyset\circ (F^\lambda_w)^{-1}$, noticing that $\partial_n^\uparrow h^{(j)}_\emptyset(q_0)=0$ by Theorem 3.4. Moreover,  $\{\psi_w^{(j)}\}_{w\in W^\lambda_*,j\leq \iota_{|w|+1}}\cup\{1\}$ form an orthogonal basis of $L^2(X,\mu^\lambda)$.  For convenience of calculation, we do not normalize these basis functions. See Figure 3.4 for an illustration of $h^{(1)}$ and $h^{(2)}$. 
\begin{figure}[h]
\begin{center}
\includegraphics[width=3.5cm]{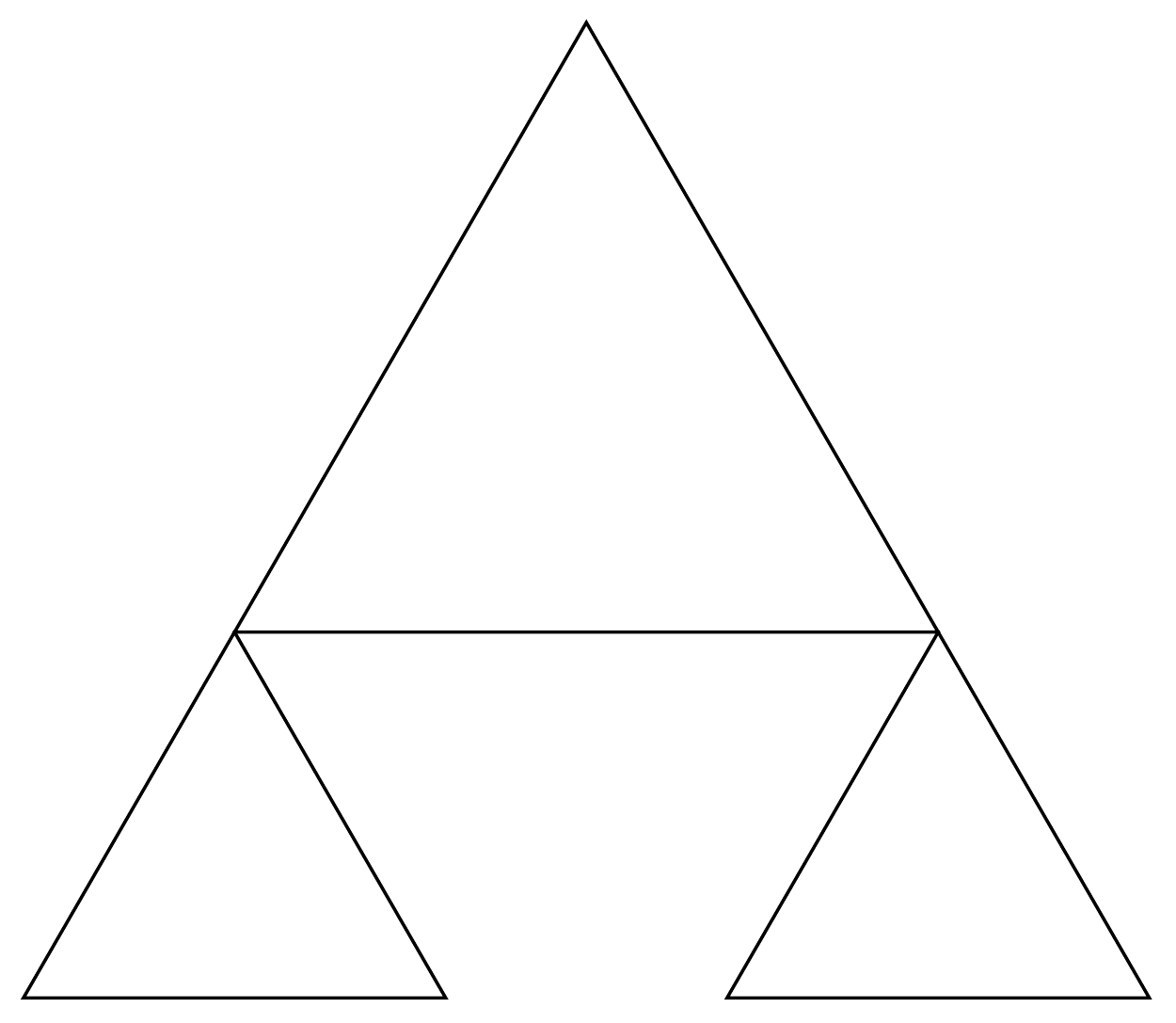}\hspace{0.3cm}
\includegraphics[width=3.5cm]{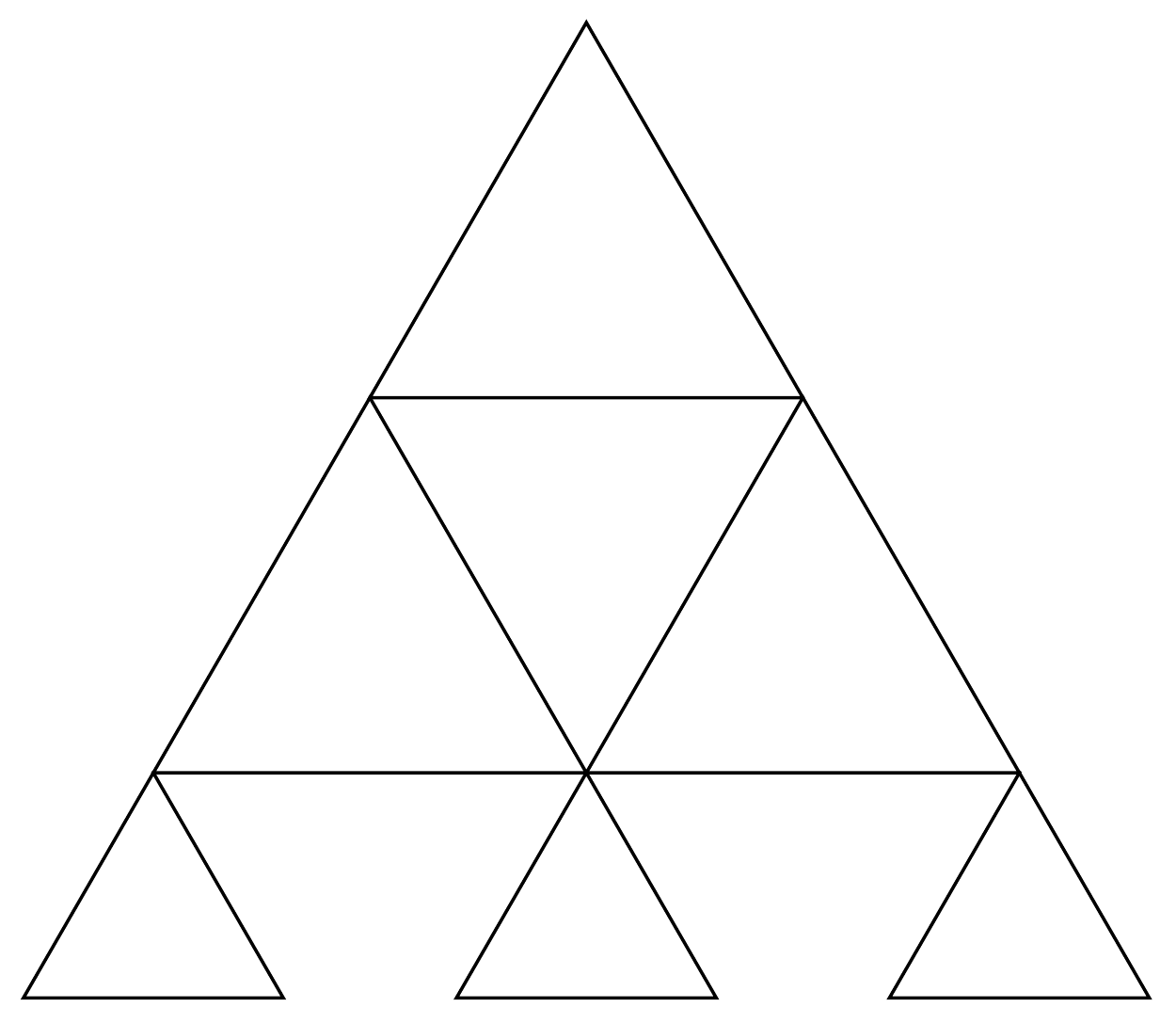}\hspace{0.3cm}
\includegraphics[width=3.5cm]{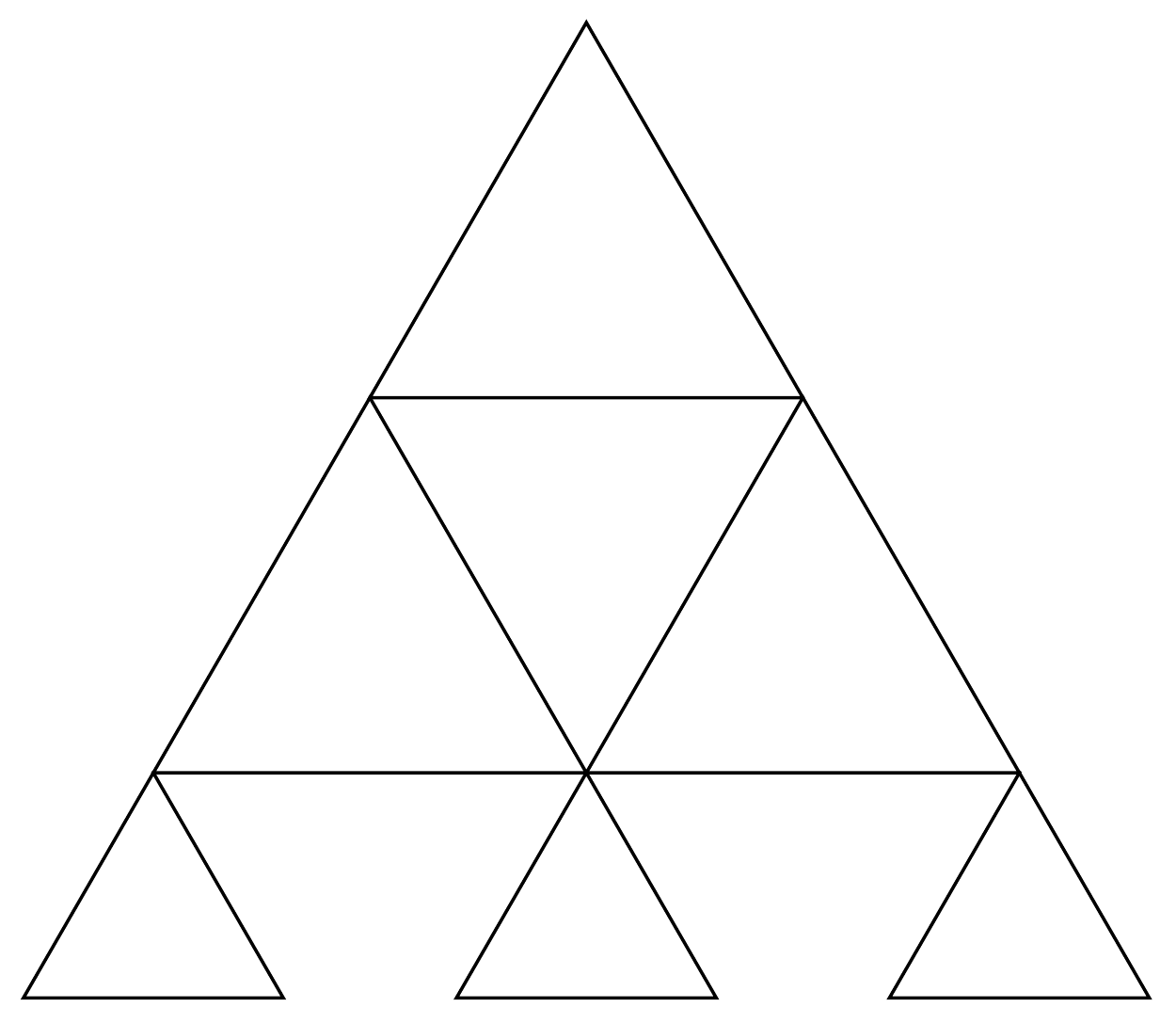}
\setlength{\unitlength}{1cm}
\begin{picture}(0,0) \thicklines
\put(-9.9,3.05){$0$}
\put(-10.95,-0.23){$1$}
\put(-9.1,-0.23){$-1$}

\put(-5.98,3.05){$0$}
\put(-7.25,-0.23){$1$}
\put(-4.89,-0.23){$-1$}
\put(-6.0,-0.23){$0$}

\put(-2.07,3.05){$0$}
\put(-2.80,1.90){$0$}
\put(-1.37,1.90){$0$}
\put(-2.50,-0.23){$-2\mu^\lambda_1$}
\put(-3.45,-0.23){$\mu^\lambda_3$}
\put(-0.90,-0.23){$\mu^\lambda_3$}

\put(-2.9,-0.6){\textbf{$h^{(2)}(\iota_1=2)$}}
\put(-10.7,-0.6){\textbf{$h^{(1)}(\iota_1=1)$}}
\put(-6.8,-0.6){\textbf{$h^{(1)}(\iota_1=2)$}}
\end{picture}
\vspace{0.6cm}
\begin{center}
\textbf{Figure 3.4. Boundary values of $h^{(1)}$, $h^{(2)}$. }
\end{center}
\end{center}
\end{figure}

Before performing the energy estimate, we list two basic lemmas. 

\textbf{Lemma 3.8.}
\textit{There exist two positive constants $C_1$ and $C_2$, such that}
\[C_1(\frac{15}{7})^{m_{|w|+1}}\leq \mathcal{E}_{\Omega_\lambda}(h^{(j)}_w)\leq C_2(\frac{15}{7})^{m_{|w|+1}},
\]
\textit{for all $0<\lambda\leq 1$ and $w\in W^\lambda_*, j\leq \iota_{|w|+1}$.}

\textit{Proof.} We only need to prove that $\mathcal{E}_{\Omega_\lambda}(h^{(j)})$ is bounded above and below by multiples of $(\frac{15}{7})^{m_1}$, as $
h^{(j),\lambda}_w=h^{(j),\lambda_n}_\emptyset\circ (F^\lambda_w)^{-1}$. In particular, we will restrict our consideration to $\frac{1}{3}<\lambda\leq 1$, since for general $0<\lambda\leq 1$, $h^{(j),\lambda}$ is supported in $F_0^{m_1-1}(\mathcal{SG}_3)$ with $h^{(j),\lambda}=h^{(j),3^{m_1-1}\lambda}\circ F_0^{-m_1+1}$. So it is sufficient to  assume $m_1=1$ and prove that $\mathcal{E}_{\Omega_\lambda}(h^{(j)})$ is bounded above and below by two positive constants.

First, we consider $\frac{2}{3}<\lambda\leq 1$. In this case, $\iota_1=2$. Let $c_1, c_2, c_3$ be some selected constants independent of $\lambda$. For each $\frac{2}{3}<\lambda\leq 1$, write $v^\lambda$ the harmonic function on $\bar{\Omega}_\lambda$ which assumes $0$ at $q_0$ and takes constant $c_i$ along $X_i$ for $i=1,2,3$. 
We claim that $$\mathcal{E}_{\Omega_{\lambda_a}}(v^{\lambda_a})\geq \mathcal{E}_{\Omega_{\lambda_b}}(v^{\lambda_b}),\text{ if }\lambda_a\leq \lambda_b.$$

To prove the claim, we construct another function $\tilde{v}$ on $ \bar{\Omega}_{\lambda_b}$ such that $\tilde{v}|_{\partial \Omega_{\lambda_b}}=v^{\lambda_b}|_{\partial \Omega_{\lambda_b}}, \tilde{v}(p_i)=v^{\lambda_a}(p_i)$ for each $1\leq i\leq 5$ and $\tilde{v}$ is harmonic in remaining region. So $\tilde{v}|_A=v^{\lambda_a}|_A$, for $A=F_0\mathcal{SG}_3\cup F_4\mathcal{SG}_3\cup F_5\mathcal{SG}_3$, which says that

\[\begin{aligned}
\mathcal{E}_{\Omega_{\lambda_b}}(\tilde{v})\leq\mathcal{E}_{\Omega_{\lambda_a}}(v^{\lambda_a}),
\end{aligned}\]
by using Corollary 3.5 and the fact that $\eta(R\lambda)$ is decreasing on $\frac{2}{3}<\lambda\leq 1$. In addition, we have $\mathcal{E}_{\Omega_{\lambda_b}}(\tilde{v})\geq \mathcal{E}_{\Omega_{\lambda_b}}(v^{\lambda_b})$, since harmonic functions minimize the energy. Combining the two inequalities, we obtain the claim. 

Thus, $\mathcal{E}_{\Omega_{\lambda}}(v^{\lambda})\geq \mathcal{E}_{\Omega_{1}}(v^1)$ \big(more precisely, $\mathcal{E}_{\Omega_{\lambda}}(v^{\lambda})\geq\lim_{\epsilon\rightarrow 0}\mathcal{E}_{\Omega_{1-\epsilon}}(v^{1-\epsilon})$\big) which provides a lower bound of $\{\mathcal{E}_{\Omega_{\lambda}}(v^{\lambda})\}$. On the other hand, to find an upper bound of $\{\mathcal{E}_{\Omega_\lambda}(v^\lambda)\}$, consider the function $\bar{v}\in C(A)$ which is harmonic in $A$ and assumes $0$ at $q_0$, $c_i$ at $p_i$ for $i=1,2,3$. It is easy to find that $\mathcal{E}_{\Omega_{\lambda}}(v^\lambda)\leq \mathcal{E}_A(\bar{v})$ by extending $\bar{v}$ to $\Omega_\lambda$ with $\bar{v}|_{F_i^\lambda\Omega_{R\lambda}}=c_i$.

The energy estimate of $h^{(1)}$ is a special case in the above discussion. To estimate the energy of $h^{(2)}$, observe that the boundary values $(h^{(2)}|_{X_1},h^{(2)}|_{X_2},h^{(2)}|_{X_3})$ vary within a compact set, denoted by $C$, since we always have $\frac{1}{4}<\mu_1^\lambda=\mu_2^\lambda<\frac{1}{3}$ and $\frac{1}{3}<\mu_3^\lambda<\frac{1}{2}$. We then have 
\[ \begin{aligned}  0<\inf\{\mathcal{E}_{\Omega_1}(v^1):( v^1|_{X_{1,1}}&,v^1|_{X_{1,2}},v^1|_{X_{1,3}})\in C\}
\\\leq \mathcal{E}_{\Omega_\lambda}(h^{(2)})
&\leq \sup\{\mathcal{E}_{A}(\bar{v}):\big(\bar{v}(p_1),\bar{v}(p_2),\bar{v}(p_3)\big)\in C\}< \infty.
\end{aligned}\]

Thus we have proved that $\mathcal{E}_{\Omega_\lambda}(h^{(j)})$ is bounded above and below by two positive constants, when $\frac{2}{3}<\lambda\leq 1$. 

A similar discussion is valid for $\frac{1}{3}<\lambda\leq\frac{2}{3}$. 
\hfill$\square$\\

The next lemma shows that all the basis functions $\{h_w^{(j)}\}$ are pairwise orthogonal in energy. 

\textbf{Lemma 3.9.} \textit{Assume $0<\lambda\leq 1$ and $w,w'\in W^\lambda_*$. Then $\mathcal{E}_{\Omega_\lambda}(h^{(j)}_w,h_{w'}^{(j')})\neq 0$ if and only if $w=w'$ and $j=j'$. }

\textit{Proof.} We discuss in different cases.

If $X_w\bigcap X_{w'}=\emptyset$, then obviously $\mathcal{E}_{\Omega_\lambda}(h^{(j)}_w,h_{w'}^{(j')})=0$, as $h_{w}^{(j)}$ and $h^{(j')}_{w'}$ support in disjoint regions.

If $X_w\subsetneq X_{w'}$, then
\[
\begin{aligned}
\mathcal{E}_{\Omega_\lambda}(h^{(j)}_w,h_{w'}^{(j')})=(\frac{15}{7})^{m_{|w|}}\big(\partial_n^\uparrow (h^{(j)}_w\circ F_{w}^\lambda)(q_0)\big)\big(h_{w'}^{(j')}(F_{w}^\lambda q_0)-h_{w'}^{(j')}|_{X_{w}}\big)=0,
\end{aligned}
\] 
by Theorem 3.4 and the fact that $\partial_n^\uparrow (h^{(j)}_w\circ F_{w}^\lambda)(q_0)=0$.

If $w=w'$ and $j\neq j'$, then $h^{(j)}_w$ and $h^{(j')}_{w'}$ assume different symmetries. So we also have $\mathcal{E}(h^{(j)}_w,h_{w'}^{(j')})=0$. \hfill$\square$\\

\textbf{Remark.} The proof of Lemma 3.9 also implies that $\mathcal{E}_{\Omega_\lambda}(h^{(j)}_w, h_0)= 0$ for each $w\in W_*^\lambda$, $j\leq \iota_{|w|+1}$.

Now for the  harmonic solution $u$ of the Dirichlet problem (1.1), we have the following estimate in energy in terms of its boundary data $f$.

\textbf{Theorem 3.10.} \textit{Let $u$ be the harmonic function in $\Omega_\lambda$ with boundary values $u(q_0)=a$ and $u|_{X_\lambda}=f$, where 
\[f=b+\sum_{w\in W^\lambda_*}\sum_{j\leq \iota_{|w|+1}} c^{(j)}_w\psi_w^{(j)}\]
with
\[
c_w^{(j)}=\big(\int_{X_\lambda} (\psi^{(j)}_w)^2 d\mu^\lambda\big)^{-1}\int_{X_\lambda} f\psi^{(j)}_wd\mu^\lambda.
\]
Then $\mathcal{E}_{\Omega_\lambda}(u)$ is bounded above and below by multiples of
\begin{equation}
(\frac{15}{7})^{m_1}(a-b)^2+\sum_{n=0}^{\infty}\sum_{w\in W^\lambda_n}\sum_{j=1}^{\iota_{n+1}}(\frac{15}{7})^{m_{n+1}}|c^{(j)}_w|^2.
\end{equation}
In particular, $u$ has finite energy if and only if (3.14) is finite.}

\textit{Proof.} We have

\[u=b+(a-b)h_0+\sum_{n=0}^{\infty}\sum_{w\in W^\lambda_n}\sum_{j=1}^{\iota_{n+1}}c^{(j)}_wh^{(j)}_w.\]
Since we have shown in Lemma 3.9 that the functions $h_0\bigcup \{h^{(j)}_w\}$ are orthogonal in energy,
\[\mathcal{E}_{\Omega_\lambda} (u)=(a-b)^2 \mathcal{E}_{\Omega_\lambda} (h_0)+\sum_{n=0}^{\infty}\sum_{w\in W^\lambda_n}\sum_{j=1}^{\iota_{n+1}} |c^{(j)}_w|^2\mathcal{E}_{\Omega_\lambda}(h^{(j)}_w).\]
Then $(3.14)$ follows from Lemma 3.8 and Corollary 3.5.\hfill$\square$\\

\textbf{Remark.} The $L^2$ norm of each Haar function $\psi_w^{(j)}$ is bounded above and below by multiples of $\big(\mu^\lambda(X_w)\big)^{\frac{1}{2}}$. If we normalize the functions $\psi_w^{(j)}$ into
\[
\bar{\psi}_w^{(j)}=\big(\int_X (\psi_w^{(j)})^2d\mu^\lambda\big)^{-\frac{1}{2}}\psi_w^{(j)},
\]
and write $\bar{c}_w^{(j)}=\int_X f\bar{\psi}^{(j)}_wd\mu^\lambda.$ Then, $\mathcal{E}_{\Omega_\lambda}(u)$ is bounded above and below by multiples of
\[
(\frac{15}{7})^{m_1}(a-b)^2+\sum_{n=0}^{\infty}\sum_{w\in W^\lambda_n}\sum_{j=1}^{\iota_{n+1}}(\frac{15}{7})^{m_{n+1}}\big(\mu^\lambda(X_w)\big)^{-1}|\bar{c}_w^{(j)}|^2.
\]

\section{
Dirichlet problems  on lower domains of $\mathcal{SG}$}

In this section, we consider the Dirichlet problem on lower domains of $\mathcal{SG}$. Similar to the last section, we assume $\mathcal{SG}$ is contained in $\mathbb{R}^2$ with boundary vertices $q_0=(\frac{1}{\sqrt{3}},1),\text{ }q_1=(0,0),\text{ }q_2=(\frac{2}{\sqrt{3}},0)$. Let $0\leq \lambda<1$, 
the lower domain of $\mathcal{SG}$, denoted by $\Omega_\lambda^-$, is defined as
\[\Omega^-_\lambda=\{(x,y)\in \mathcal{SG}\setminus V_0|y<1-\lambda\}.\]
The boundary of $\Omega_\lambda^-$ is
\[\partial\Omega_\lambda^-=X^-_\lambda\cup \{q_1,q_2\},\]
with 
\[X^-_\lambda=
\begin{cases}
V_{d(\lambda)}\cap X_\lambda, \text{ if $\lambda$ is a  dyadic rational,}\\
X_\lambda,\hspace{1.13cm} \text{ if $\lambda$ is not a dyadic rational,}
\end{cases}\]
where $d(\lambda)$ is the smallest integer such that $\lambda$ is a multiple of  $2^{-d(\lambda)}$ and $X_\lambda=\{(x,y)\in \mathcal{SG}|y=1-\lambda\}$. We still abbreviate $X_\lambda^-$ to $X$ throughout this section for convenience. Write $\bar{\Omega}_\lambda^{-}=\Omega_\lambda^-\cap \partial \Omega_\lambda^-$ the closure of $\Omega_\lambda^-$. See Figure 4.1 for two typical domains.

\begin{figure}[h]
\begin{center}
\includegraphics[width=4cm]{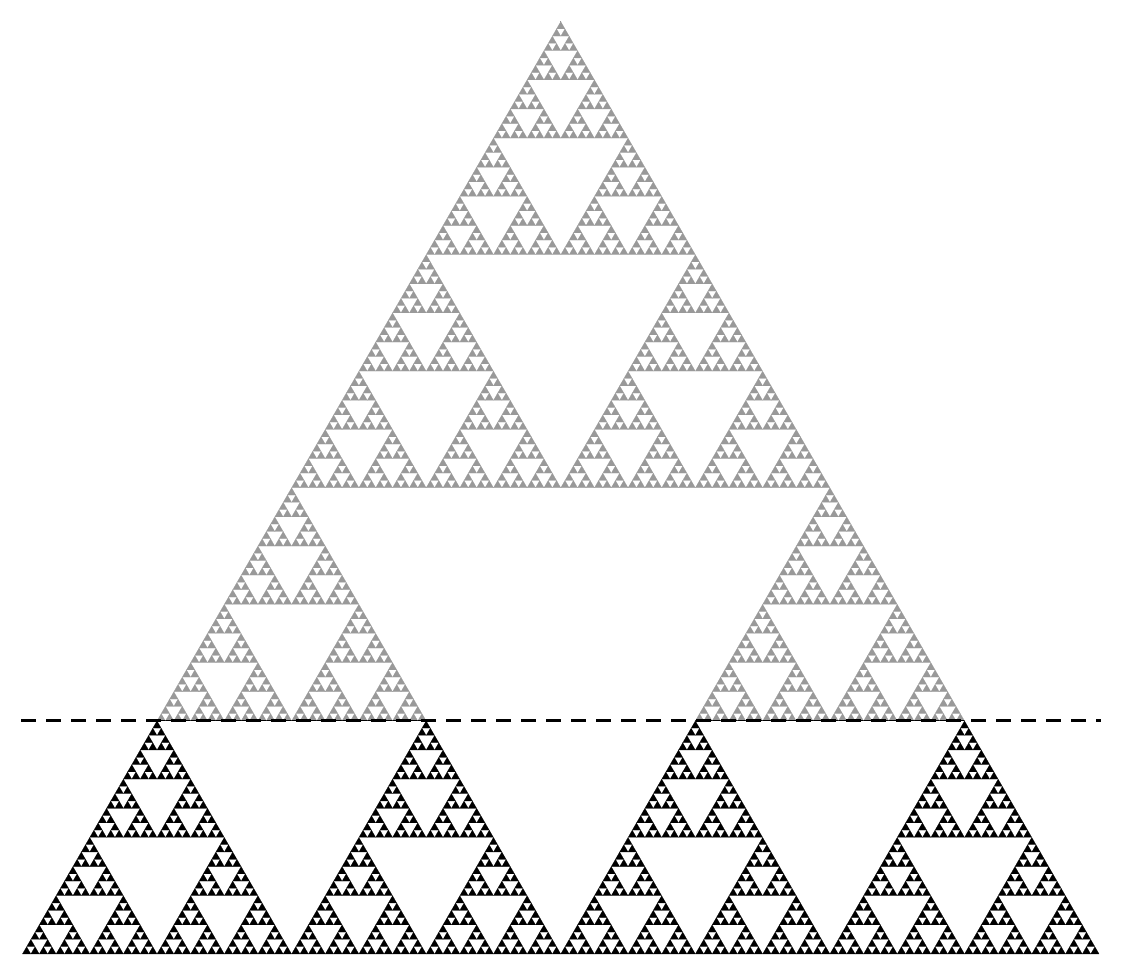}\hspace{1cm}
\includegraphics[width=4cm]{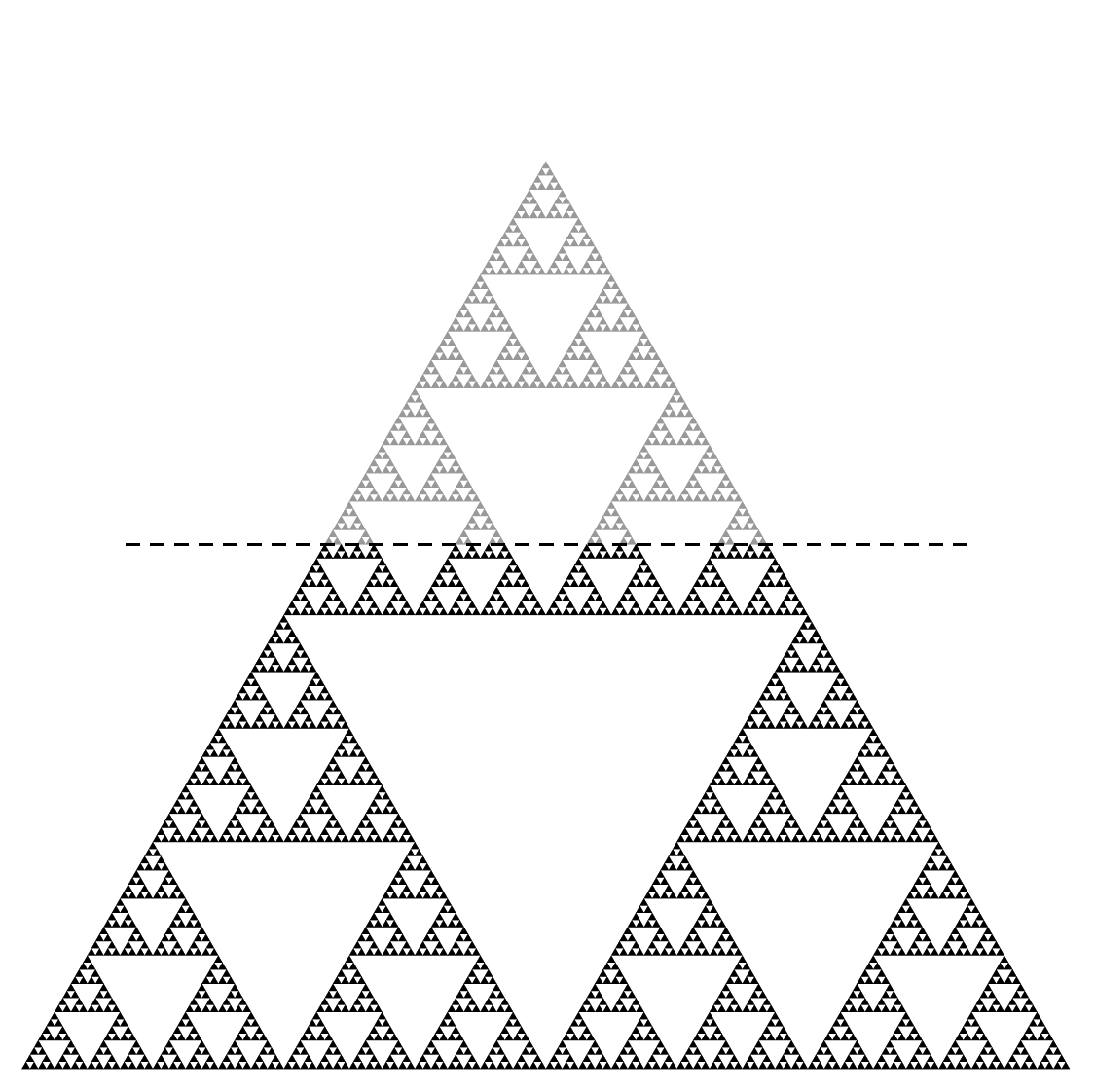}
\begin{center}
\textbf{Figure 4.1. Two typical domain $\Omega_\lambda^-$'s. $\lambda$ is (or not) a dyadic rational.}
\end{center}
\end{center}
\end{figure}

For $0\leq \lambda<1$, we write $\lambda$ in its binary expansion,
\[\lambda=\sum_{k=1}^{\infty} e_k(\lambda)2^{-k},e_k(\lambda)=0,1 \text{ for } k\geq 1.\]
We \textbf {forbid infinitely consecutive $1$'s} to make the expansion unique.  Denote
\[S\lambda=\sum_{k=1}^{\infty} e_{k+1}(\lambda)2^{-k}.\]
It is easy to check the relationship between  $\Omega_\lambda^-$ and $\Omega_{S\lambda}^-$ as following, 
\begin{equation}
\bar{\Omega}_\lambda^-=
\begin{cases}
F_0\bar{\Omega}_{S\lambda}^-\bigcup F_1\mathcal{SG}\bigcup F_2\mathcal{SG},\text{ if }e_1(\lambda)=0,\\
F_1\bar{\Omega}_{S\lambda}^-\bigcup F_2\bar{\Omega}_{S\lambda}^-,\hspace{1.1cm}\text{  if }e_1(\lambda)=1.
\end{cases}
\end{equation}
See Figure 4.2 for an illustration.

\begin{figure}[h]
\begin{center}
\includegraphics[width=4.3cm]{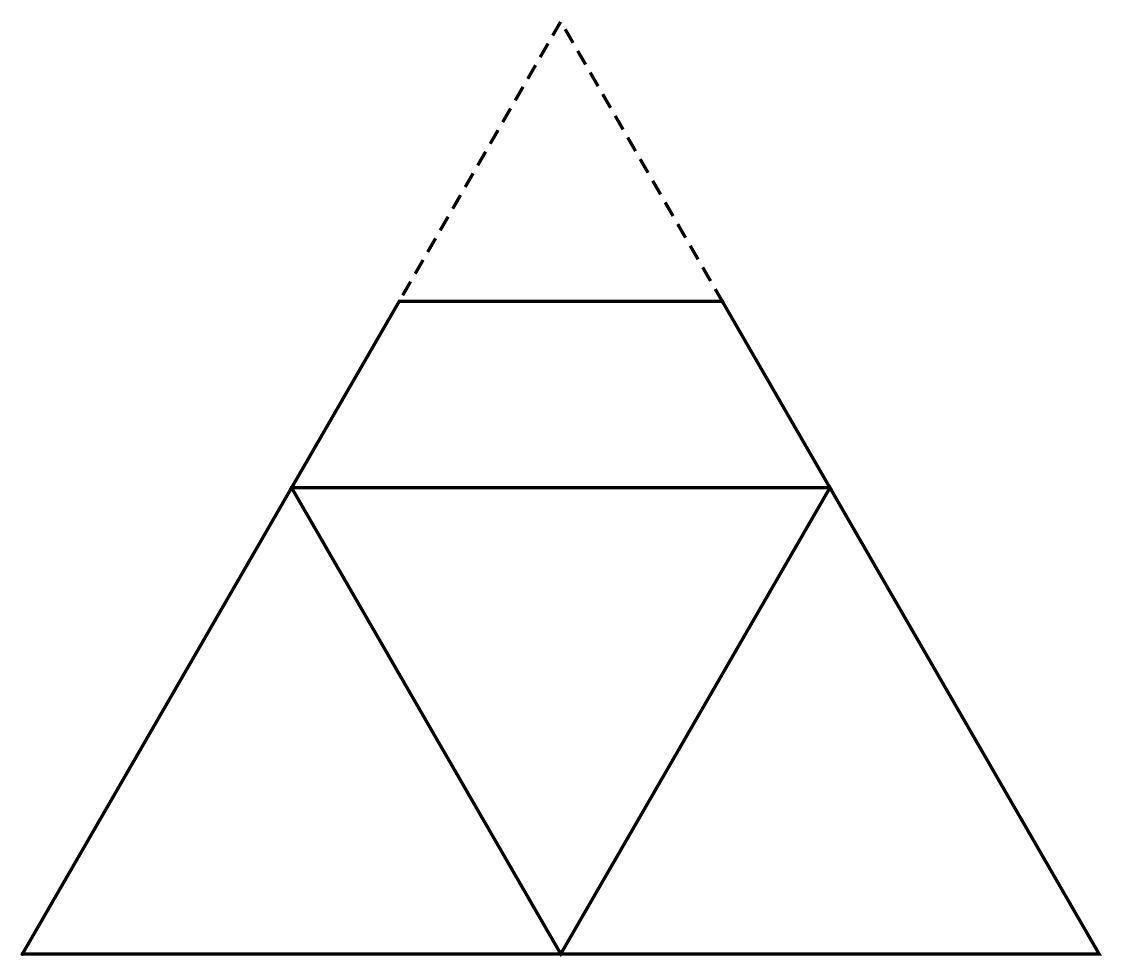}\hspace{1cm}
\includegraphics[width=4.3cm]{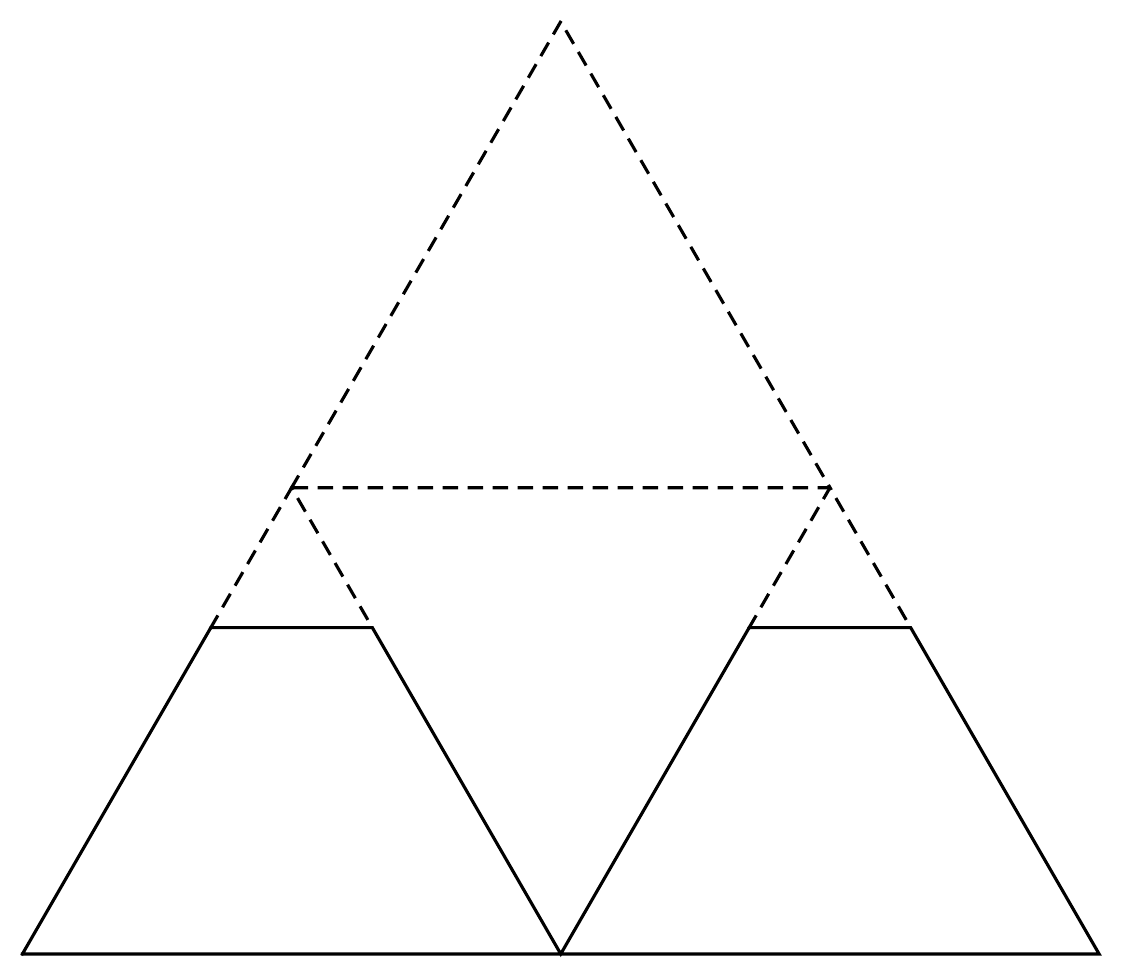}
\setlength{\unitlength}{1cm}
\begin{picture}(0,0) \thicklines
\put(-8.3,2.1){$F_0\bar{\Omega}_{S\lambda}^-$}
\put(-9.3,0.5){$F_1\mathcal{SG}$}
\put(-7.2,0.5){$F_2\mathcal{SG}$}
\put(-7.95,3.8){$q_0$}
\put(-10.1,-0.1){$q_1$}
\put(-5.8,-0.1){$q_2$}
\put(-8.15,-0.2){$F_1q_2$}
\put(-9.6,1.85){$F_0q_1$}
\put(-6.7,1.85){$F_0q_2$}

\put(-3.9,0.5){$F_1\bar{\Omega}_{S\lambda}^-$}
\put(-1.85,0.5){$F_2\bar{\Omega}_{S\lambda}^-$}
\put(-2.55,3.8){$q_0$}
\put(-4.6,-0.1){$q_1$}
\put(-0.4,-0.1){$q_2$}
\put(-2.7,-0.2){$F_1q_2$}
\end{picture}
\vspace{0.2cm}
\begin{center}
\textbf{Figure 4.2. The relationship between $\Omega_\lambda^-$ and $\Omega_{S\lambda}^-$. $e_1(\lambda)=0$ in the left one  and   $e_1(\lambda)=1$ in the right one.}
\end{center}
\end{center}
\end{figure}

 It is natural to introduce the following sets of words
\[\tilde{W}^\lambda_m=\{w\in W_m|w_k=0 \text{ if }e_k(\lambda)=0 ,\text{ and }  w_k=1 \text { or }2 \text{ if }e_k(\lambda)=1, \forall 1\leq k\leq m\},\] so that 
$$X=\bigcup _{w\in \tilde{W}_m^\lambda}X_w, \forall m\geq 0,$$ where $X_w=F_w \mathcal{SG}\cap X$. Write $\tilde{W}_*^\lambda=\bigcup_{m\geq 0}\tilde{W}_m^\lambda$.
In addition, denote $A_{\lambda,m}$ to be the closure of 
$\Omega_\lambda^-\setminus\bigcup\limits_{w\in \tilde{W}^\lambda_m} F_w\bar{\Omega}_{S^m\lambda}^-,$
which is the union of all $m$-cells contained in $\Omega_\lambda^-$. For example, 
$A_{\lambda,1}=F_1\mathcal{SG}\bigcup F_2\mathcal{SG}$ if $e_1(\lambda)=0$, and $A_{\lambda,1}=\emptyset$ if $e_1(\lambda)=1$.

\subsection{Extension Algorithm}

The Dirichlet problem on $\Omega_\lambda^-$ is stated as
\begin{equation}
\begin{cases}
\Delta u=0 \text{ in }  \Omega^-_\lambda,\\
u|_{\partial\Omega_\lambda^-}=f, f\in C(\partial\Omega^-_\lambda).
\end{cases}
\end{equation}
Analogous to the previous two sections, we only need to find an explicit algorithm for $u|_{\Omega_\lambda^-\bigcap V_1}$ in terms of the boundary data $f$, since then the problem of finding values of $u$ elsewhere in $\Omega_\lambda^-$ is essentially the same after dilation. 

From the matching condition at each vertex in  $V_1\bigcup \Omega_\lambda^-$, there exist the equations,
\begin{equation}
\begin{cases}
4u(F_1q_2)-f(q_1)-f(q_2)-u(F_0q_1)-u(F_0q_2)=0,\\
\frac{3}{5}\partial_n^\leftarrow u(F_0q_1)+\big(2u(F_0q_1)-u(F_1q_2)-f(q_1)\big)=0,\\
\frac{3}{5}\partial_n^\rightarrow u(F_0q_2)+\big(2u(F_0q_2)-u(F_1q_2)-f(q_2)\big)=0,\\
\end{cases}\text{ if }e_1(\lambda)=0,
\end{equation}
and
\begin{equation}
\partial_n^\rightarrow u(F_1q_2)+\partial_n^\leftarrow u(F_1q_2)=0,\quad\text{ if }e_1(\lambda)=1.
\end{equation}

We need to express the involved normal derivatives in terms of $u(F_0q_1)$, $u(F_0q_2)$, $u(F_1q_2)$ and the boundary data $f$. For this purpose, we introduce two important coefficients 
\begin{equation}
\eta_1(\lambda)=\partial_n^\leftarrow h^\lambda_1(q_1),\eta_2(\lambda)=-\partial_n^\rightarrow h^\lambda_1(q_2),
\end{equation}
where $h^\lambda_1$ is the harmonic function with values $h^\lambda_1(q_1)=1,h^\lambda_1(q_2)=0$ and $h^\lambda_1|_X=0$. Symmetrically, we also have
\[\eta_2(\lambda)=-\partial_n^\leftarrow h^\lambda_2(q_1),\eta_1(\lambda)=\partial_n^\rightarrow h^\lambda_2(q_2),\]
where $h^\lambda_2$ is the harmonic function with values $h^\lambda_2(q_1)=0,h^\lambda_2(q_2)=1$ and $h^\lambda_2|_X=0$. We omit the superscript $\lambda$ of $h^\lambda_i$ when there is no confusion caused. We will discuss on how to calculate these coefficients in the second part of this section. Here we only mention the following property.

\textbf{Lemma 4.1.} \textit{For $0\leq\lambda<1$, we have $\eta_1(\lambda)\geq 2$,  $0\leq \eta_2(\lambda)\leq 1$.}

\textit{Proof.} By using the maximum principle for harmonic functions, it is easy to see that $\eta_1$ is an increasing function of $\lambda$, and $\eta_2$ is a decreasing function of $\lambda$. Thus we have
\[\eta_1(\lambda)\geq \eta_1(0)=2,\quad\eta_2(\lambda)\leq \eta_2(0)=1.\]
Obviously, $\eta_2(\lambda)=-\partial_n^\rightarrow h_1(q_2)$ is nonnegative.\hfill$\square$

\textbf{Remark. } More precisely, we have $\eta_1(\lambda)+\eta_2(\lambda)\geq 3$. In fact, we just need to consider the antisymmetric harmonic function $h_1-h_2$ whose normal derivative at $q_1$ is $\eta_1(\lambda)+\eta_2(\lambda)$. Using the maximum principle on the left half part of $\Omega_\lambda^-$, one can check that $\eta_1+\eta_2$ is an increasing function of $\lambda$.

Using the coefficients $\eta_1(\lambda)$, $\eta_2(\lambda)$, we have the following two lemmas.

\textbf{Lemma 4.2.} \textit{Assume $h=h(q_1)h_1+h(q_2)h_2$. If $e_1(\lambda)=0$, then }
\[
\begin{pmatrix}
\partial_n^\leftarrow h(F_0q_1)\\
\partial_n^\rightarrow h(F_0q_2)
\end{pmatrix}=M_0^\lambda
\begin{pmatrix}
\partial_n^\leftarrow h(q_1)\\
\partial_n^\rightarrow h(q_2)
\end{pmatrix}.
\]
\textit{If $e_1(\lambda)=1$, then }
\[
\begin{pmatrix}
\partial_n^\leftarrow h(F_1q_1)\\
\partial_n^\rightarrow h(F_1q_2)
\end{pmatrix}=M_1^\lambda
\begin{pmatrix}
\partial_n^\leftarrow h(q_1)\\
\partial_n^\rightarrow h(q_2)
\end{pmatrix}, \textit{\text{ and }}
\begin{pmatrix}
\partial_n^\leftarrow h(F_2q_1)\\
\partial_n^\rightarrow h(F_2q_2)
\end{pmatrix}=M_2^\lambda
\begin{pmatrix}
\partial_n^\leftarrow h(q_1)\\
\partial_n^\rightarrow h(q_2)
\end{pmatrix}.
\]
\textit{Here the matrices $M^\lambda_i$ are }
\begin{equation}
\begin{aligned}
&M_0^\lambda=\begin{pmatrix}
\frac{3+3\eta_1(S\lambda)+3\eta_2(S\lambda)}{6+4\eta_1(S\lambda)+4\eta_2(S\lambda)} & \frac{3+\eta_1(S\lambda)+\eta_2(S\lambda)}  {6+4\eta_1(S\lambda)+4\eta_2(S\lambda)}\\
\frac{3+\eta_1(S\lambda)+\eta_2(S\lambda)}{6+4\eta_1(S\lambda)+4\eta_2(S\lambda)} &
\frac{3+3\eta_1(S\lambda)+3\eta_2(S\lambda)}{6+4\eta_1(S\lambda)+4\eta_2(S\lambda)}
\end{pmatrix},\\
 M_1^\lambda=&\begin{pmatrix}
1 & 0\\-\frac{\eta_2(S\lambda)}{2\eta_1(S\lambda)} & \frac{\eta_2(S\lambda)}{2\eta_1(S\lambda)}
\end{pmatrix}\textit{\text{ and }}
M_2^\lambda=\begin{pmatrix}
\frac{\eta_2(S\lambda)}{2\eta_1(S\lambda)} & -\frac{\eta_2(S\lambda)}{2\eta_1(S\lambda)}\\0 & 1
\end{pmatrix}.
\end{aligned}
\end{equation}

\textit{Proof.} Obviously, the normal derivatives $\partial_n^\leftarrow h(F_iq_1),\partial_n^\rightarrow h(F_iq_2)$ for $i\in \tilde{W}_1^\lambda$ are determined linearly by $\partial_n^\leftarrow h(q_1),\partial_n^\rightarrow h(q_2)$. We only need to calculate the matrices $M^\lambda_i$. In Figure 4.3, we picture  the values of the symmetric function $h_1+h_2$ and the antisymmetric function $h_1-h_2$, solved from equations (4.3) and (4.4), involving $\eta_1(S\lambda)$ and $\eta_2(S\lambda)$. 

\begin{figure}[h]
\begin{center}
\includegraphics[width=4.5cm]{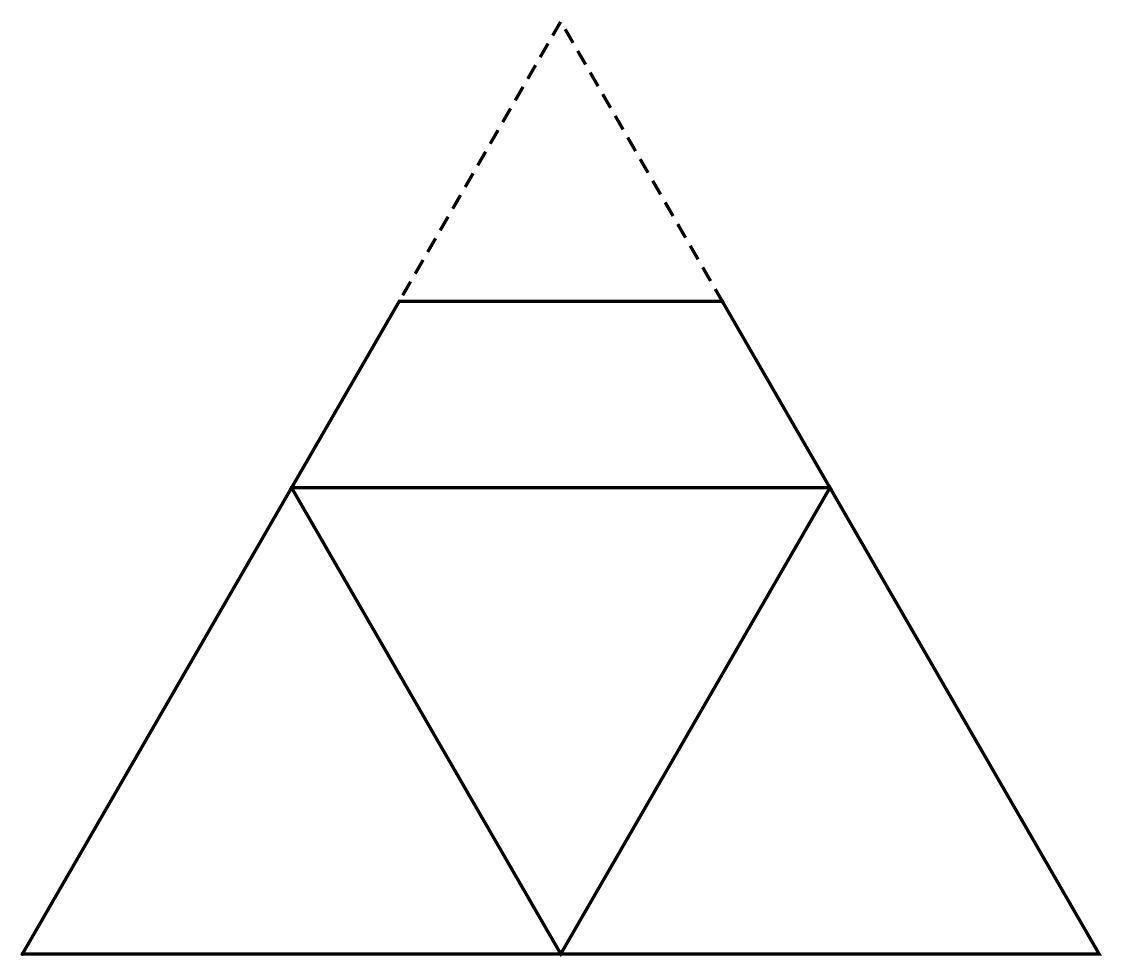}\hspace{1cm}
\includegraphics[width=4.5cm]{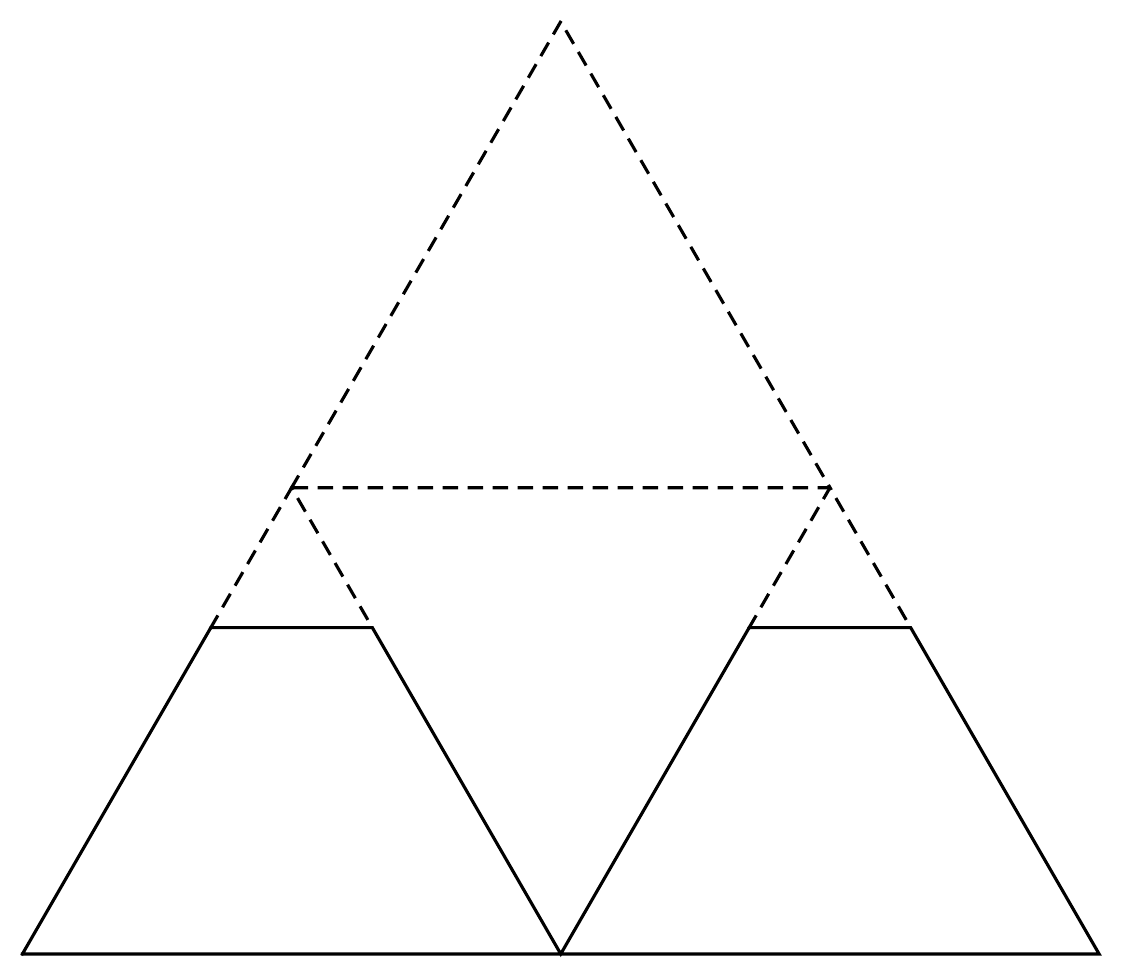}\\
\vspace{0.5cm}
\quad\includegraphics[width=4.5cm]{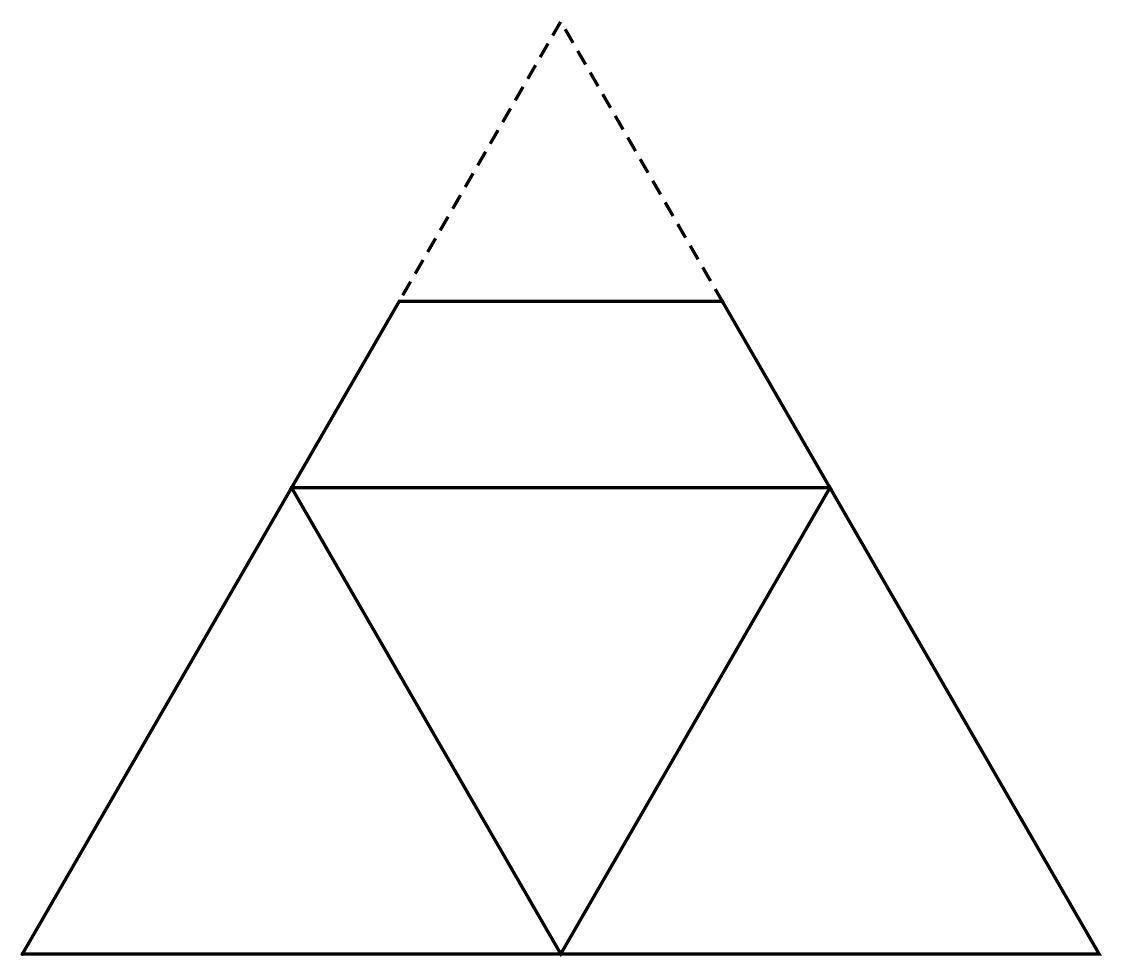}\hspace{1cm}
\includegraphics[width=4.5cm]{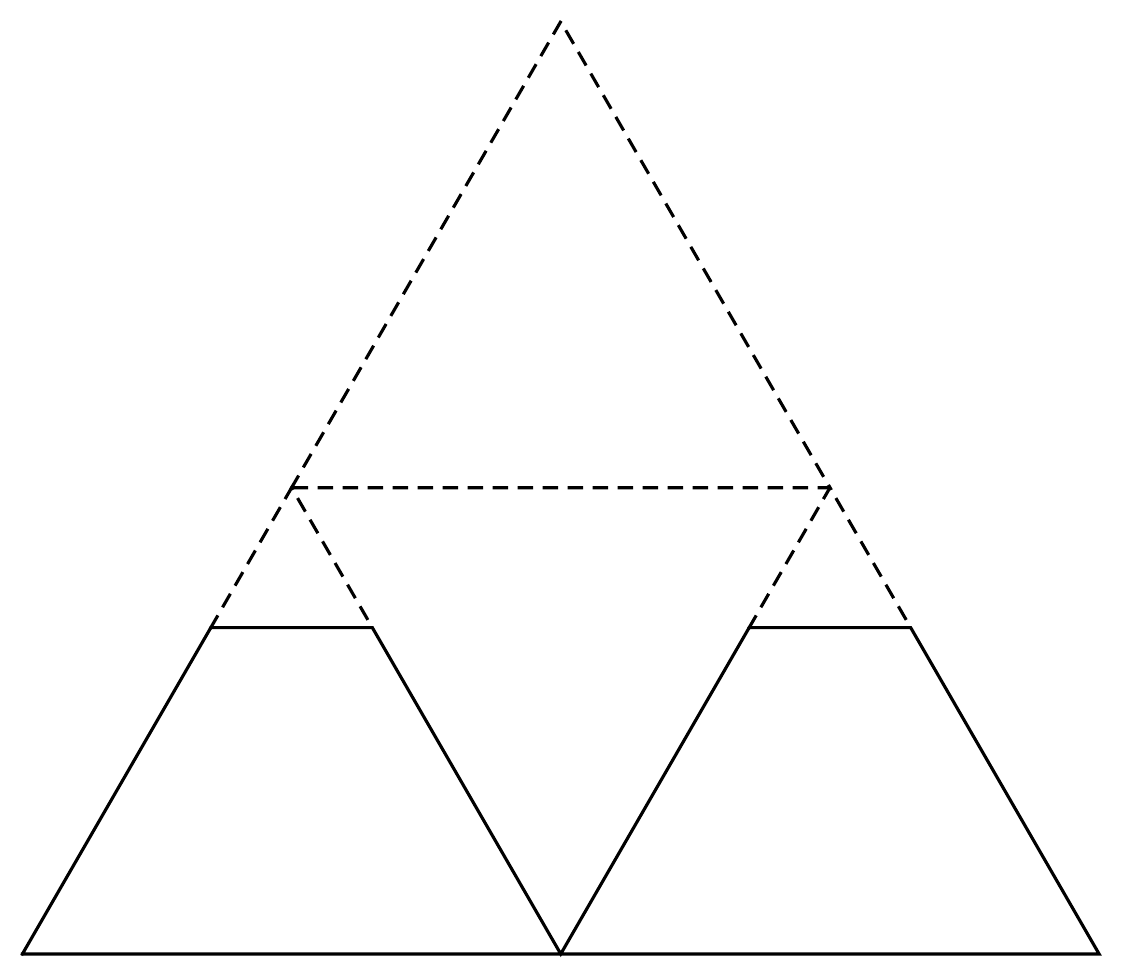}
\setlength{\unitlength}{1cm}
\begin{picture}(0,0) 
\put(-10.8,2.2){$\frac{1}{2+\eta_1(S\lambda)+\eta_2(S\lambda)}$}
\put(-7.9,2.2){$-\frac{1}{2+\eta_1(S\lambda)+\eta_2(S\lambda)}$}
\put(-8.2,-0.2){$0$}
\put(-10.5,-0.15){$1$}
\put(-6.1,-0.25){$-1$}
\put(-8.2,2.8){$0$}

\put(-10.8,6.6){$\frac{3}{3+2\eta_1(S\lambda)-2\eta_2(S\lambda)}$}
\put(-7.9,6.6){$\frac{3}{3+2\eta_1(S\lambda)-2\eta_2(S\lambda)}$}
\put(-9.3,4.2){$\frac{3+\eta_1(S\lambda)-\eta_2(S\lambda)}{3+2\eta_1(S\lambda)-2\eta_2(S\lambda)}$}
\put(-10.4,4.3){$1$}
\put(-5.9,4.3){$1$}
\put(-8.2,7.15){$0$}

\put(-4.9,-0.2){$1$}
\put(-0.4,-0.2){$-1$}
\put(-2.55,-0.2){$0$}
\put(-1.5,1.45){$0$}
\put(-3.65,1.45){$0$}

\put(-4.8,4.3){$1$}
\put(-0.3,4.3){$1$}
\put(-2.9,4.15){$\frac{\eta_2(S\lambda)}{\eta_1(S\lambda)}$}
\put(-1.5,5.9){$0$}
\put(-3.65,5.9){$0$}
\end{picture}
\vspace{0.2cm}
\begin{center}
\textbf{Figure 4.3.  $h_1+h_2$ and $h_1-h_2$ in two cases.}
\end{center}
\end{center}
\end{figure}

Hence if $e_1(\lambda)=0$, we have
\[M_0^\lambda\begin{pmatrix}
1\\1
\end{pmatrix}=
\begin{pmatrix}
1\\1
\end{pmatrix},
M_0^\lambda
\begin{pmatrix}
\frac{3+2\eta_1(S\lambda)+2\eta_2(S\lambda)}{2+\eta_1(S\lambda)+\eta_2(S\lambda)}\\-\frac{3+2\eta_1(S\lambda)+2\eta_2(S\lambda)}{2+\eta_1(S\lambda)+\eta_2(S\lambda)}
\end{pmatrix}=
\begin{pmatrix}
\frac{\eta_1(S\lambda)+\eta_2(S\lambda)}{2+\eta_1(S\lambda)+\eta_2(S\lambda)}\\-\frac{\eta_1(S\lambda)+\eta_2(S\lambda)}{2+\eta_1(S\lambda)+\eta_2(S\lambda)}
\end{pmatrix},
\]
and if $e_1(\lambda)=1$, we have
\[
M_1^\lambda \begin{pmatrix}
1\\1
\end{pmatrix}=
\begin{pmatrix}
1\\0
\end{pmatrix},
M_1^\lambda \begin{pmatrix}
\eta_1(S\lambda)\\-\eta_1(S\lambda)
\end{pmatrix}=
\begin{pmatrix}
\eta_1(S\lambda)\\-\eta_2(S\lambda)
\end{pmatrix},
\]
and similarly for $M_2^\lambda$ by symmetry considerations.\hfill$\square$

By iteratively using the above matrices, we can calculate the normal derivatives of $h$ at vertices close to $X$,
\begin{equation}
\begin{pmatrix}
\partial_n^\leftarrow  h(F_wq_1)\\
\partial_n^\rightarrow h(F_wq_2)
\end{pmatrix}
=M_w^\lambda
\begin{pmatrix}
\partial_n^\leftarrow  h(q_1)\\
\partial_n^\rightarrow h(q_2)
\end{pmatrix},
\end{equation}
for $w\in \tilde{W}^\lambda_*$ with $M_w^\lambda=M_{w_{|w|}}^{S^{|w|-1}\lambda}\cdots M_{w_2}^{S\lambda} M_{w_1}^\lambda$.  Moreover, we have

\textbf{Lemma 4.3.} \textit{Assume $h=h(q_1)h_1+h(q_2)h_2$. Then $\forall m\geq 0$,
\begin{equation}
\sum_{w\in\tilde{W}^\lambda_m} \big(\partial_n^\leftarrow  h(F_wq_1)+\partial_n^\rightarrow h(F_wq_2)\big)=\partial_n^\leftarrow  h(q_1)+\partial_n^\rightarrow h(q_2),
\end{equation}
and $\forall w\in \tilde{W}_*^\lambda$,
\begin{equation}
\partial_n^\leftarrow  h(F_wq_1)+\partial_n^\rightarrow h(F_wq_2)\geq 0,\text{ if }h(q_1)\geq 0\text{ and }h(q_2)\geq 0.\\
\end{equation}}

\textit{Proof.} By using the local Gauss-Green's formula on $A_{\lambda,m}$ we get (4.8) holds for each $m\geq 0$. Notice that for $w\in\tilde{W}_*^\lambda$, $$\partial_n^{\leftarrow}h(F_wq_1)+\partial_n^\rightarrow h(F_wq_2)=(\frac{5}{3})^{|w|}\big(h(F_wq_1)+h(F_wq_2)\big)\big(\eta_1(S^{|w|}\lambda)-\eta_2(S^{|w|}\lambda)\big).$$ By using Lemma 4.1 and the maximum principle for harmonic functions, we have both  $\eta_1(S^{|w|}\lambda)-\eta_2(S^{|w|}\lambda)\geq 0$ and $h(F_wq_1)+h(F_wq_2)\geq 0$, in case of $h(q_1)\geq 0$ and $h(q_2)\geq 0$, which gives (4.9).\hfill$\square$

According to Lemma 4.2 and 4.3, we introduce two measures on $X$.

\textbf{Definition 4.4.} \textit{Define $\mu^\lambda_1$ to be the unique probability measure on $X$ satisfying 
\[\mu^\lambda_1(X_w)=\frac{1}{\eta_1(\lambda)-\eta_2(\lambda)}\begin{pmatrix}
1,1
\end{pmatrix}M_w^\lambda 
\begin{pmatrix}
\eta_1(\lambda)\\-\eta_2(\lambda)
\end{pmatrix},\forall w\in \tilde{W}^\lambda_*,
\] 
Symmetrically, define $\mu^\lambda_2$ by
\[\mu^\lambda_2(X_w)=\frac{1}{\eta_1(\lambda)-\eta_2(\lambda)}\begin{pmatrix}
1,1
\end{pmatrix}M_w^\lambda 
\begin{pmatrix}
-\eta_2(\lambda)\\\eta_1(\lambda)
\end{pmatrix},\forall w\in \tilde{W}^\lambda_*.
\] 
}

Note that for $i=1,2$, we have
$\partial_n^\leftarrow  h_i(F_wq_1)+\partial_n^\rightarrow h_i(F_wq_2)=\big(\eta_1(\lambda)-\eta_2(\lambda)\big)\mu_i^\lambda(X_w)$.

\textbf{Theorem 4.5.} \textit{Let $u$ be a solution of the Dirichlet problem (4.2). Then
\begin{eqnarray}
\partial_n^\leftarrow u(q_1)=\eta_1(\lambda)f(q_1)-\eta_2(\lambda)f(q_2)-\big(\eta_1(\lambda)-\eta_2(\lambda)\big)\int_{X} fd\mu^\lambda_1,\\
\partial_n^\rightarrow u(q_2)=\eta_1(\lambda)f(q_2)-\eta_2(\lambda)f(q_1)-\big(\eta_1(\lambda)-\eta_2(\lambda)\big)\int_{X} fd\mu^\lambda_2.
\end{eqnarray}
In addition, if $\mathcal{E}_{\Omega_\lambda^-}(u)<\infty$, then
\begin{equation}
\mathcal{E}_{\Omega_\lambda^-}(h_1,u)=\partial_n^\leftarrow u(q_1),\text{ }\mathcal{E}_{\Omega_\lambda^-}(h_2,u)=\partial_n^\rightarrow u(q_2).
\end{equation}}

\textit{Proof.}  Using the local Gauss-Green's formula on $A_{\lambda,m}$, we have
\[
\begin{aligned}
\mathcal{E}_{A_{\lambda,m}}(h_1,u)=&\partial_n^\leftarrow h_1(q_1)f(q_1)+\partial_n^\rightarrow h_1(q_2)f(q_2)\\&-\sum_{w\in \tilde{W}_m^\lambda}\big(\partial_n^\leftarrow h_1(F_wq_1)u(F_wq_1)+\partial_n^\rightarrow h_1(F_wq_2)u(F_wq_2)\big)\\
=&\eta_1(\lambda)f(q_1)-\eta_2(\lambda)f(q_2)\\
&-\sum_{w\in \tilde{W}_m^\lambda}\frac{u(F_wq_1)+u(F_wq_2)}{2}\cdot \big(\partial_n^\leftarrow h_1(F_wq_1)+\partial_n^\rightarrow h_1(F_wq_2)\big)\\
&-\sum_{w\in \tilde{W}_m^\lambda}\frac{u(F_wq_1)-u(F_wq_2)}{2}\cdot \big(\partial_n^\leftarrow h_1(F_wq_1)-\partial_n^\rightarrow h_1(F_wq_2)\big).
\end{aligned}
\]

Analogous to the proof of Theorem 3.4, it is easy to see  that
\[\lim_{m\to\infty}\sum_{w\in \tilde{W}_m^\lambda}\frac{u(F_wq_1)+u(F_wq_2)}{2}\cdot \big(\partial_n^\leftarrow h_1(F_wq_1)+\partial_n^\rightarrow h_1(F_wq_2)\big)=\big(\eta_1(\lambda)-\eta_2(\lambda)\big)\int_{X} fd\mu^\lambda_1.\]
 On the other hand, by Lemma 4.1, we have
\[\begin{aligned}
|\partial_n^\leftarrow h_1(F_wq_1)-\partial_n^\rightarrow h_1(F_wq_2)|&=\Big|(\frac{5}{3})^{m}\big(h_1(F_wq_1)-h_1(F_wq_2)\big)\big(\eta_1(S^m\lambda)+\eta_2(S^m\lambda)\big)\Big|\\
&\leq 3(\frac{5}{3})^{m}\big(h_1(F_wq_1)+h_1(F_wq_2)\big)\big(\eta_1(S^m\lambda)-\eta_2(S^m\lambda)\big)\\&=3\big(\partial_n^\leftarrow h_1(F_wq_1)+\partial_n^\rightarrow h_1(F_wq_2)\big).
\end{aligned}\]
So by Lemma 4.3,  $\sum_{w\in \tilde{W}_m^\lambda}|\partial_n^\leftarrow h_1(F_wq_1)-\partial_n^\rightarrow h_1(F_wq_2)|\leq 3\big(\eta_1(\lambda)-\eta_2(\lambda)\big)$. Noticing that $u$ is uniformly continuous on $\bar{\Omega}_\lambda^-$, we have
\[
\lim_{m\to\infty}\sum_{w\in \tilde{W}_m^\lambda}\frac{u(F_wq_1)-u(F_wq_2)}{2}\cdot \big(\partial_n^\leftarrow h_1(F_wq_1)-\partial_n^\rightarrow h_1(F_wq_2)\big)=0.
\]

Combining the above facts, we get
\begin{equation}
\lim_{m\to\infty}\mathcal{E}_{A_{\lambda,m}}(h_1, u)= \eta_1(\lambda)f(q_1)-\eta_2(\lambda)f(q_2)-\big(\eta_1(\lambda)-\eta_2(\lambda)\big)\int_{X} fd\mu^\lambda_1.
\end{equation}
Similarly, 
\begin{equation}
\lim_{m\to\infty}\mathcal{E}_{A_{\lambda,m}}(h_2, u)= \eta_1(\lambda)f(q_2)-\eta_2(\lambda)f(q_1)-\big(\eta_1(\lambda)-\eta_2(\lambda)\big)\int_{X} fd\mu^\lambda_2.
\end{equation}

For the rest part of the theorem, define a sequence of harmonic functions $u_n$ such that $u_n|_{X_\tau}$ is constant for each $\tau\in \tilde{W}_n^\lambda$, and $u_n$ converges to $u$ uniformly as $n\rightarrow\infty$. Then
$$
\mathcal{E}_{A_{\lambda,m}}(h_1,u_n)=\mathcal{E}_{A_{\lambda,n}}(h_1, u_n)+\sum_{\tau\in \tilde{W}_n^\lambda} \mathcal{E}_{F_\tau A_{S^n\lambda,m-n}}(u_n, h_1)\to \partial_n^\leftarrow u_n(q_1)
$$
as $m\to\infty$, where we use (4.13) and (4.14) on $A_{S^n\lambda,m-n}$, noticing that $u_n\circ F_\tau$ is a linear combination of $h_1^{S^n\lambda}$ and $h_2^{S^n\lambda}$. Combining with (4.13) again, we get 
\[\partial_n^\leftarrow u_n(q_1)=\eta_1(\lambda)u_n(q_1)-\eta_2(\lambda)u_n(q_2)-\big(\eta_1(\lambda)-\eta_2(\lambda)\big)\int_{X} u_nd\mu^\lambda_1.\]
Taking $n\to\infty$, we have proved (4.10). Similarly, we also have (4.11).
\hfill$\square$

Combining Theorem 4.5 with (4.3) and (4.4), after an easy calculation, we finally get the following extension algorithm.

\textbf{Theorem 4.6.(Extension Algorithm)} \textit{There exists a unique solution of the Dirichlet problem (4.2). In addition, we have the following formulas for $u|_{V_1\bigcap \Omega_\lambda^-}$.}

\textit{If $e_1(\lambda)=0$, then}
\[
\begin{aligned}
u(F_0q_1)=&\frac{9+5\eta_1(S\lambda)+\eta_2(S\lambda)}{4\eta_1(S\lambda)^2+14\eta_1(S\lambda)-2\eta_2(S\lambda)-4\eta_2(S\lambda)^2+12}f(q_1)\\&+
\frac{3+\eta_1(S\lambda)+5\eta_2(S\lambda)}{4\eta_1(S\lambda)^2+14\eta_1(S\lambda)-2\eta_2(S\lambda)-4\eta_2(S\lambda)^2+12}f(q_2)\\&+
\frac{\big(7+4\eta_1(S\lambda)\big)\big(\eta_1(S\lambda)-\eta_2(S\lambda)\big)}{4\eta_1(S\lambda)^2+14\eta_1(S\lambda)-2\eta_2(S\lambda)-4\eta_2(S\lambda)^2+12}\int_{X_{S\lambda}^-} f\circ F_0 d\mu^{S\lambda}_1\\
&+\frac{\big(1+4\eta_2(S\lambda)\big)\big(\eta_1(S\lambda)-\eta_2(S\lambda)\big)}{4\eta_1(S\lambda)^2+14\eta_1(S\lambda)-2\eta_2(S\lambda)-4\eta_2(S\lambda)^2+12}\int_{X_{S\lambda}^-} f\circ F_0 d\mu^{S\lambda}_2,
\end{aligned}
\]
\textit{and $u(F_1q_2)=\frac{1}{4}\big(u(F_0q_1)+u(F_0q_2)+f(q_1)+f(q_2)\big)$. The formula of $u(F_0q_2)$ is symmetrical to that of $u(F_0q_1)$.}

\textit{If $e_2(\lambda)=1$, then}
\[u(F_1q_2)=\frac{\eta_2(S\lambda)}{2\eta_1(S\lambda)}\big(f(q_1)+f(q_2)\big)+\frac{\eta_1(S\lambda)-\eta_2(S\lambda)}{2\eta_1(S\lambda)}\big(\int_{X_{S\lambda}^-}f\circ F_1d\mu^{S\lambda}_2+\int_{X_{S\lambda}^-} f\circ F_2d\mu^{S\lambda}_1\big).\]

\subsection{The calculation of $\eta$}
In this section, we focus on the calculation of $\eta_1(\lambda),\eta_2(\lambda)$. In particular, we will prove the following theorem.

\textbf{Theorem 4.7.} \textit{(a) $\eta_1$ is a continuous increasing function of $\lambda$ on $[0,1)$, while $\eta_2$ is a continuous decreasing function of $\lambda$.\\
(b) For $0\leq\lambda<1$, $\big(\eta_1(\lambda), \eta_2(\lambda)\big)=T_{e_1(\lambda)}\big(\eta_1(S\lambda), \eta_2(S\lambda)\big)$ with
 \begin{equation}T_0(x,y)=\big(\frac{5}{6}\cdot\frac{3+2x+2y}{2+x+y}+\frac{5}{2}\cdot\frac{x-y}{3+2x-2y},   \frac{5}{6}\cdot\frac{3+2x+2y}{2+x+y}-\frac{5}{2}\cdot\frac{x-y}{3+2x-2y}\big),
 \end{equation}
and
\begin{equation}T_1(x,y)=\big(\frac{5}{3}(x-\frac{y^2}{2x}), \frac{5y^2}{6x}\big).
\end{equation}
In addition, for any fixed positive numbers $c_1>c_2$, we have
\begin{equation}
\big(\eta_1(\lambda),\eta_2(\lambda)\big)=\lim_{m\to\infty} T_{e_1(\lambda)}\circ T_{e_2(\lambda)}\circ\cdots\circ T_{e_m(\lambda)}(c_1,c_2).
\end{equation}}

\textit{Proof of Theorem 4.7(a).}  By Lemma 4.1, $\eta_1$ is increasing and $\eta_2$ is decreasing for $ \lambda\in[0,1)$. So we only need to show $\eta_1$ and $\eta_2$ are continuous functions. 

\textbf{Claim 1.}  \textit{For $0\leq\lambda<1$, we have
\[0\leq h_1(F_wq_i)\leq (\frac{3}{5})^{|w|}\big(\eta_1(\lambda)-\eta_2(\lambda)\big),\forall w\in \tilde{W}^\lambda_*,\forall i=1,2.\]}

\textit{Proof of Claim 1.} Let $m=|w|$. Noticing that $h_1(F_wq_i)\geq 0$ by the maximum principle for harmonic functions, by Lemma 4.1, we have
$$
\partial_n^\leftarrow h_1(F_wq_1)+\partial_n^\rightarrow h_1(F_wq_2)=(\frac{5}{3})^m\big(\eta_1(S^m\lambda)-\eta_2(S^m\lambda)\big)\big(h_1(F_wq_1)+h_1(F_wq_2)\big)\geq 0.$$
On the other hand, by Lemma 4.3, 
$$
\partial_n^\leftarrow h_1(F_wq_1)+\partial_n^\rightarrow h_1(F_wq_2)\leq\partial_n^\leftarrow h_1(q_1)+\partial_n^\rightarrow h_1(q_2)=\eta_1(\lambda)-\eta_2(\lambda).
$$

Combining the above two inequalities, using Lemma 4.1 again, we get
$$
h_1(F_wq_i)\leq (\frac{3}{5})^m\frac{\eta_1(\lambda)-\eta_2(\lambda)}{\eta_1(S^m\lambda)-\eta_2(S^m\lambda)}\leq (\frac{3}{5})^{m}\big(\eta_1(\lambda)-\eta_2(\lambda)\big).\hspace{3cm}\square
$$

Without loss of generality, we assume $1-\lambda> 2^{-j}$ for some positive integer $j$. Then we have $1\leq \eta_1(\lambda)-\eta_2(\lambda)\leq \eta_1(1-2^{-j})-\eta_2(1-2^{-j})$ and $F_1^j\mathcal{SG}\subset\Omega_\lambda^-$, $F_2^j\mathcal{SG}\subset \Omega_\lambda^-$. 

\textbf{Claim 2.} \textit{If $|\lambda_a-\lambda_b|\leq 2^{-m}$ with $0\leq\lambda_a,\lambda_b< 1-2^{-j}$, we have
$$
|\eta_i(\lambda_a)-\eta_i(\lambda_b)|\leq 4(\frac{3}{5})^{m-j}\big(\eta_1(1-2^{-j})-\eta_2(1-2^{-j})\big) \text{ for } i=1,2.$$}

\textit{Proof of Claim 2.} Assume $\lambda_a\leq\lambda_b$. 
First, consider the special case that
\[\lambda_b=2^{-m}+\sum_{k=1}^{m}e_k(\lambda_a)2^{-k}.\]
In this case, $\bar{\Omega}_{\lambda_b}^-=A_{\lambda_a,m}.$ See Figure 4.4 for an illustration.
\begin{figure}[h]
\begin{center}
\includegraphics[width=5cm]{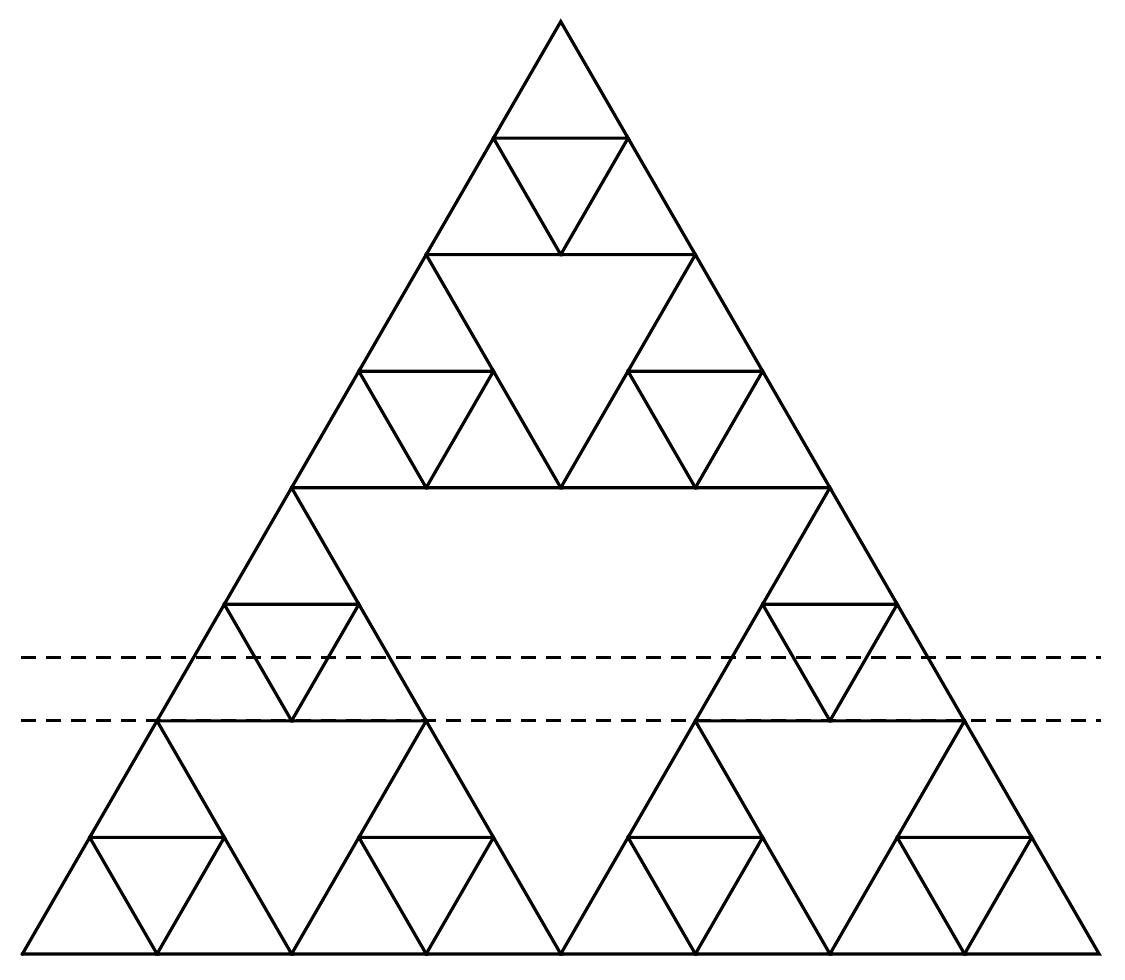}
\setlength{\unitlength}{1cm}
\begin{picture}(0,0) 
\put(-5.7,1.4){$X_{\lambda_a}$}
\put(-5.7,1.0){$X_{\lambda_b}$}
\end{picture}

\begin{center}
\textbf{Figure 4.4. A special case for $\lambda_a,\lambda_b$. ($m=3$)}
\end{center}
\end{center}
\end{figure}
By Claim 1, we have
\[\|(h_1^{\lambda_a}-h_1^{\lambda_b})|_{\partial \Omega_{\lambda_b}^-}\|_\infty\leq (\frac{3}{5})^{m}\big(\eta_1(\lambda_a)-\eta_2(\lambda_a)\big)\]
as $\partial \Omega_{\lambda_b}^-\subset\{q_1,q_2\}\bigcup\{F_wq_1,F_wq_2\}_{w\in \tilde{W}^{\lambda_a}_m}$. Thus, 
\[\|(h_1^{\lambda_a}-h_1^{\lambda_b})|_{\Omega_{\lambda_b}^-}\|_\infty\leq (\frac{3}{5})^m\big(\eta_1(\lambda_a)-\eta_2(\lambda_a)\big)\leq \big(\frac{3}{5})^m(\eta_1(1-2^{-j})-\eta_2(1-2^{-j})\big)\]
by the maximum principle for harmonic functions.

Next, for the general case, we introduce $\tilde{\lambda}_a,\tilde{\lambda}_b$ as follows,
\[\tilde{\lambda}_a=2^{-m}+\sum_{k=1}^{m}e_k(\lambda_a)2^{-k},\quad \tilde{\lambda}_b=2^{-m}+\sum_{k=1}^{m}e_k(\lambda_b)2^{-k}.\]
See Figure 4.5 for an illustration.
\begin{figure}[h]
\begin{center}
\includegraphics[width=5cm]{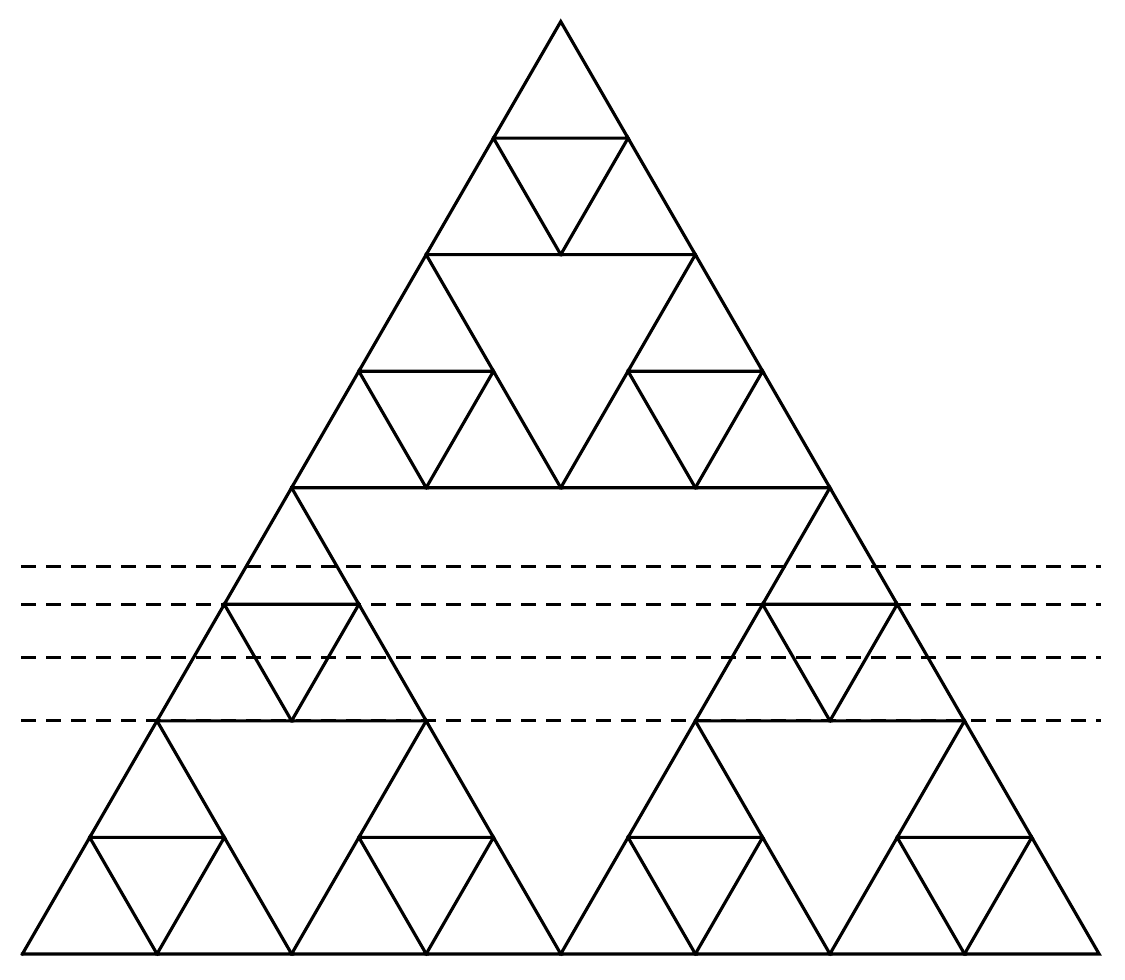}
\setlength{\unitlength}{1cm}
\begin{picture}(0,0) 
\put(-5.7,1.8){$X_{\lambda_a}$}
\put(-5.7,1.4){$X_{\lambda_b}$}

\put(-0.3,1.55){$X_{\tilde{\lambda}_a}$}
\put(-0.3,1.05){$X_{\tilde{\lambda}_b}$}
\end{picture}

\begin{center}
\textbf{Figure 4.5.  A general case for $\lambda_a,\lambda_b$. ($m=3$)}\end{center}
\end{center}
\end{figure}
It is easy to check that $\tilde{\lambda}_a=\tilde{\lambda}_b$ or 
$\tilde{\lambda}_b=\tilde{\lambda}_a+2^{-m}=2^{-m}+\sum_{k=1}^{m}e_k(\tilde{\lambda}_a)2^{-k}$, and $\lambda_a\leq\lambda_b\leq\tilde{\lambda}_b$. Thus we have
\[\begin{aligned}
\|(h_1^{\lambda_a}-h_1^{\lambda_b})|_{\Omega_{\tilde{\lambda}_b}^-}\|_\infty  &\leq \|(h_1^{\lambda_a}-h_1^{\tilde{\lambda}_b})|_{\Omega_{\tilde{\lambda}_b}^-}\|_\infty\\
&\leq\|(h_1^{\lambda_a}-h_1^{\tilde{\lambda}_a})|_{\Omega_{\tilde{\lambda}_b}^-}\|_\infty+\|(h_1^{\tilde{\lambda}_a}-h_1^{\tilde{\lambda}_b})|_{\Omega_{\tilde{\lambda}_b}^-}\|_\infty\\
&\leq 2(\frac{3}{5})^m\big(\eta_1(1-2^{-j})-\eta_2(1-2^{-j})\big).
\end{aligned}\]

Using the above estimates, we finally have
\[\begin{aligned}
|\eta_1(\lambda_a)-\eta_1(\lambda_b)|=&|\partial_n^\leftarrow h_1^{\lambda_a}(q_1)-\partial_n^\leftarrow h_1^{\lambda_b}(q_1)|\\
=&\Big|(\frac{5}{3})^j\big(2-h_1^{\lambda_a}(F_1^jq_0)-h_1^{\lambda_a}(F_1^jq_2)\big)-(\frac{5}{3})^j\big(2-h_1^{\lambda_b}(F_1^jq_0)-h_1^{\lambda_b}(F_1^jq_2)\big)\Big|\\
\leq& 4(\frac{3}{5})^{m-j}\big(\eta_1(1-2^{-j})-\eta_2(1-2^{-j})\big).
\end{aligned}\]
Similarly, we have the same inequality for $\eta_2$.\hfill$\square$

From Claim 2,  we have that $\eta_1$ and $\eta_2$ are continuous on $[0,1-2^{-j})$ for any positive integer $j$. Thus, we have proved (a).\hfill$\square$

\textit{Proof of Theorem 4.7(b).} Looking at the functions $h_1+h_2$ and $h_1-h_2$ pictured in 
Figure 4.3, by computing their normal derivatives at $q_1$, we get
\[\begin{cases}
\eta_1(\lambda)+\eta_2(\lambda)=\frac{5}{3}\cdot\frac{3+2\eta_1(S\lambda)+2\eta_2(S\lambda)}{2+\eta_1(S\lambda)+\eta_2(S\lambda)},\\
\eta_1(\lambda)-\eta_2(\lambda)=\frac{5}{3}\cdot\frac{3\eta_1(S\lambda)-3\eta_2(S\lambda)}{3+2\eta_1(S\lambda)-2\eta_2(S\lambda)},
\end{cases}\text{ if } e_1(\lambda)=0,\]

and

\[\begin{cases}
\eta_1(\lambda)+\eta_2(\lambda)=\frac{5}{3}\eta_1(S\lambda),\\
\eta_1(\lambda)-\eta_2(\lambda)= \frac{5}{3}\big(\eta_1(S\lambda)-\frac{\eta_2(S\lambda)^2}{\eta_1(S\lambda)}\big),
\end{cases}\text{ if } e_1(\lambda)=1.\]
The above equations lead to 
\begin{equation}\big(\eta_1(\lambda), \eta_2(\lambda)\big)=T_{e_1(\lambda)}\big(\eta_1(S\lambda), \eta_2(S\lambda)\big).
\end{equation}

For the rest of Theorem 4.7(b), we need the following definition. 
 
\textbf{Notation.} \textit{(a) For $0\leq \lambda<1$ and fixed positive numbers $c_1>c_2$, define a sequence of resistance forms $\mathcal{E}_m^{(c_1,c_2)}(\cdot,\cdot)$ on $V_m^\lambda=(V_m\bigcap \bar{\Omega}_\lambda^-)\bigcup\{F_wq_0\}_{w\in \tilde{W}^\lambda_m}$ with the conductances
\[c^{\lambda,m}_{x,y}(c_1,c_2)=
\begin{cases}
0,\qquad\qquad\qquad\text{if }x\nsim_m y,\\
(\frac{5}{3})^m(c_1-c_2)\text{, if }\{x,y\}=\{F_wq_0,F_wq_1\} \text{ or }\{x,y\}=\{F_wq_0,F_wq_2\}
\\\qquad\qquad\qquad\qquad\qquad\qquad\text{ for some } w\in \tilde{W}_m^\lambda,\\
(\frac{5}{3})^mc_2,\qquad\quad\text{ if }\{x,y\}=\{F_wq_1,F_wq_2\} \text{ for some } w\in \tilde{W}_m^\lambda,\\
(\frac{5}{3})^m,\qquad\qquad \text{ Otherwise.}
\end{cases}
\]}

\textit{(b) Let $h^{\lambda,(c_1,c_2)}_{1,m}$ be a sequence of functions  on $V_m^\lambda$ harmonic with respect to $\mathcal{E}_m^{(c_1,c_2)}(\cdot,\cdot)$, assuming boundary values} 
\[h^{\lambda,(c_1,c_2)}_{1,m}(q_1)=1, h^{\lambda,(c_1,c_2)}_{1,m}(q_2)=0, h^{\lambda,(c_1,c_2)}_{1,m}(F_wq_0)=0,\forall w\in\tilde{W}^\lambda_m.\]

\textit{ Still denote by $h^{\lambda,(c_1,c_2)}_{1,m}$ the piecewise harmonic function which assumes the same value on $V_m^\lambda$ and takes harmonic extension  elsewhere in $\Omega_\lambda^-\bigcup\{F_w\mathcal{SG}\}_{w\in \tilde{W}_m^\lambda}$.}\\

In Figure 4.6, we give an example of $V_m^\lambda$ together with some conductances. We abbreviate $h^{\lambda,(c_1,c_2)}_{1,m}$ to $h_{1,m}$ when there is no confusion caused. By Theorem 4.5, one can easily check that $h_{1,m}|_{V_m\cap\bar{\Omega}_\lambda^-}=h_1|_{V_m\cap\bar{\Omega}_\lambda^-}$, when $(c_1,c_2)=\big(\eta_1(S^m\lambda), \eta_2(S^m\lambda)\big)$.

\begin{figure}[h]
\begin{center}
\includegraphics[width=11cm]{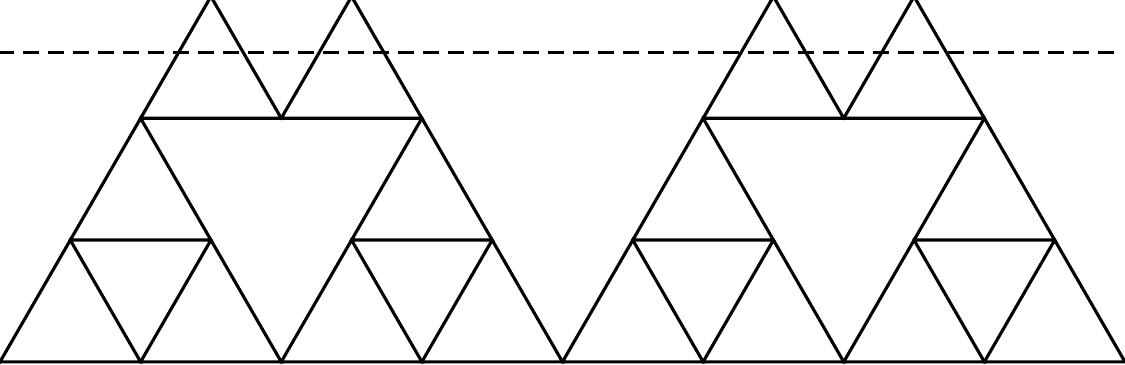}
\setlength{\unitlength}{1cm}
\begin{picture}(0,0) 
\put(-8.25,2.55){$(\frac{5}{3})^3c_2$}
\put(-8.25,2.55){$(\frac{5}{3})^3c_2$}

\put(-2.1,3.3){$(\frac{5}{3})^3(c_1-c_2)$}
\put(-2.1,3.3){$(\frac{5}{3})^3(c_1-c_2)$}

\put(-5.3,0.2){$(\frac{5}{3})^3$}
\put(-5.3,0.2){$(\frac{5}{3})^3$}

\put(-6.2,0.8){$(\frac{5}{3})^3$}
\put(-6.2,0.8){$(\frac{5}{3})^3$}

\put(-10.9,1.6){$(\frac{5}{3})^3$}
\put(-10.9,1.6){$(\frac{5}{3})^3$}

\put(-11.7,2.9){$L_\lambda$}
\end{picture}
\vspace{0.1cm}
\begin{center}
\textbf{Figure 4.6. $V_m^\lambda$ and some conductances. ($\frac{5}{8}<\lambda<\frac{3}{4},m=3$) }
\end{center}
\end{center}
\end{figure}

We still assume $1-\lambda>2^{-j}$ for some positive integer $j$.  

\textbf{Claim 3.} \textit{For  $c_1>c_2>0$, we have \[\begin{aligned}
&\Big((\frac{5}{3})^j\big(2-h_{1,m}(F_1^jq_0)-h_{1,m}(F_1^jq_2)\big),(\frac{5}{3})^j\big(h_{1,m}(F_2^jq_0)+h_{1,m}(F_2^jq_1)\big)\Big)
\\=& T_{e_1(\lambda)}\circ T_{e_2(\lambda)}\circ\cdots\circ T_{e_m(\lambda)}(c_1,c_2).
\end{aligned}\]}

\textit{Proof of Claim 3.} It is easy to see the claim holds by inductively using (4.18) when $(c_1,c_2)=\big(\eta_1(S^m\lambda),\eta_2(S^m\lambda)\big)$, since then $h_{1,m}|_{V_m\cap\bar{\Omega}_\lambda^-}=h_1|_{V_m\cap\bar{\Omega}_\lambda^-}$. For general $c_1,c_2$,  the claim still holds in a completely similar way. \hfill$\square$

\textbf{Claim 4.} \textit{Write $(c_{1,m},c_{2,m})=T_{e_1(\lambda)}\circ T_{e_2(\lambda)}\circ\cdots\circ T_{e_m(\lambda)}(c_1,c_2)$ for short. Then
\[h_{1,m}(F_wq_i)\leq (\frac{3}{5})^m\frac{c_{1,m}-c_{2,m}}{c_1-c_2}\leq 4(\frac{3}{5})^{m-j}\frac{1}{c_1-c_2}, \forall w\in \tilde{W}^\lambda_m, \forall i=1,2.\]}

\textit{Proof of Claim 4.} The first inequality is obtained analogously to the proof of Claim 1.
The second inequality follows from Claim 3 and the fact $\|h_{1,m}\|_\infty\leq 1$. \hfill$\square$

Notice that $h^{\lambda,(c_1,c_2)}_{1,m}$ and $h^{\lambda,(\eta_1(S^m\lambda),\eta_2(S^m\lambda))}_{1,m}$ satisfy the same mean value equations on $V_m^\lambda\setminus(\bigcup_{w\in \tilde{W}^\lambda_m} F_wV_0)$. By Claim 4, we have the following estimate due to the maximum principle for harmonic functions,
\[\big|(h^{\lambda,(c_1,c_2)}_{1,m}-h^{\lambda,(\eta_1(S^m\lambda),\eta_2(S^m\lambda))}_{1,m})(F_i^jq_{i'})\big|\leq 4(\frac{3}{5})^{m-j}(\frac{1}{c_1-c_2}+1),
\]
for $(i,i')\in\{(1,0),(1,2),(2,0),(2,1)\}$, where we use the fact that $\eta_1(S^m\lambda)-\eta_2(S^m\lambda)\geq 1$.

 Thus, by Claim 3, we have
\[\begin{aligned}
&\Big\|T_{e_1(\lambda)}\circ T_{e_2(\lambda)}\circ\cdots\circ T_{e_m(\lambda)}(c_1,c_2)-\big(\eta_1(\lambda),\eta_2(\lambda)\big)\Big\|\\
=&\Big\|T_{e_1(\lambda)}\circ T_{e_2(\lambda)}\circ\cdots\circ T_{e_m(\lambda)}(c_1,c_2)-T_{e_1(\lambda)}\circ T_{e_2(\lambda)}\circ\cdots\circ T_{e_m(\lambda)}\big(\eta_1(S^m\lambda),\eta_2(S^m\lambda)\big)\Big\|\\
\leq & 8(\frac{3}{5})^{m-2j}(\frac{1}{c_1-c_2}+1)
\end{aligned}\]
with $\|(a,b)\|=\max\{|a|,|b|\}$, which yields (4.17) as $m\rightarrow\infty$ for $0\leq\lambda<1-2^{-j}$. Noticing that $j$ is arbitrary, we complete the proof of Theorem 4.7(b).\hfill$\square$

\subsection{Some calculations on $M_w^\lambda$}
In this subsection, we consider two special cases,  $e_1(\lambda)=e_2(\lambda)=\cdots=e_m(\lambda)=0$ or $e_1(\lambda)=e_2(\lambda)=\cdots=e_m(\lambda)=1$ for some $m$. It would be convenient to get a direct expression for $M_w^\lambda,w\in \tilde{W}_m^\lambda$ and $\underbrace{T_0\circ T_0\cdots\circ T_0}_{m\text{'s}}$ or $\underbrace{T_1\circ T_1\cdots\circ T_1}_{m\text{'s}}$. 

\textbf{Proposition 4.8.} \textit{Assume $e_1(\lambda)=e_2(\lambda)=\cdots=e_m(\lambda)=0$. Then we have}

 $$\textit{\text{(a)}} \hspace{0.5cm} \eta_1(\lambda)+\eta_2(\lambda)=3+\frac{14\big(\eta_1(S^m\lambda)+\eta_2(S^m\lambda)-3\big)}{3(15^m-1)\big(\eta_1(S^m\lambda)+\eta_2(S^m\lambda)-3\big)+14\cdot 15^m},$$
$$\eta_1(\lambda)-\eta_2(\lambda)=\frac{\eta_1(S^m\lambda)-\eta_2(S^m\lambda)}{(1-(\frac{3}{5})^m)\big(\eta_1(S^m\lambda)-\eta_2(S^m\lambda)\big)+(\frac{3}{5})^m},\hspace{0.7cm}$$

\textit{
\[\textit{\text{(b)}} \hspace{0.5cm} M^\lambda_{0^{m}}=\begin{pmatrix}a&b\\b&a\end{pmatrix}\text{ with } a+b=1 \text{ and }\hspace{4.9cm}\]
}
\[a-b=\frac{14\cdot 5^m\big(\eta_1(S^m\lambda)+\eta_2(S^m\lambda)\big)}{(9\cdot 15^m+5)\big(\eta_1(S^m\lambda)+\eta_2(S^m\lambda)\big)+15(15^m-1)},\hspace{0.4cm}\] \textit{where we use the notation $0^{m}$ to represent the word $00\cdots 0$ of length $m$. }

\textit{Proof.} (a) Without loss of generality, assume $\lambda\neq 0$. By (4.15), we have
\[\begin{cases}
\eta_1(S^{k-1}\lambda)+\eta_2(S^{k-1}\lambda)=\frac{5}{3}\cdot\frac{3+2\eta_1(S^k\lambda)+2\eta_2(S^k\lambda)}{2+\eta_1(S^k\lambda)+\eta_2(S^k\lambda)},\\
\eta_1(S^{k-1}\lambda)-\eta_2(S^{k-1}\lambda)=\frac{5}{3}\cdot\frac{3\eta_1(S^k\lambda)-3\eta_2(S^k\lambda)}{3+2\eta_1(S^k\lambda)-2\eta_2(S^k\lambda)},
\end{cases}\forall 1\leq k\leq m.\]
Then
\[
\eta_1(S^{k-1}\lambda)+\eta_2(S^{k-1}\lambda)-3=\frac{\eta_1(S^k\lambda)+\eta_2(S^k\lambda)-3}{15+3\big(\eta_1(S^k\lambda)+\eta_2(S^k\lambda)-3\big)}, \forall 1\leq k\leq m.\]
Noticing that by the Remark below Lemma 4.1, $\eta_1(\lambda)+\eta_2(\lambda)>3$ whenever $0<\lambda<1$, we have
\[\begin{aligned}
\frac{1}{\eta_1(\lambda)+\eta_2(\lambda)-3}&=3+\frac{15}{\eta_1(S\lambda)+\eta_2(S\lambda)-3}\\
&=3(1+15+\cdots+15^{m-1})+\frac{15^m}{\eta_1(S^m\lambda)+\eta_2(S^m\lambda)-3}.
\end{aligned}\]
So
\[\eta_1(\lambda)+\eta_2(\lambda)=3+\frac{14\big(\eta_1(S^m\lambda)+\eta_2(S^m\lambda)-3\big)}{3(15^m-1)\big(\eta_1(S^m\lambda)+\eta_2(S^m\lambda)-3\big)+14\cdot 15^m}.\]

Similarly, we also have
\[\eta_1(\lambda)-\eta_2(\lambda)=\frac{\eta_1(S^m\lambda)-\eta_2(S^m\lambda)}{\big(1-(\frac{3}{5})^m\big)\big(\eta_1(S^m\lambda)-\eta_2(S^m\lambda)\big)+(\frac{3}{5})^m}.\]

(b) We only need to calculate the eigenvalues of $M^\lambda_{0^{m}}$ corresponding to the vectors $\begin{pmatrix}1\\1\end{pmatrix}$ and $\begin{pmatrix}1\\-1\end{pmatrix}$. Obviously,
$M^\lambda_{0^{m}}\begin{pmatrix}
1\\1
\end{pmatrix}=\begin{pmatrix}
1\\1
\end{pmatrix}.$
Let $h$  be an antisymmetric harmonic function with $h|_X=0$. Set $a_k=h(F_0^{m-k}q_1)=-h(F_0^{m-k}q_2).$ See Figure 4.7 for the values of $h$.

\begin{figure}[h]
\begin{center}
\includegraphics[width=5.5cm]{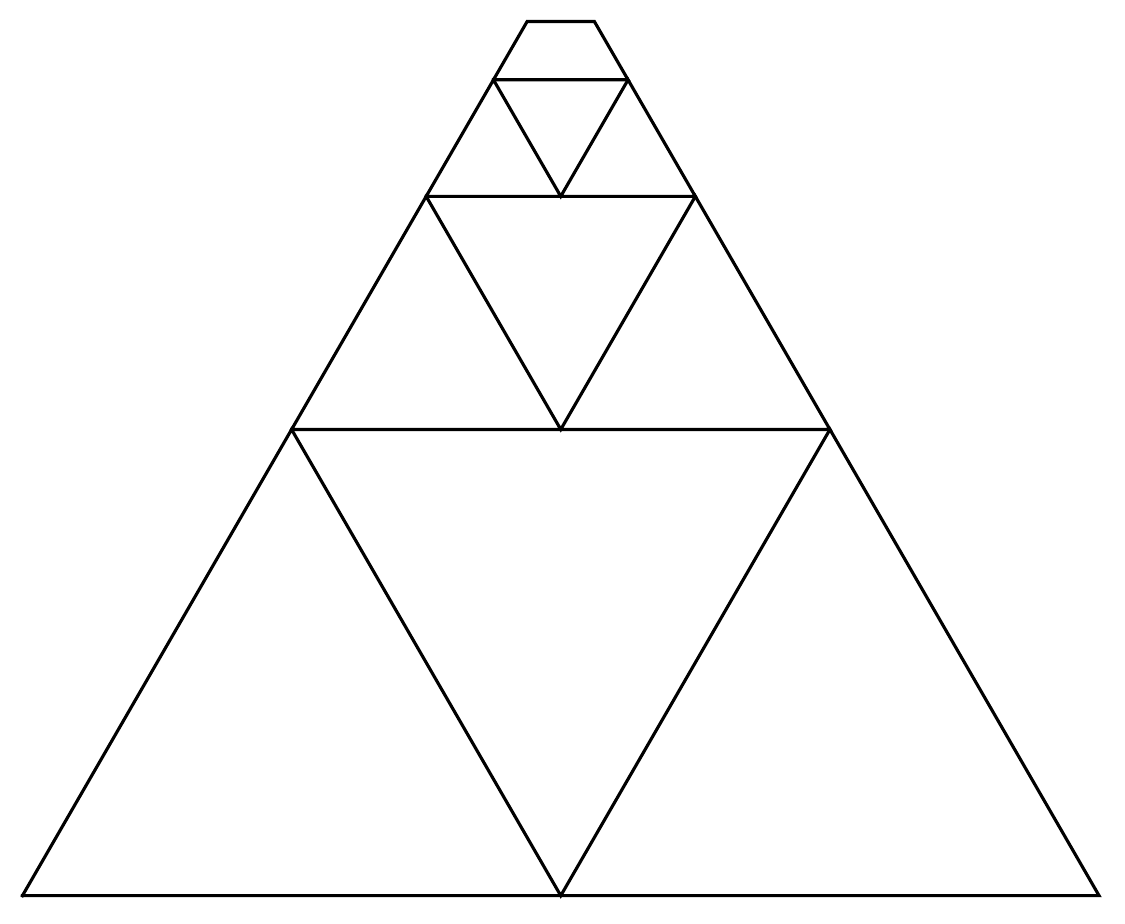}
\setlength{\unitlength}{1cm}
\begin{picture}(0,0) 
\put(-3.7,4.1){$a_0$}
\put(-4.0,3.55){$a_1$}
\put(-4.7,2.4){$a_2$}
\put(-5.0,1.5){$\cdot$}
\put(-5.05,1.4){$\cdot$}
\put(-5.10,1.3){$\cdot$}
\put(-5.9,-0.1){$a_m$}

\put(-2.7,4.1){$-a_0$}
\put(-2.35,3.55){$-a_1$}
\put(-1.7,2.4){$-a_2$}
\put(-1.1,1.5){$\cdot$}
\put(-1.05,1.4){$\cdot$}
\put(-1.0,1.3){$\cdot$}
\put(-0.4,-0.1){$-a_m$}

\put(-3.05,3.3){$0$}
\put(-3.05,2.15){$0$}
\put(-3.05,-0.15){$0$}
\put(-3.05,4.45){$0$}
\end{picture}
\vspace{0.1cm}
\begin{center}
\textbf{Figure 4.7. The values of $h$.}
\end{center}
\end{center}
\end{figure}

Using the matching conditions at $F_0^{m-k}q_1$ for $1\leq k<m$, we have
\[\frac{16}{3}a_{k}-\frac{5}{3}a_{k-1}-a_{k+1}=0.\]
So there exist two constants $c_1, c_2$ such that $$a_k=c_1 5^k+\frac{c_2}{3^k}.$$ Noticing there remains the equation at $F_0^mq_1$,
\[\big(\eta_1(S^m\lambda)+\eta_2(S^m\lambda)\big)a_0+2a_0-a_1=0,\] 
we have
\[\big(\eta_1(S^m\lambda)+\eta_2(S^m\lambda)+2\big)(c_1+c_2)-(5c_1+\frac{1}{3}c_2)=0.\]
This gives that 
\[\frac{c_2}{c_1}=-\frac{3\eta_1(S^m\lambda)+3\eta_2(S^m\lambda)-9}{3\eta_1(S^m\lambda)+3\eta_2(S^m\lambda)+5}.\]

On the other hand, we have
\[\begin{cases}
\partial_n^\leftarrow h(q_1)=\frac{5}{3}(2a_m-a_{m-1})=3\cdot 5^mc_1-5\cdot 3^{-m-1}c_2,\\
\partial_n^\leftarrow h(F_0^mq_1)=(\frac{5}{3})^m(a_1-2a_0)=3(\frac{5}{3})^m c_1-(\frac{5}{3})^{m+1}c_2.
\end{cases}\]

Thus the eigenvalue of $M^\lambda_{0^{m}}$ corresponding  to $\begin{pmatrix}1\\-1\end{pmatrix}$ is
\[\frac{\partial_n^\leftarrow h(F_0^mq_1)}{\partial_n^\leftarrow h(q_1)}=
\frac{14\cdot 5^m\big(\eta_1(S^m\lambda)+\eta_2(S^m\lambda)\big)}{(9\cdot 15^m+5)\big(\eta_1(S^m\lambda)+\eta_2(S^m\lambda)\big)+15(15^m-1)}.\hspace{2.23cm}\square\]

\textbf{Proposition 4.9.} \textit{Assume $e_1(\lambda)=e_2(\lambda)=\cdots=e_m(\lambda)=1$. Denote $0<x<1$  the solution of}
$
\frac{\eta_1(S^m\lambda)}{\eta_2(S^m\lambda)}=\frac{x+x^{-1}}{2}.
$
\textit{Then we have}

\textit{(a)}
 \[ \quad \eta_1(\lambda)=(\frac{5}{3})^m\frac{x-x^{-1}}{x+x^{-1}}\cdot\frac{x^{2^m}+x^{-2^m}}{x^{2^m}-x^{-2^m}}\eta_1(S^m\lambda)\textit{\text{ and }}  \eta_2(\lambda)=(\frac{5}{3})^m\frac{x-x^{-1}}{x^{2^m}-x^{-2^m}}\eta_2(S^m\lambda).\]

\textit{(b) For $w\in \tilde{W}_m^\lambda$, the matrix $M^\lambda_w$ is given by}
\[M^\lambda_w=\begin{pmatrix}
x^j & -x^{2^m-j}\\
-x^{j+1} & x^{2^m-j-1}
\end{pmatrix}\cdot
\begin{pmatrix}
1 & -x^{2^m}\\
-x^{2^m}  & 1
\end{pmatrix}^{-1},\]
\textit{with $j$ being the integer such that $0\leq j\leq 2^m-1$ satisfying 
\begin{equation}
j=\sum_{k=1}^{m}(w_k-1)\cdot 2^{m-k}.
\end{equation}}

\textit{Proof. } (a) By (4.16), we have
$\frac{\eta_1(S^{m-1}\lambda)}{\eta_2(S^{m-1}\lambda)}=2\big(\frac{\eta_1(S^m\lambda)}{\eta_2(S^m\lambda)}\big)^2-1=\frac{x^2+x^{-2}}{2}.$ 
Inductively, we have $$\frac{\eta_1(S^k\lambda)}{\eta_2(S^k\lambda)}=\frac{x^{2^{m-k}}+x^{-2^{m-k}}}{2}, \forall 0\leq k\leq m.$$

Still by (4.16), we have
\[\eta_2(\lambda)=\frac{5}{6}\cdot\frac{\eta_2(S\lambda)}{\eta_1(S\lambda)}\eta_2(S\lambda)=\cdots=(\frac{5}{6})^m\eta_2(S^m\lambda)\prod_{k=1}^{m}\frac{\eta_2(S^k\lambda)}{\eta_1(S^k\lambda)}.\]
Since
$$
\prod_{k=1}^{m}\frac{\eta_1(S^k\lambda)}{\eta_2(S^k\lambda)}=\prod_{k=1}^m\frac{x^{2^{m-k}}+x^{-2^{m-k}}}{2}=2^{-m}\frac{x^{2^m}-x^{-2^m}}{x-x^{-1}},
$$
we get
\[\eta_2(\lambda)=(\frac{5}{3})^m\frac{x-x^{-1}}{x^{2^m}-x^{-2^m}}\eta_2(S^m\lambda).\]
Taking $\frac{\eta_1(\lambda)}{\eta_2(\lambda)}=\frac{x^{2^m}+x^{-2^m}}{2}$ and $\frac{\eta_1(S^m\lambda)}{\eta_2(S^m\lambda)}=\frac{x+x^{-1}}{2}$ into the above result, we get the representation for $\eta_1(\lambda)$.

(b) Denote $p_j=F_wq_1$ with $j,w$ related by (4.19).   Write $p_{2^m}=q_2$. See Figure 4.8 for the location of $p_j$'s.

\begin{figure}[h]
\begin{center}
\includegraphics[width=10cm]{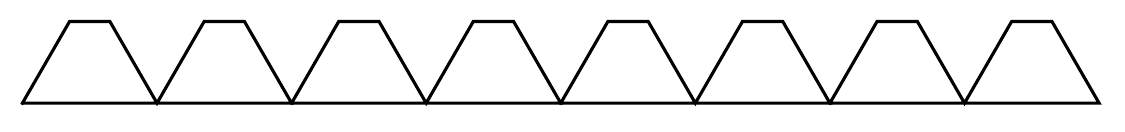}
\setlength{\unitlength}{1cm}
\begin{picture}(0,0) 
\put(-11.1,0){$q_1=p_0$}
\put(-8.95,0){$p_1$}
\put(-7.8,0){$p_2$}
\put(-6.55,0){$p_3$}
\put(-5.35,0){$p_4$}
\put(-4.2,0){$p_5$}
\put(-3.0,0){$p_6$}
\put(-1.8,0){$p_7$}
\put(-0.5,0){$p_8=q_2$}

\put(-9.8,1.1){$X_{111}$}
\put(-8.6,1.1){$X_{112}$}
\put(-7.4,1.1){$X_{121}$}
\put(-6.2,1.1){$X_{122}$}
\put(-5.0,1.1){$X_{211}$}
\put(-3.8,1.1){$X_{212}$}
\put(-2.6,1.1){$X_{221}$}
\put(-1.4,1.1){$X_{222}$}
\end{picture}
\vspace{0.1cm}
\begin{center}
\textbf{Figure 4.8.  The locations of $p_j$'s. ($m=3$)}
\end{center}
\end{center}
\end{figure}

Then for any harmonic function $h$ satisfying $h|_X=0$, we have
\[\eta_2(S^m\lambda) h(p_{j+2})-2\eta_1(S^m\lambda) h(p_{j+1})+\eta_2(S^m\lambda)h(p_{j})=0,\forall  0\leq j\leq 2^m-2.\]   By solving equations,  there exist two constants $c_1, c_2$ such that
\[h(p_j)=c_1x^j+c_2x^{2^m-j}.\]

Then we have
\[\begin{aligned}
\partial_n^\leftarrow h(p_j)&=(\frac{5}{3})^m\big(\eta_1(S^m\lambda) h(p_j)-\eta_2(S^m\lambda) h(p_{j+1})\big)\\
&=c_1'x^j-c_2'x^{2^m-j},\text{\quad for } 0\leq j\leq 2^m-1,
\end{aligned}\]
where $c_1'=c_1(\frac{5}{3})^m\big(\eta_1(S^m\lambda)-x\eta_2(S^m\lambda )\big)$ and $-c_2'=c_2(\frac{5}{3})^m\big(\eta_1(S^m\lambda)-x^{-1}\eta_2(S^m\lambda)\big).$ Similarly,  
\[
\partial_n^\rightarrow h(p_j)=-c_1'x^j+c_2'x^{2^m-j},\text{ for } 1\leq j\leq 2^m.
\]
Thus,
\[
\begin{pmatrix}
\partial_n^\leftarrow h(q_1)\\\partial_n^\rightarrow h(q_2)\end{pmatrix}=
\begin{pmatrix}1 & -x^{2^m}\\-x^{2^m} & 1 \end{pmatrix}
\begin{pmatrix}c_1'\\c_2'\end{pmatrix},
\begin{pmatrix}
\partial_n^\leftarrow h(F_wq_1)\\\partial_n^\rightarrow h(F_wq_2)\end{pmatrix}=
\begin{pmatrix} x^j & -x^{2^{m-j}}\\-x^{j+1} & x^{2^m-j-1} \end{pmatrix}
\begin{pmatrix}c_1'\\c_2'\end{pmatrix}.
\]
Comparing with (4.7), we have proved the part (b) of Proposition 4.9.        \hfill$\square$\\

Proposition 4.9 is motivated by Theorem 5.4 of [S1], which solves the special case that $\lambda=1-2^{-m}$.

\section{Extension to $\mathcal{SG}_l$}

In this section, we will briefly discuss how to extend previous results to general level-$l$ Sierpinski gasket $\mathcal{SG}_l$. 

\subsection{Dirichlet problem on half domains}  We still use $\Omega$ to denote the half domain and $X$ its Cantor set boundary. As shown in Section 2, to solve the Dirichlet problem on the half domain of $\mathcal{SG}_l$, it suffices to obtain the extension algorithm for $u|_{V_1\cap\Omega}$ in terms of the boundary data $f$. We summarize it into  following two steps.

\textbf{Step 1.} \textit{Find a formula for $\partial_n^\leftarrow u(q_1)$, in the form of 
\[\partial_n^\leftarrow u(q_1)=3f(q_1)-3\int_X fd\mu.\]}

Here the measure $\mu$, as shown in Theorem 2.1 for $\mathcal{SG}_3$, is a multiple of the normal derivative of the antisymmetric harmonic function $h_a$ on $X$. To work out $\mu$, we introduce some notations.

Let
$\{p_{j,\emptyset}\}_{j=1}^{[\frac{l}{2}]}=\{p\in V_1\cap X: p=F_iq_2 \text{ for some } 0\leq i\leq \frac{l^2+l-2}{2}\},$
and prescribe that $p_{j,\emptyset}$ locates above $p_{j+1,\emptyset}$ for each $1\leq j\leq [\frac{l}{2}]-1$. Write
\[\tilde{W}_m=\{i: \#(F_{i}\mathcal{SG}_l\cap X)=\infty\}^m,\text{ and } \tilde{W}_*=\bigcup_{m=0}^\infty\tilde{W}_m.\]
For $w\in \tilde{W}_*$, $1\leq j\leq[\frac{l}{2}]$, denote
$p_{j,w}=F_wp_{j,\emptyset}.$
Obviously, $\{p_{j,w}\}_{1\leq j\leq[\frac{l}{2}], w\in \tilde{W}_*}=V_*\cap X\setminus\{q_0\}$ is dense in $X$.
See the following example for an illustration.

\textbf{Example 5.1.}  In Figure 5.1 (a), we label the contraction mappings of $\mathcal{SG}_4$. So we have
$\tilde{W}_m=\{0,6\}^m.$ The vertices $\{p_{j,\emptyset}\}_{j=1}^2$ are plotted in Figure 5.1 (b).

\begin{figure}[h]
\begin{center}
\includegraphics[width=4.5cm]{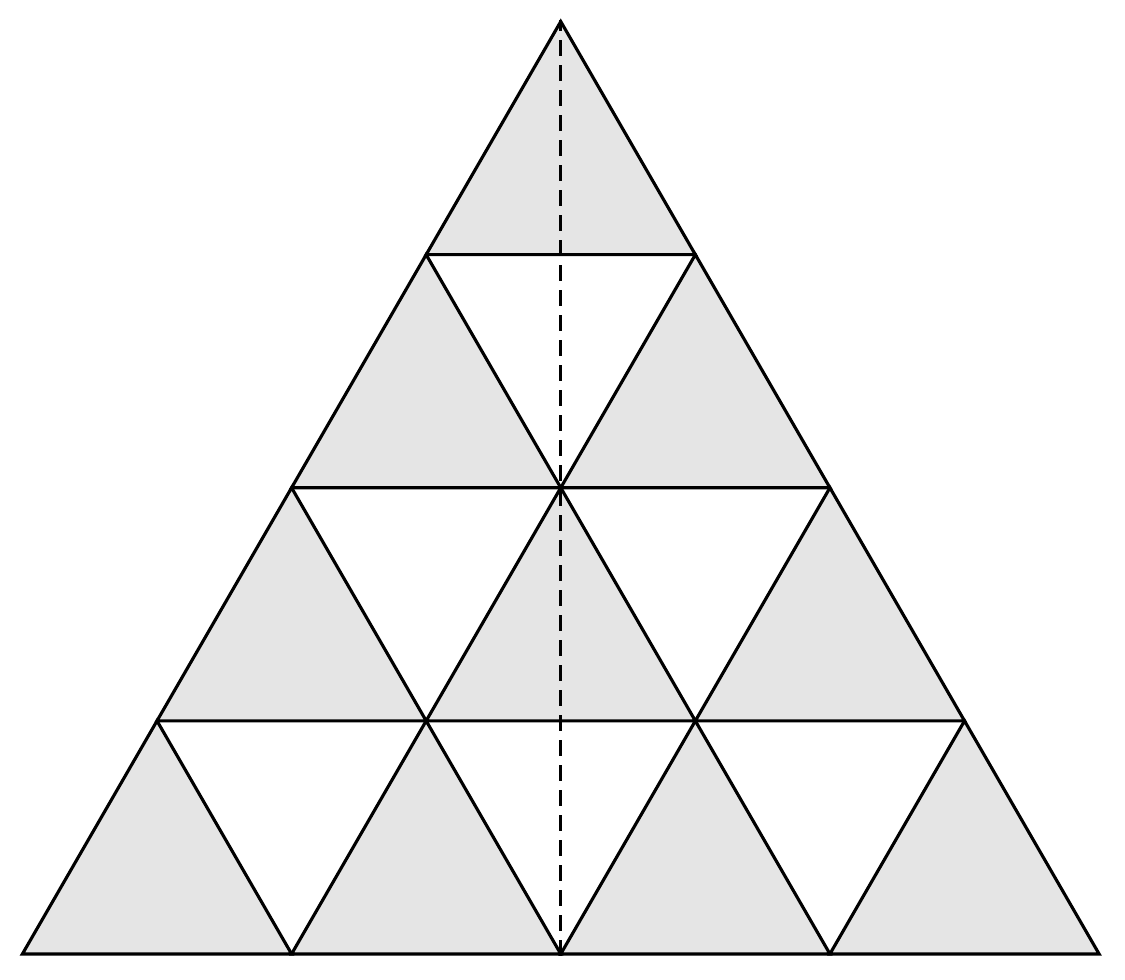}\hspace{1cm}
\includegraphics[width=2.3cm]{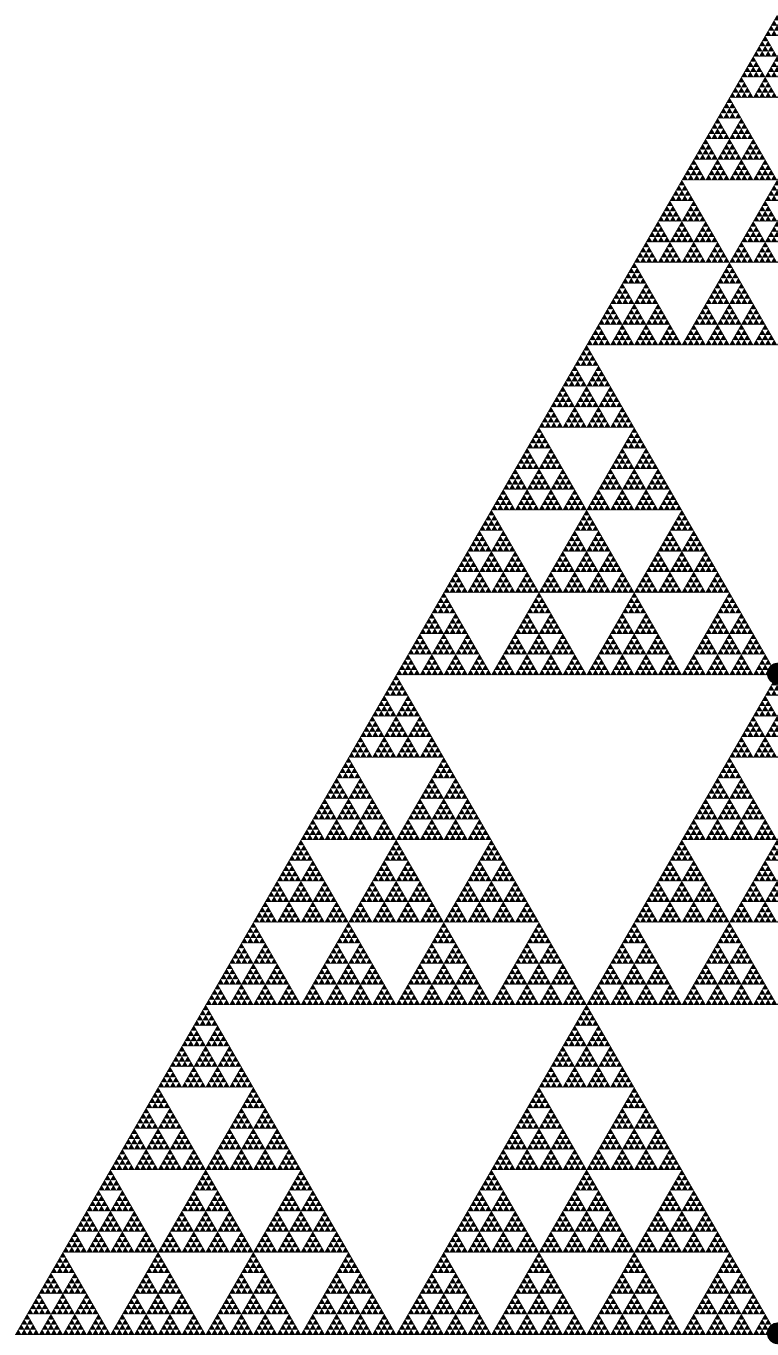}
\setlength{\unitlength}{1cm}
\begin{picture}(0,0) \thicklines
\put(-6.1,3.1){$F_0$}
\put(-6.65,2.2){$F_3$}
\put(-5.55,2.2){$F_4$}
\put(-7.15,1.25){$F_5$}
\put(-6.1,1.25){$F_6$}
\put(-5.0,1.25){$F_7$}
\put(-7.7,0.35){$F_1$}
\put(-6.65,0.35){$F_8$}
\put(-5.55,0.35){$F_9$}
\put(-4.45,0.35){$F_2$}
\put(-6.1,-0.4){(a)}

\put(-0.15,2.0){$p_{1,\emptyset}$}
\put(-0.15,0){$p_{2,\emptyset}$}
\put(-1.2,-0.4){(b)}
\end{picture}
\begin{center}
\vspace{0.4cm}
\textbf{Figure 5.1. The half domain of $\mathcal{SG}_4$.}
\end{center}
\end{center}
\end{figure}

With the above notations, introduce the pure atomic probability measure $\mu$ on $X$ satisfying 
\begin{equation}
\mu(\{p_{j,w}\})=-\frac{1}{3} \partial_n^\rightarrow h_a(p_{j,w}), \text{ for } w\in \tilde{W}_*\text{ and }1\leq j\leq [\frac{l}{2}].
\end{equation}

For $i\in \tilde{W}_1$, denote $\mu_i=r^{-1}h_a(F_iq_1)$ and write $\mu_w=\mu_{w_1}\mu_{w_2}\cdots\mu_{w_{|w|}}$, where $r$ is the renormalization factor of the energy on $\mathcal{SG}_l$. It is easy to check that $\partial_n^\rightarrow h_a(p_{j,w})=\mu_w \partial_n^\rightarrow h_a(p_{j,\emptyset})$, so that $$\mu(\{p_{j,w}\})=-\frac{\mu_w}{3} \partial_n^\rightarrow h_a(p_{j,\emptyset}).$$

\textbf{Step 2.} \textit{Solve linear equations determined by the matching conditions of normal derivatives on $V_1\bigcap \Omega$.}

That is
\begin{equation}
\begin{cases}
r\partial_n^\leftarrow u(F_iq_1)+\sum_{y\sim_1 F_iq_1,y\notin F_iV_0}\big(u(F_iq_1)-u(y)\big)=0,\quad i\in \tilde{W}_1,\\
\sum_{y\sim_1 x}\big(u(x)-u(y)\big)=0, \qquad\qquad\qquad\qquad \text{ for other } x\in V_1\bigcap \Omega.
\end{cases}
\end{equation}
Taking
$\partial_n^\leftarrow u(F_iq_1)=3r^{-1}\big(u(F_iq_1)-\int_X f\circ F_id\mu\big)$ into (5.2), the remaining problem is solving the linear equations. However,  even for the values of $h_a$, the calculation  becomes much more complicated, so we would not provide a general solution of (5.2) here. Nevertheless, let's look at the simplest case.

\textbf{Example 5.2.} Consider the half domain of $\mathcal{SG}$. In this case, $\tilde{W}_*=\{\emptyset,0,00,...\}$ and  $p_{1,\emptyset}=F_1q_2$. For convenience, we write $p_k=F_0^kp_{1,\emptyset}$. See Figure 5.2 below for the notations.

\begin{figure}[h]
\begin{center}
\includegraphics[width=2.2cm]{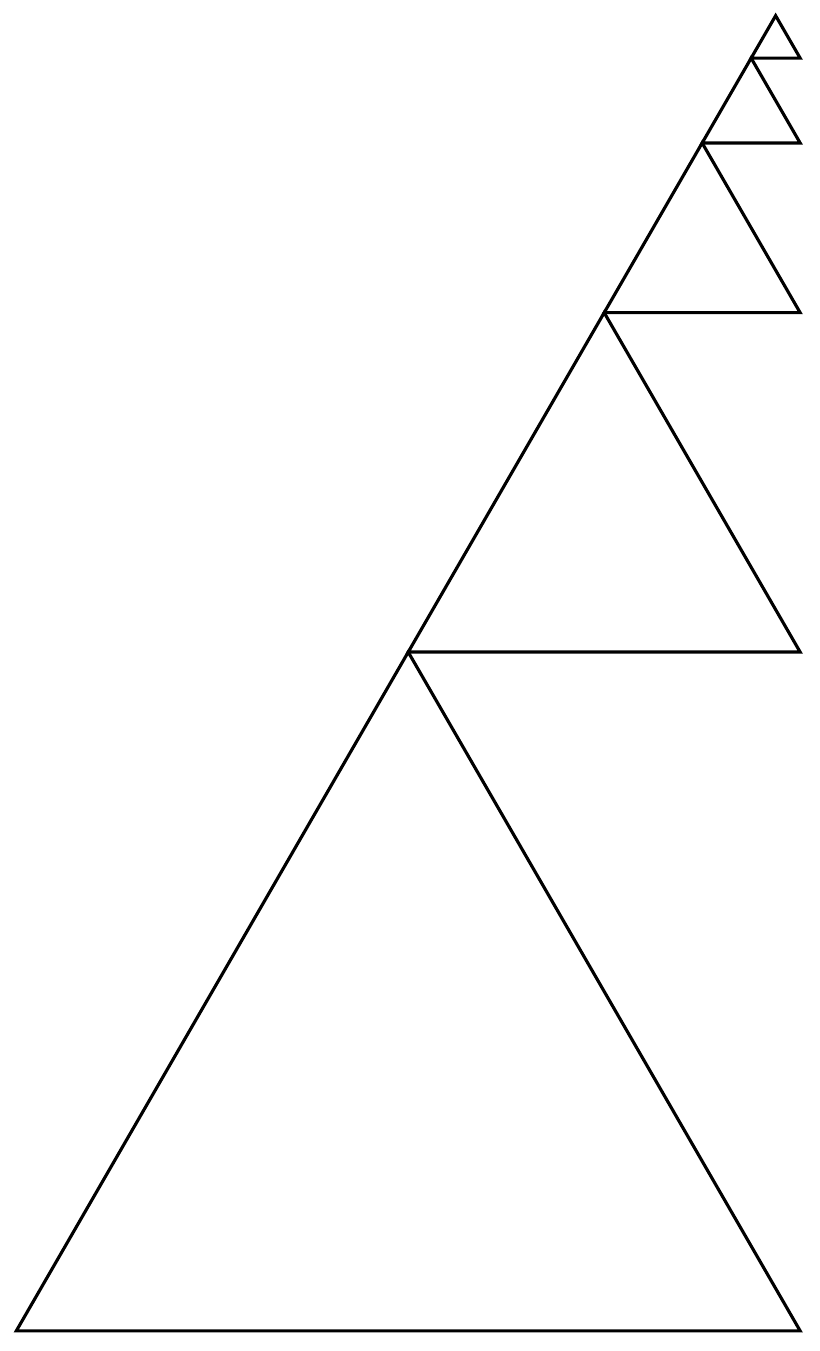}
\setlength{\unitlength}{1cm}
\begin{picture}(0,0) \thicklines
\put(-0.25,0){$p_0$}
\put(-0.25,1.8){$p_1$}
\put(-0.25,2.7){$p_2$}

\put(-2.7,0){$q_1$}
\put(-2.0,1.85){$F_0q_1$}
\put(-1.55,2.75){$F^2_0q_1$}

\put(-0.35,3.55){\LARGE$\cdot$}
\put(-0.35,3.5){\LARGE$\cdot$}
\put(-0.35,3.45){\LARGE$\cdot$}
\end{picture}
\begin{center}
\textbf{Figure 5.2. The vertices on a half $\mathcal{SG}$.}
\end{center}
\end{center}
\end{figure}

The measure $\mu$ is given by
\[\mu(\{p_k\})=2\cdot 3^{-k-1}.\]
There is only one equation in (5.2),
\[
2u(F_0q_1)-f(q_1)-f(p_0)+3\cdot \big(u(F_0q_1)-\int_X f\circ F_0 d\mu\big)=0.
\]
So
\begin{equation}
u(F_0q_1)=\frac{1}{5}f(q_1)+\frac{1}{5}f(p_0)+\frac{3}{5}\int_X f\circ F_0 d\mu
\end{equation}
which yields that
$$u(F_0q_1)=\frac{1}{5}f(q_1)+\frac{1}{5}f(p_0)+\frac{6}{5}\sum_{k=1}^{\infty} 3^{-k}f(p_k),$$
which is Corollary 2.5 in [LS].

We mention that the  Dirichlet to Neumann map was studied in [LS], and it was shown that $\|(\frac{3}{5})^{k+1}\partial_n ^\rightarrow u(p_k)\|_\infty<\infty$ if $\{f(p_k)\}_{k=0}^\infty\in l^\infty$. Here, we present another interesting observation, which describes where $\{\partial_n^\rightarrow u(p_k)\}_{k\geq 0}$ live in when
 $f\in C(\partial \Omega)$.

\textbf{Theorem 5.3.} \textit{Let $u$ be the solution of the Dirichlet problem (1.1) on the half domain of $\mathcal{SG}$, then the series  
$\sum_{k=0}^\infty (\frac{3}{5})^{k+1}\partial_n^\rightarrow u(p_k)$ converges. }

\textit{Conversely, given a convergence series $\{\eta_k\}_{k=-1}^\infty$, there exists a unique $f\in C(\partial \Omega)$ such that $f(q_1)=0$, $\partial_n^\leftarrow u(q_1)=\eta_{-1}$ and $\partial_n^\rightarrow u(p_k)=(\frac{5}{3})^{k+1}\eta_k$ for $k\geq 0$, where $u$ is the harmonic function with boundary data $f$.}

\textit{Proof.} By using (5.3) and the fact that $\int_X f\circ F_0^k d\mu=\frac{2}{3}f(p_k)+\frac{1}{3}\int_X f\circ F_0^{k+1} d\mu$, we have
\[\begin{aligned}
(\frac{3}{5})^{k+1}\partial_n^\rightarrow u(p_k)=&2f(p_k)-u(F_0^kq_1)-u(F_0^{k+1}q_1)\\
=&2f(p_k)-u(F_0^kq_1)-u(F_0^{k+1}q_1)+\frac{5}{2}\big(u(F_0^{k+1}q_1)-\frac{1}{5}u(F_0^kq_1)\\&-\frac{1}{5}u(p_k)-\frac{3}{5}\int_X f\circ F_0^{k+1}d\mu\big)\\
=&\frac{3}{2}\big(u(F_0^{k+1}q_1)-u(F_0^kq_1)\big)+\frac{3}{2}f(p_k)-\frac{3}{2}\int_X f\circ F_0^{k+1}d\mu\\
=&\frac{3}{2}\big(u(F_0^{k+1}q_1)-u(F_0^kq_1)\big)+\frac{9}{4}\big(\int_X f\circ F_0^kd\mu-\int_X f\circ F_0^{k+1}d\mu\big).
\end{aligned}\]
Thus we have $\sum_{k=0}^\infty (\frac{3}{5})^{k+1}\partial_n^\rightarrow u(p_k)$ converges to $\frac{9}{4}\int_X fd\mu-\frac{3}{4}f(q_0)-\frac{3}{2}f(q_1).$

Conversely, suppose there exists such a $f\in C(\partial\Omega)$,  it must satisfies
\[\begin{cases}
2f(p_k)-u(F_0^kq_1)-u(F_0^{k+1}q_1)=\eta_k,\\
f(p_k)+u(F_0^kq_1)-2u(F_0^{k+1}q_1)=\sum_{i=-1}^{k} (\frac{3}{5})^{k-i}\eta_i,
\end{cases} \forall k\geq 0,\]
which arise from the definition of normal derivatives at $p_k$ and $F_0^{k+1}q_1$, using the fact that $-\partial_n^\uparrow u(F_0^{k+1}q_1)=\partial_n^\leftarrow u(q_1)+\sum_{i=0}^{k}\partial_n^\rightarrow u(p_i)$. The solution of the above equations is
\[\begin{cases}
u(F_0^{k+1}q_1)=-\eta_{-1}-\frac{4}{3}\sum_{i=0}^{k}\eta_i+\sum_{i=-1}^{k}(\frac{3}{5})^{k-i}\eta_i,\\
f(p_k)=-\eta_{-1}-\frac{4}{3}\sum_{i=0}^{k}\eta_i+\frac{4}{3}\sum_{i=-1}^{k}(\frac{3}{5})^{k-i}\eta_i+\frac{1}{3}\eta_k,
\end{cases} \forall k\geq 0,\]
where $f(p_k)$ and $u(F_0^kq_1)$ converge to $f(q_0)=-\eta_{-1}-\frac{4}{3}\sum_{i=0}^{\infty}\eta_i$. Thus we find a unique function $f$ which satisfies the prescribed conditions. 

\hfill$\square$

However, it is not clear where $\{\partial_n^\rightarrow p_{j,w}\}_{w\in \tilde{W}_*,1\leq j\leq [\frac{l}{2}]}$ live in for general cases, even for $\mathcal{SG}_3$.

Next, we deal with the energy estimate of harmonic functions on $\bar{\Omega}$. For general $\mathcal{SG}_l$ case, we have 
\begin{equation}
C_1Q(f)\leq \mathcal{E}_{\Omega}(u)\leq C_2Q(f)
\end{equation}
for some positive constants $C_1,C_2$, with
\begin{equation}
Q(f)=\sum_{j=1}^{[\frac{l}{2}]}\big(f(q_1)-f(p_{j,\emptyset})\big)^2+\sum_{w\in \tilde{W}_*}\sum_{i\in \tilde{W}_1}\sum_{j,j'}r^{-|w|}\big(f(p_{j,w})-f(p_{j',wi})\big)^2.
\end{equation}
The method is essentially the same as that for $\mathcal{SG}_3$ case. We split the domain into countable pieces, namely $F_wO_1,w\in\tilde{W}_*$, with $O_1$ being $\bigcup\{F_i\mathcal{SG}_l|F_i\mathcal{SG}_l\subset \bar{\Omega},0\leq i\leq\frac{l^2+l-2}{2}\}$. Then $\mathcal{E}_\Omega (u)=\sum_{w\in \tilde{W}_*}\mathcal{E}_{F_wO_1} (u)$. To get the first inequality of (5.4), we rearrange the summation of $\mathcal{E}_\Omega(u)$ and estimate it directly, following the argument in the proof of Theorem 2.3. To get the second inequality, we just need to construct a piecewise harmonic function whose energy is less than a multiple of $Q(f)$. Readers please see the example below in the case of $\mathcal{SG}_4$.

\textbf{Example 5.4.} In Figure 5.3(a), we present the simple set $O_1$ of the half domain of $\mathcal{SG}_4$. Figure 5.3(b) gives the values of the constructed piecewise harmonic function $v$ on $O_2$, assuming 
\[
v|_{\partial \Omega}=f, \text{ and }v(F_wF_1q_0)=f(p_{1,w}), v(F_wF_1q_2)=f(p_{2,w}),\forall w\in \tilde{W}_*,
\]
and taking harmonic extension elsewhere. 

\begin{figure}[h]
\begin{center}
\includegraphics[width=3cm]{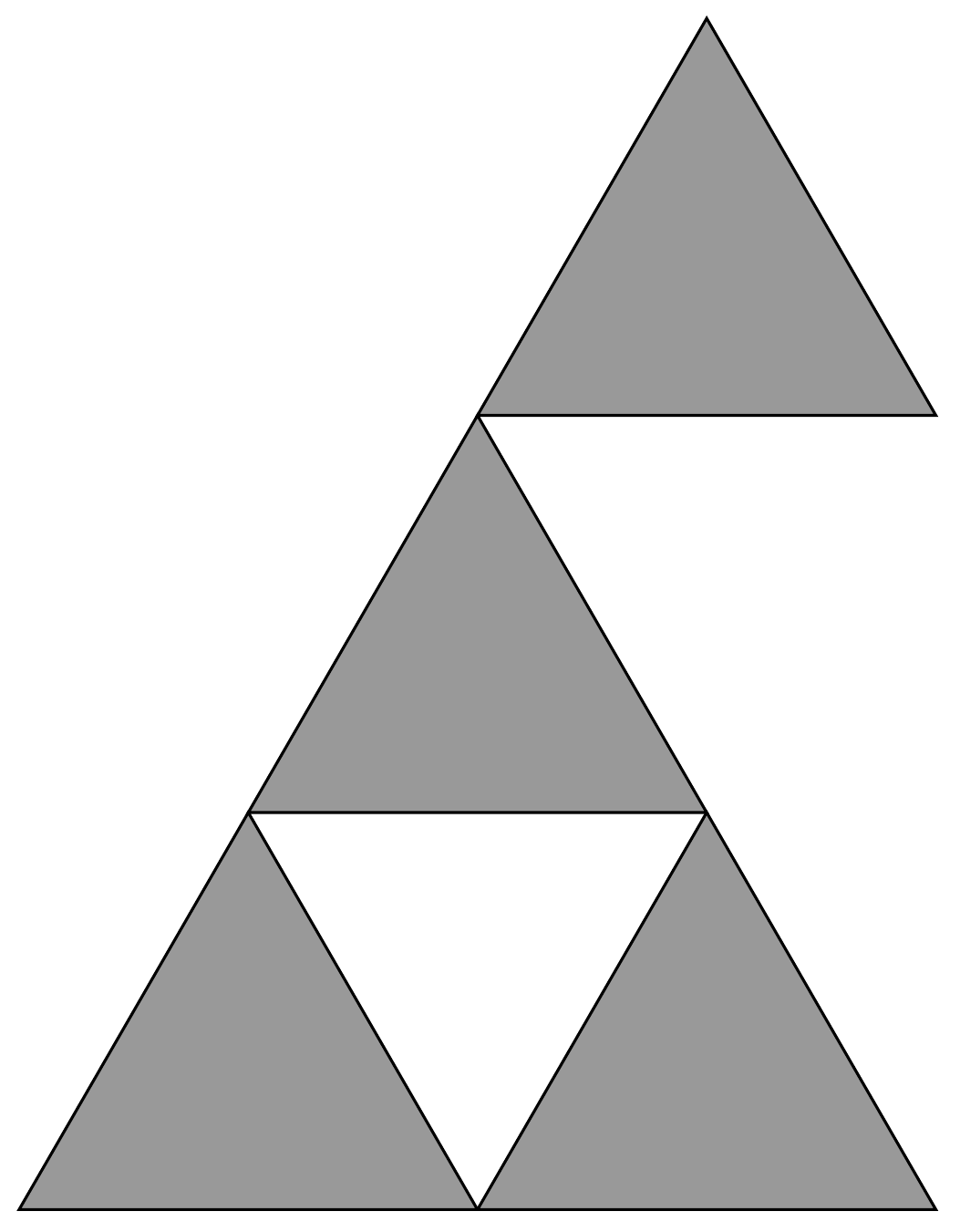}\hspace{1cm}
\includegraphics[width=3cm]{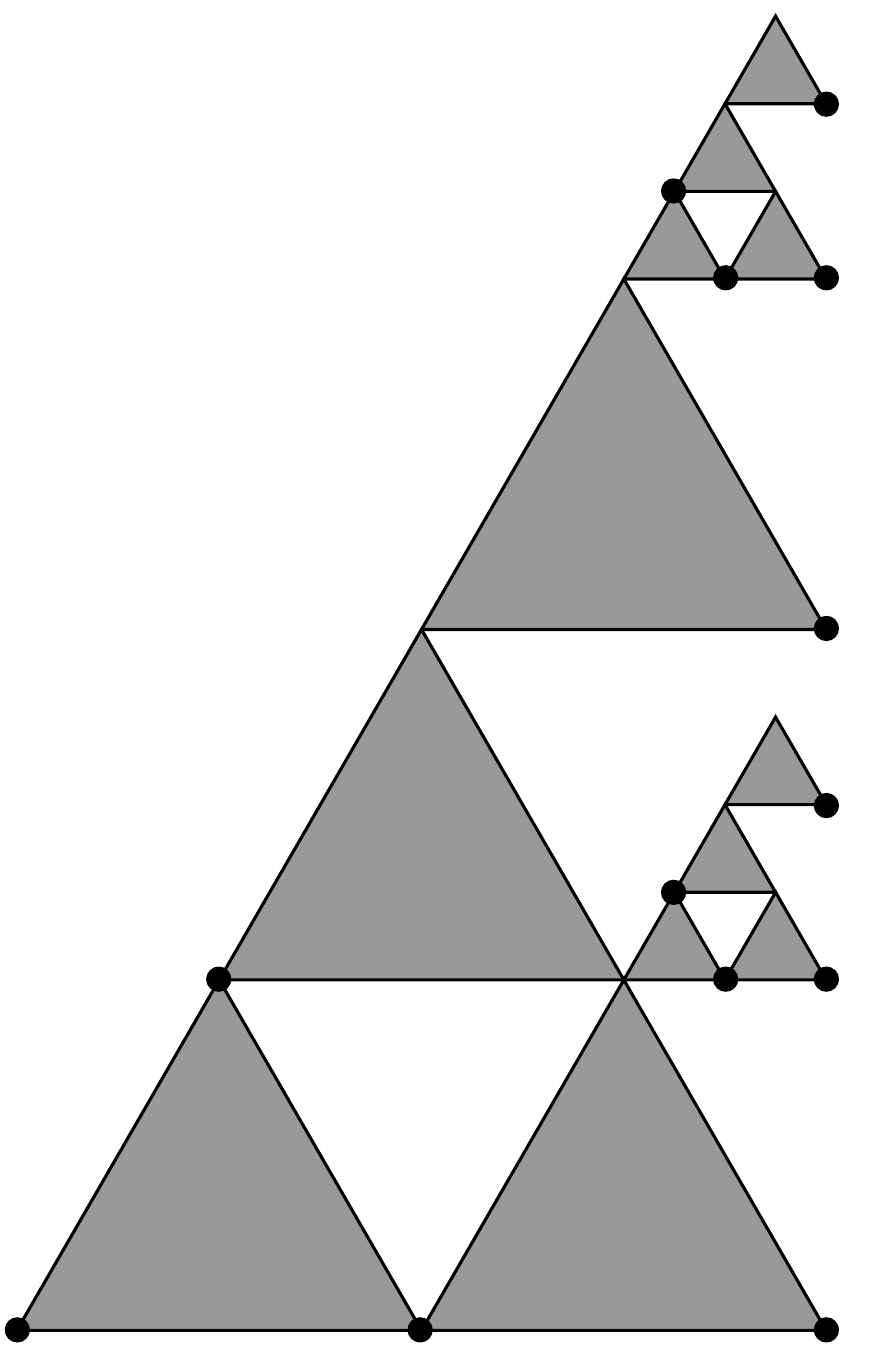}
\setlength{\unitlength}{1cm}
\begin{picture}(0,0) \thicklines
\put(-3.6,-0.1){\tiny$f(q_1)$}
\put(-0.3,0){\tiny$f(p_{2,\emptyset})$}
\put(-0.3,2.5){\tiny$f(p_{1,\emptyset})$}
\put(-2.2,-0.1){\tiny$f(p_{2,\emptyset})$}
\put(-3.4,1.3){\tiny$f(p_{1,\emptyset})$}

\put(-1.0,3.5){\tiny$f(p_{2,0})$}
\put(-0.3,3.7){\tiny$f(p_{2,0})$}
\put(-1.75,4){\tiny$f(p_{1,0})$}
\put(-0.3,4.2){\tiny$f(p_{1,0})$}

\put(-1.0,1.1){\tiny$f(p_{2,6})$}
\put(-0.3,1.3){\tiny$f(p_{2,6})$}
\put(-1.7,1.6){\tiny$f(p_{1,6})$}
\put(-0.3,1.9){\tiny$f(p_{1,6})$}

\put(-6.1,-0.5){(a)}
\put(-1.9,-0.5){(b)}
\end{picture}
\vspace{0.5cm}
\begin{center}
\textbf{Figure 5.3. $O_1,O_2$ and the values of $v$.}
\end{center}
\end{center}
\end{figure}

\subsection{Dirichlet problem on upper or lower domains} The Dirichlet problem on the upper and lower domains in general $\mathcal{SG}_l$ are much more complicated, as we need to discuss different cases of $\lambda$.

For the upper domains, we use the infinite expansion 
\[\lambda=\sum_{k=1}^\infty \iota_k\cdot l^{-m_k}\]
to characterize $\Omega_\lambda$, where $\{m_k\}_{k\geq 1}$ is an increasing sequence of positive integers and $\iota_k$ take values from $\{1,2,\cdots,l-1\}$. The number $\iota_n$ decides the relationship between $\Omega_{R^{n-1}\lambda}$ and $\Omega_{R^n\lambda}$, where $R\lambda=\sum_{k=2}^\infty \iota_k\cdot l^{-(m_k-m_1)}$, and there are $l-1$ choices  in the $\mathcal{SG}_l$ setting. 

As for the lower domains, we refer to a different expansion 
\[\lambda=\sum_{k=1}^{\infty}e_k(\lambda)l^{-k},\] 
with $e_k(\lambda)$ taking values from $\{0,1,\cdots,l-1\}$. 
We forbid infinitely consecutive $(l-1)$'s to make the expansion unique.  
Similarly, different $e_n(\lambda)$ determines different type of relationships between $\Omega_{S^{n-1}\lambda}^-$ and $\Omega_{S^n\lambda}^-$, where $S\lambda=\sum_{k=1}^{\infty}e_{k+1}(\lambda)l^{-k}$. 

The approaches in Section 3 and Section 4 to solve the Dirichlet problem on upper or lower domains still work, although the calculations involved turn to be  rather complicated. We list the main steps.
 
\textbf{Step 1.} Denote $\eta(\lambda)=\partial_n^\uparrow h_0(q_0)$ \big(or $\eta_1(\lambda)=\partial_n^\leftarrow h_1(q_1),\eta_2(\lambda)=-\partial_n^\rightarrow h_1(q_2)$\big). Represent $\eta(\lambda)$ in terms of $\eta(R\lambda)$ with the relationship between $\Omega_\lambda$ and $\Omega_{R\lambda}$ \big(or represent $\eta_1(\lambda)$, $\eta_2(\lambda)$ in terms of $\eta_1(S\lambda)$ and $\eta_2(S\lambda)$\big).  Use the above representations iteratively to approximate $\eta(\lambda)$ \big(or $\eta_1(\lambda),\eta_2(\lambda)$\big), and the proof is essentially the same as Lemma 3.1 (or Theorem 4.7).

\textbf{Step 2.} Calculate the normal derivatives of $h_0$ (or $h_1,h_2$) along the Cantor set $X$, using the crucial  coefficients $\eta(\lambda)$ \big(or $\eta_1(\lambda),\eta_2(\lambda)$\big). The normal derivatives of $h_0$(or $h_1,h_2$) hold the key to the representation of $\partial^\uparrow_n u(q_0)$ \big(or $\partial^\leftarrow_n u(q_1)$,$\partial^\rightarrow_n u(q_2)$\big) in terms of  the boundary data $f$.

\textbf{Step 3.} Solve the linear equations determined by the matching conditions of normal derivatives on the crucial points.

 Lastly, the Haar series expansion used in the energy estimate still works in general $\mathcal{SG}_l$ cases. The key observation is that we can still use the analogue of Theorem 3.4 to show that we can decompose the harmonic solution associated with a square integrable boundary value data (with respect to a suitable choice of measure), into a summation of countably infinite, pairwise orthogonal in energy, locally supported harmonic functions with suitable piecewise constant boundary values.

\end{document}